\documentclass[12pt,reqno]{article}
\usepackage[top = 2.5cm, bottom = 2.5cm, left = 2.5cm, right = 2.5cm]{geometry}
 \usepackage{url}
 \usepackage[usenames,dvipsnames]{color}
\usepackage[normalem]{ulem}

\title{Local algorithms, regular graphs of large girth, and random regular graphs}%
\author{
  Carlos Hoppen\thanks{Supported by FAPERGS~(Proc.\,11/1436-1) and CNPq (Proc.~486108/2012-0 and~304510/2012-2).} \\
 {\small Instituto de Matem\'{a}tica}\\ 
 {\small  Universidade
 Federal do Rio Grande do Sul}\\{\small Brazil}\\
 {\small {\tt  choppen@ufrgs.br }} 
\and
Nicholas  Wormald\thanks{Research supported by the  Canada Research
Chairs Program, NSERC and the Australian Laureate Fellowship program of the ARC} \\ {\small
School of Mathematical Sciences}
\\  {\small Monash University}\\ {\small  Australia}\\ {\small {\tt  nick.wormald@monash.edu }} 
}

  \newcommand{\lab}[1]{\label{#1}}                % hides labels

\newcommand{\remove}[1]{}

\usepackage{amssymb}
\usepackage{epsfig}
\usepackage{amsfonts}
\usepackage{amsmath}
\usepackage[mathscr]{euscript}
\usepackage{amsthm}
\usepackage{graphicx}

\newcounter{discounter}
 0

\newcommand{\bee}{\begin{equation}}
\newcommand{\ee}{\end{equation}}
\newcommand{\bea}{\begin{eqnarray}}
\newcommand{\eea}{\end{eqnarray}}

\newcommand{\bean}{\begin{eqnarray*}}
\newcommand{\eean}{\end{eqnarray*}}
\newcommand\eqn[1]{(\ref{#1})}
\newcommand{\bel}[1]{\bee\lab{#1}}
\newcommand{\halt}{*}

\newcommand{\Yv}{{\bf Y}}

\newcommand{\Op}{{\rm Op}}

\newcommand{\G}{{\cal G}}
\newcommand{\C}{{\cal C}}
\newcommand{\dv}{{\bf r}}
\newcommand{\zv}{{\bf z}}
\newcommand{\yv}{{\bf y}}

\newcommand{\nv}{{\bf n}}

\newcommand{\D}{{\cal D}}
\newcommand{\dk}{{d(k)}}

 \def\proof{\noindent{\bf Proof}\ \ }

\def\qed{~~\vrule height8pt width4pt depth0pt}

\def\A{{\cal A}}
\def\Gnr{{\cal G}_{n,r}}

\def\pr{{\bf P}}
\def\ex{{\bf E}}
\def\U{\mathscr{O}}
\def\output{\mathscr{E}}
\def\vu{\omega}

\def\eps{\epsilon}

\def\E{\mathbf{E}}

\def\d{\tau}
\def\r{R}

\def\y{\mathbf{y}}
\def\Y{\mathbf{Y}}

\def\vs{\mathbf{s}}

\def\se{\diamond}
\DeclareMathOperator{\Var}{Var}
\DeclareMathOperator{\Cov}{Cov}
\newcommand{\etal}{\textsl{et al.}}
\newcommand{\reals}{\ensuremath {\mathbb R} }
\newtheorem{thm}{Theorem}[section]
\newtheorem{cor}[thm]{Corollary}
\newtheorem{lemma}[thm]{Lemma}

\theoremstyle{definition}
\newtheorem*{defn}{Definition}
\theoremstyle{remark}

 \numberwithin{equation}{section}
 \numberwithin{figure}{section}
 \numberwithin{table}{section}
\makeatletter
\def\cleardoublepage{
  \clearpage
  \if@twoside\ifodd\c@page\else
  \hbox{}
  \thispagestyle{empty}
  \newpage
  \if@twocolumn\hbox{}\newpage\fi
  \fi\fi
  }

\makeatother

\begin{document}
\date{}

\maketitle

\begin{abstract}
We introduce a general class of algorithms and supply a number of general results useful for analysing these algorithms when applied to regular graphs of large girth. As a result, we can transfer a number of results proved for random regular graphs into (deterministic) results about all regular graphs with sufficiently large girth. This is an uncommon  direction of   transfer of results, which is usually from the deterministic setting to the random one. In particular, this approach enables, for the first time, the achievement of results equivalent to those obtained on random regular graphs by a powerful class of algorithms which contain prioritised actions.   As examples, we obtain new upper or lower bounds on the size of maximum independent sets, minimum dominating sets, maximum and minimum bisection, maximum $k$-independent sets, minimum $k$-dominating sets and minimum connected and weakly-connected dominating sets in $r$-regular graphs with large girth. 
\end{abstract}

\section{Introduction}
\lab{s:intro}

Graphs of large girth have been of interest since the pioneering result in probabilistic graph theory of Erd\H{o}s~\cite{erd59} in the late 1950's showing that, for any given positive integers $k$ and $g$, there exists a graph with girth at least $g$ and chromatic number at least $k$, evincing the global character of the chromatic number of a graph. 

Since then, the interaction between the local and global behaviour of graphs has come under much more scrutiny. We concentrate in this article on $r$-regular graphs, where $r$ is fixed, which provide interesting examples of sparse graphs. For such graphs, Borodin and Kostochka~\cite{BK} and Lawrence~\cite{ncw:lawrence} (see also Catlin~\cite{ncw:catlin}) obtained a kind of antithesis of  Erd\H{o}s' result, by showing that an $r$-regular  graph with girth at least 4 has chromatic number at most $\frac34 (r+2)$, thereby disproving a conjecture of Gr{\"u}nbaum. Improvements came later, in particular for large $r$, with Johansson~\cite{J} showing that the chromatic number of triangle-free graphs of maximum degree $r$ is at most $9r/\log r$. More recently, a number of results have appeared proving that some particular property is satisfied by all $r$-regular graphs having a sufficiently large girth (see Section~\ref{s:related} for a number of examples).
 
In this paper, we prove  properties of regular graphs with large girth by studying certain algorithms (some pre-existing, some new) and showing a relationship between their behaviour on random regular graphs and their behaviour on regular graphs of large girth. 

By suitably choosing the algorithm to produce an appropriate structure, we obtain upper or lower bounds on a variety of well studied graph parameters such as the size of the maximum independent set, the minimum dominating set, the minimum bisection or the minimum feedback vertex set, in all regular graphs of sufficiently large girth. Broadly speaking, we show that an algorithm $\A$ belonging to a quite large class of algorithms, the size of the structure produced by $\A$ is almost the same for $r$-regular graphs of very large girth, as it is for a random $r$-regular graph. This can be viewed as a partial converse of existing one-directional results which   translate properties of regular graphs of large girth into highly likely properties of random regular graphs. Throughout this paper,  the degree in a regular graph is a fixed number $r$, and `sufficiently large girth' means that there exists $g$ such that the statement holds when the girth of the graphs under consideration is at least $g$.

Lauer and second author~\cite{LW} introduced a new approach to proving properties of large girth regular graphs, by analysing a greedy algorithm for independent sets. The method used simple expectation, in contrast to more sophisticated sharp concentration techniques that some of the results on triangle-free graphs employed.  The particular greedy algorithm analysed there bore a strong resemblance to the greedy algorithm analysed in~\cite{des} for random regular graphs, and resulted in exactly the same asymptotic bound on the size of the independent sets.  This opened up the possibility that other algorithms which have previously been analysed on random regular graphs might have analogues for regular graphs of large girth which achieve the same bounds. In particular, there are many more sophisticated  algorithms, producing powerful bounds, which are based on iterating an operation involving a random vertex of minimum degree (and then deleting it and some neighbours). Such algorithms are called {\em prioritised} in~\cite{deprio}, where a general approach to analysing them was introduced that involved {\em deprioritising} the algorithms. This assists the analysis when the algorithm passes through a number of phases. In the present paper we analyse a general class of algorithms, {which we call {\em local deletion algorithms} because they repeatedly carry out explorations of neighbourhoods of vertices, and delete the edges that are exposed in the exploration. In particular, we provide an explicit system of differential equations whose solutions track the typical value of several natural variables as the algorithm is applied.  The system is the same as for  corresponding   algorithms applied to random regular graphs. Using this, we are able to transfer the results on many algorithms, including the deprioritised ones mentioned above, for random regular graphs into  (deterministic) results about all regular graphs with sufficiently large girth. 

We obtain many new bounds on functions of graphs with large girth in this work, but we see the main contribution as supplying a number of general results useful for analysing local deletion  algorithms, and in particular showing new relations between random regular graphs and regular graphs of large girth.  Many algorithms we analyse give results on global properties of random graphs, despite the local nature of the algorithms. Moreover, the relevant local structure of the graphs we consider is fixed --- the $r$-regular tree. Further investigation of this area may be interesting in light of the recent interest in graph limits in the sparse case. See, for instance, Lov{\'a}sz~\cite{Lbook}, Hatami, Lov\'{a}sz and Szegedy~\cite{HLS}, and the references therein. Greedy algorithms are considered in connection with the development of constant-time approximation algorithms to compute a series of parameters of graphs with bounded degree sequence, as in Elek and Lippner~\cite{EL}. In particular, Gamarnik and Sudan~\cite{GS} showed that, for sufficiently large $r$, local algorithms (which include our local deletion algorithms) cannot approximate the size of the largest independent set in an $r$-regular graph of large girth by an arbitrarily small multiplicative error. However, this appears to say nothing about small $r$, say $r=3$.

A preliminary version of some of the ideas in this paper first appeared in the doctoral thesis of the first author~\cite{hoppen_thesis}, of which the second author was the supervisor. In particular, this included an early version of our incorporation of the results achievable by prioritised algorithms acting on random regular graphs into the framework of~\cite{LW} using ideas from~\cite{deprio}.

\section{Related work and examples of improved bounds}\lab{s:related}

A number of graph parameters have been considered for $r$-regular graphs with large girth. Since these are standard parameters, we refrain from defining them until their definition is important. Kostochka~\cite{K} showed that a $d$-regular graph ($d\ge 5$ and fixed) with sufficiently large girth has chromatic number at most $2+\lfloor d/2\rfloor$. Hopkins and Staton~\cite{hopkins_staton} gave lower bounds on the size of independent sets in cubic (i.e.\ 3-regular) graphs with large girth, superseded by Shearer's bounds~\cite{shearer1,shearer2} for $r$-regular graphs in general. For cubic graphs on $n$ vertices (and sufficiently large girth), this bound was $0.4139n$. Kr{\'a}l,  \v{S}koda and   Volec~\cite{KSV} showed that a cubic graph of sufficiently large girth has a dominating set of size at most $0.299871n$. For definitions of the graph parameters  mentioned here, the reader may refer to Sections~\ref{sec:applications} and~\ref{s:more}, where we give a detailed account of our results mentioned in this section.

The above-mentioned paper~\cite{LW} improved Shearer's bounds on independent sets for all $r\ge 7$. The method of~\cite{LW} was also extended by the present authors in~\cite{hoppen_wormald1} to obtain lower bounds on the size of a maximum induced forest in large girth $r$-regular graphs. These are equivalent to upper bounds on vertex feedback sets, or decycling sets.

As mentioned earlier, the first author's thesis~\cite{hoppen_thesis} contains an embryonic version of this work. There, the above-mentioned bounds of Shearer on independent sets in large girth $r$-regular graphs were all improved. The theory in the present paper permits these results to be derived much more economically (see Section~\ref{Isets thesis}). (A significant part of the improvement  is that a general result has replaced the need to prove the `independence lemma' of~\cite{LW} for each algorithm considered.)   In the cubic case the new bound was $0.4328n$. After the appearance of~\cite{hoppen_thesis}, Kardo\v{s},    Kr{\'a}l and  Volec~\cite{VolecFrac} used an adaptation of its technique to analyse a better algorithm which improved the bound to 
$0.4352n$. (We note that their method in this paper seems much closer to the earlier~\cite{LW}, since the algorithm does not correspond to a prioritised one: the probabilities chosen as parameters of the algorithm are fixed.)  This result was yet again improved, to $0.4361n$, very recently by Cs{\'o}ka,  Gerencs{\'e}r, Harangi and Vir{\'a}g~\cite{CGHV}  using invariant Gaussian processes on the (infinite) $d$-regular tree.
 
Also evolving from the early version of our  method in~\cite{hoppen_thesis}, Kardo\v{s}, Kr{\'a}l and Volec~\cite{Volec} found strong lower bounds on the size of maximum cuts in cubic graphs with large girth,  which improved on earlier bounds of Z\'{y}ka~\cite{zyka},  and the first author~\cite{forests_extended_abs} found lower bounds on the largest induced forest in regular graphs with large girth.  

We apply our general results to many problems and improve, in various ways, all of the results on sufficiently large girth graphs mentioned above. As special cases, we obtain new lower bounds on the size of maximum independent sets, minimum dominating sets, maximum and minimum bisection, maximum $k$-independent sets and minimum $k$-dominating sets in $r$-regular graphs with large girth. In particular, for independent sets we improve the known results: both the results in~\cite{VolecFrac} and~\cite{CGHV} (which only considered the cubic case), and in~\cite{LW} and~\cite{shearer2} (for $r\ge 4$). For dominating sets we improve the result in~\cite{KSV} (which considered only the cubic case), and for higher degrees we give the first bounds specifically obtained for  large girth. These improve on a general upper bound due to Reed~\cite{reed} and on refinements thereof obtained by Kawarabayashi, Plummer and Saito~\cite{KPS} for graphs with a 2-factor. Additionally, our upper bounds are stronger in the sense that they also hold for minimum dominating sets which also have the additional property of being independent sets. Strangely perhaps, relaxing the independence condition does not produce any easy significant improvement. For maximum cuts in cubic graphs, we derive analytic lower bounds, based on differential equations, such that the bounds obtained in~\cite{Volec} can be interpreted as long recurrences for computing numerical approximations of our analytic result.

The method introduced in~\cite{LW} produces explicit bounds on the parameter involved for graphs of given girth. Using adaptations of this approach, Gamarnik and Goldberg~\cite{GG} obtained lower bounds for various paramaters in $r$-regular graphs with an explicit large girth.  Similarly, in~\cite{VolecFrac} and \cite{Volec} whose method evolved from~\cite{LW}, explicit values of girth are given. In the present paper we have introduced, amongst other things, the use of sharp concentration, which has resulted in the loss of a direct connection with explict girths. If desired it would still be possible to extract a version of the bounds that we obtain as a function of girth. 

We observe here that greedy algorithms employed to obtain results on the chromatic number of random regular graphs, as by Shi and Wormald~\cite{colouring},  are not in the scope of our methods at present.  It is a conjecture of the second author that 4-regular graphs with sufficiently large girth are all 3-colourable.

\section{Introduction to the general results}
\lab{s:introresults}

In this work, we often consider  some given property $P$ that a set of vertices of an input graph $G$ might have, and consider the function 
$$f_P(G)=\max \{|U| \colon U \subseteq V(G) \textrm{ and }U \textrm{ satisfies }P\}.$$
For instance, if $P$ is the property of being an independent set, then $f_P(G)$ is the size of a maximum independent set in $G$. Given a property $P$ and a positive integer $r$, our main objective is to find a constant $c(r,P)$ such that, if $G$ is an arbitrary $n$-vertex $r$-regular with sufficiently large girth, then $f_P(G) \geq c(r,P)\, n$. We also consider upper bounds similarly (though we could by complementation), and we may consider subsets of the edge set $E(G)$, or of $V(G)\cup E(G)$.

We say that a sequence of events $A_n$ defined in probability spaces $\Omega_n$ occurs asymptotically almost surely (a.a.s.) if $\lim_{n \rightarrow \infty} \pr(A_n)=1$. We sometimes mix this notation with other asymptotic notation such as $O()$ and $o()$ in the same statement. There is a natural interpretation of this, for formal definitions see~\cite{GTH}.
Given positive integers $n> r$ (where $nr$ is even for feasibility), consider the probability space $\mathcal{G}_{n,r}$ of all $r$-regular graphs with vertex set $V=\{1,\ldots,n\}$ with uniform probability distribution. It is a well known fact (see for instance Bollob\'{a}s~\cite{bollobas0} and Wormald~\cite{wormald1}) that, for fixed integers $r$ and $g$, the expected number of vertices in $\mathcal{G}_{n,r}$ that lie in cycles of length at most $g$ is bounded. An immediate consequence of this, which has often been used, is that if a random $r$-regular graph a.a.s.\ has no sets $U$ satisfying $P$ of cardinality at least $c_u(r,P)n$, then any bound $c(r,P)$ as described above must be less than $c_u(r,P)$. 

Another (coarse) consequence of the small number of short cycles in random regular graphs is the following.  Assume that $f_P(G)$ can only change by a bounded amount if a bounded number of vertices and edges are deleted from $G$  and that  $f_P(G) \geq c(r,P) \, n$ for every $n$-vertex $r$-regular graph $G$ with sufficiently large girth. Then, for every $\delta>0$, an $n$-vertex random $r$-regular graph $G'$ asymptotically almost surely satisfies $f_P(G') \geq (c(r,P)-\delta)n$. This will imply that all the results we obtain in the present paper for $r$-regular graphs of large girth also apply essentially unchanged to random $r$-regular graphs, although our method derives them via random $r$-regular graphs and makes such conclusions essentially redundant in this case. 

As a partial converse, if a random $r$-regular graph a.a.s.\ satisfies some property, then for every fixed $g>0$, a random $r$-regular graph with girth at least $g$ a.a.s.\ satisfies the same
property. This also follows from the results in~\cite{bollobas0} and~\cite{wormald1}.  However, not all asymptotic properties of random regular graphs hold for {\em all} $r$-regular graphs with girth sufficiently large. For instance, connectedness in the case $r\ge 3$.  We refer the reader to~\cite{models} for this and other basic results about random regular graphs. 

The interest in properties of random regular graphs has contributed to the development of a powerful method for analysing random processes, which establishes a connection between the process and an associated differential equation or system of differential equations, and is consequently known as the differential equation method. This was introduced by the second author in~\cite{des,desurvey} and given a general framework for random regular graphs in~\cite{deprio}, which gave bounds on graph functions obtained a.a.s.\ from a certain class of algorithms. We show in the present paper that the same bounds apply deterministically for regular graphs of sufficiently large girth, for  algorithms satisfying certain conditions. These conditions are not very restrictive, being satisfied in all the cases where the result in~\cite{deprio} has been used so far.

A good deal of research has been devoted to analysing the performance of algorithms on graphs with large girth.  An approach introduced by Lauer and Wormald~\cite{LW} is based on judiciously defining an algorithm that outputs a set $S$ with a property $P$, when applied to a fixed $r$-regular input graph $G$ of   girth at least $G$, such that the expected size of $S$ does not depend on the choice of $G$. This expectation is of course a lower bound on the maximum cardinality of a   set satisfying $P$, and also an upper bound on the  minimum set satisfying~$P$. 
  
In~\cite{LW}, and the ensuing papers based on its approach, by Hoppen and Wormald~\cite{hoppen_wormald1} and Kardos, Kr{\'a}l and Volec~\cite{Volec}, the expected output sizes were computed precisely. To do this, exact probabilities are calculated for some particular local scenarios through the solution of a system of recurrence equations. In a graph of large girth, the close neighbourhood of a vertex is a tree. Restricting the algorithms to having a local nature ensures that their action on the neighbourhood of a vertex is identical to their action on the neighbourhood of a vertex of an infinite $r$-regular tree. To obtain amenable calculations, `independence lemmas' are established, showing that in the branches of the tree induced by vertices near a given vertex, certain events are independent. This independence breaks down after a number of iterations of the algorithm, the number increasing with the girth of the graph. In some cases,  finding suitable events that were independent has been a significant difficulty.  

The main goal of  our present work is to adapt the idea in~\cite{LW}, which essentially uses a quite simple algorithm, to the much more powerful algorithms treated in~\cite{deprio}. The intrinsic difficulty lies in the nature of these algorithms: they use different operations in single steps, with prioritisation of the operations causing the `rules' to change potentially with each step. On the other hand, the analysis in~\cite{LW} requires uniform rules for large `chunks' of time. To deal with this, we essentially `deprioritise' the algorithms as in~\cite{deprio}. Additionally, we alter the  independence arguments by proving a relationship  with random regular graphs,  where a much stronger independence property holds (see Lemma~\ref{l:survival}).    This eliminates the need to find suitable independent events on a case by case basis, and provides much of the power of our approach, enabling easy direct evaluation of the relevant probabilities. However, since random regular graphs contain some short cycles, we now require some sharp concentration arguments, rather than the simple first moment method as in~\cite{LW,hoppen_wormald1}.
\medskip

\noindent
{\bf Two examples of algorithms: independent and dominating sets}
\smallskip

Before stating the main methodological results in this paper, we illustrate the type of algorithm under consideration. To this end, we describe two randomised procedures that find objects commonly studied in graph theory: a large independent set and a small dominating set. Given a graph $G=(V,E)$, an \emph{independent set} in $G$ is a set $S \subset V$ such that the induced subgraph $G[S]$  has no edges. A \emph{dominating set} in $G$ is a set $T \subset V$ such that any $v \in V$ lies in $T$ or has a neighbour in $T$. 

The randomised procedure $P_{ind}$ looks for a large independent set: starting with $G_0=G$, consider the \emph{survival graph} $G_{t}$ defined as follows for $t \geq 1$. Assuming that $i$ is the minimum degree of the vertices of $G_{t-1}$, the procedure chooses a vertex $v$ of degree $i$ uniformly at random among all such vertices. The vertex $v$ is added to the independent set, while the new survival graph $G_t$ is obtained by removing $v$ and all its neighbours from $G_{t-1}$. The procedure continues until the survival graph is empty. Clearly, this is a greedy procedure that yields a maximal independent set. 

In the case of dominating sets, we may consider a similar step-by-step randomised procedure $P_{dom}$: starting with $G_0=G$, the following rule is applied at every step $t \geq 1$. Let $i$ be the minimum degree of $G_{t-1}$, and choose a vertex $v$ of degree $i$ uniformly at random among all such vertices in $G_{t-1}$. If $i=0$ add it to the dominating set and remove it from $G_{t-1}$. Otherwise   choose a   vertex $w$ uniformly at random among all neighbours of $v$ with largest degree and add it to the dominating set. The new survival graph $G_t$ is obtained by removing $w$ and all its neighbours from $G_{t-1}$. 

\medskip

\noindent
{\bf The general case}
\smallskip

These two algorithms can be viewed as examples of a more general class of randomised algorithms, which we call \emph{local deletion algorithms}. Before giving formal definitions in the next section, we summarise the main features of these algorithms.
\begin{itemize}
\item[1.] Starting with an input graph $G_0=G$, a local deletion algorithm is an iterative algorithm to obtain a set of interest $\U$ with the property that, at each step $t \geq 1$, there is a \emph{selection step}, in which a set of vertices is selected at random in the survival graph $G_{t-1}$ according to some probability distribution. In both $P_{ind}$ and $P_{dom}$ a single vertex is selected, and the probability distribution is uniform over all vertices with minimum degree in the respective survival graph.

\item[2.] For each selected vertex $v$, there is an \emph{exploration step} to obtain information about vertices near $v$.  This step repeats iterations which check the degrees of neighbours of vertices already reached by the exploration step, and selects new vertices to be explored using some randomised rule which takes these degrees into account.

\item[3.]  There is an \emph{insertion step}, which adds some subset of the vertices queried in the exploration step, or some subset of the edges incident with them, to the set $\U$. In this step, the survival graph is then updated by deleting all vertices queried in the exploration step. In $P_{ind}$, the set of vertices added to the set consists of the single vertex $v$ selected, whilst in $P_{dom}$ it is the vertex $w$, unless $v$ is isolated, in which case $v$ is added to the set. More generally, we also define output  functions whose values are determined by the explored neighbourhood.  

\item[4.]  Both the exploration step and the insertion step  should depend only on the isomorphism type of the explored neighbourhood. 
\end{itemize}
We also extend the concept so as to apply to vertex-coloured graphs. These are not proper colourings: each vertex is merely assigned one of a finite list of colours. The {\em type} of a vertex is the ordered pair consisting of its degree and its colour. Then, where `degree' is mentioned above, we may read `type.' The insertion step becomes a \emph{recolouring step}, where each queried vertex is assigned one of a finite set of output colours, and the colours of their neighbours in the survival graph may be changed. A local deletion algorithm with no vertex colour is called \emph{native}. Our present description, as above, is restricted to native local deletion algorithms, but our results will be proved in the more general setting.

Although one possibility is to select just one vertex in each round,  the algorithm may also select a large number of vertices in a single step, which enables large output sets to be obtained in a bounded number of steps of the algorithm. The price to pay is that `clashes' may occur, where nearby vertices are both selected in the same round, so that their query graphs intersect, leading to conflicts in the insertion step. These conflicts must be resolved somehow.

A local deletion algorithm is said to be \emph{chunky} if, at each round $t \geq 1$, the set $S_t$ of vertices selected in the survival graph $G_{t-1}$ is obtained by selecting each of the remaining vertices independently with some probability $p_{t,i}$ that depends on its degree $i$. The main results in this work tell us that chunky local deletion algorithms can be analysed for regular graphs with large girth and that their performance is very similar to the performance of algorithms that have been previously considered in the random regular setting. Indeed, we prove (see Theorem~\ref{large girth}) that the expected size of the output set of a chunky local deletion algorithm is the same for every $r$-regular graph $G$ whose girth is sufficiently large in terms of the number of steps $N$ undertaken by the algorithm. In other words, the performance of the algorithm is independent of the particular input graph as long as the girth is large enough. Moreover, we demonstrate (see~Lemma~\ref{thm conc}) that the performance of the algorithm is almost the same when the input graph contains a small number of short cycles, which is a.a.s.\ the case for a random regular graph. As a consequence, computing the expected performance of a chunky local deletion algorithm for a fixed input graph with large girth is equivalent to computing its expected performance over random regular graphs (see Theorem~\ref{conc rreg}). 

In light of this equivalence, we may concentrate on the probability space of random regular graphs, which is a convenient setting for analysing chunky local deletion algorithms. As it turns out, we may show that the values of a series of variables associated with an application of a local deletion algorithm follow  the solutions of a system of differential equations  (see Theorem~\ref{chunkyranreg}), which may be obtained explicitly using the ideas in the proof of Lemma~\ref{l:trend}. This is summarised in Theorem~\ref{chunky des}, where we describe a general strategy for analysing a local deletion algorithm that is defined through specific rules for the selection and exploration steps. This analysis proceeds as if the algorithm acted on a random graph with restrictions on its degree sequence, but its conclusions lead to deterministic results about regular graphs with sufficiently large girth.

In Section~\ref{s:deprio}, we show that algorithms in a large class considered previously in the random regular setting can be `chunkified', that is, they can be turned into chunky local deletion algorithms with negligible influence on their performance (see Theorem~\ref{chunkydeprio}). This is the final ingredient for deriving deterministic bounds for graphs with large girth from previously known results (or, more correctly, their proofs) on random regular graphs. In particular, we give an extension to the large girth context of~\cite[Theorem 1]{deprio}, which determines the performance of a general class of algorithms in terms of the solutions of differential equations (see Theorem~\ref{t:fromDeprio}). The remainder of Section~\ref{sec:applications} gives some applications of this result, whereby we translate known results on random regular graphs into results on regular graphs of sufficiently large girth. Section~\ref{s:more} provides many more applications of our general results. Finally, we make some comments on extensions of our results in Section~\ref{s:final}.

\section{Local deletion algorithms}\lab{s:local}

The aim of this section is to formally define the class of local deletion algorithms. The basic descriptions are in Section~\ref{s:localCol}. In the subsequent section, we consider a particular family of local deletion algorithms which select a relatively large number of vertices in each step{. This motivates the term ``chunky." Since our object is the topic of graphs of large girth, in this section all graphs are simple, i.e.\ have no loops or multiple edges.

\subsection{Definitions for the general case}\lab{s:localCol}

Let $D$ be a positive integer. At this point, and in the remainder of this paper, the distance between two vertices $u$ and $v$ in a graph $G$ is denoted by $d_G(u,v)$, while the distance between $u$ and a subgraph $H$ is given by $d_G(u,H)=\min\{d_G(u,v) \colon v \in H\}$. (The subscript $G$ may be omitted if the graph is clear from context.) 

In order to give a precise description of a local deletion algorithm, it is necessary to properly define each of the steps mentioned in the informal description given in the introduction. This includes defining how to choose a set of vertices in each round (selection rule), how to explore the neighbourhood of a chosen vertex (local rule), and how to update the graph based on this exploration (recolouring rule). 

\begin{defn}[Transient, output and neutral colours, coloured graph, type of a vertex]
We shall assume that we have two sets of colours $\C$ and $\output$, which denote respectively the set of \emph{transient colours}, which are assigned to the vertices in the survival graph, and the set of \emph{output colours}, which are assigned to the vertices that are deleted from the survival graph.   We will apply our results to $r$-regular graphs in which all vertices initially have the same transient colour, which we call {\em neutral}. Moreover, we assume that the sets $\C$ and $\output$ have special colours, denoted by $\propto'$ and $\propto$, respectively, 
which will be useful when dealing with clashes.
A \emph{coloured graph} refers to a graph whose vertices are assigned colours in the set $\C$ of transient colours. Given a coloured graph $G$ and a vertex $v \in V(G)$, the \emph{type} $\d_G(v)$ of $v$ is the ordered pair $(c,d)$, where $c$ is the colour and $d$ is the degree of $v$ in $G$.
\end{defn}

\begin{defn}[Selection rule]
A selection rule is a function $\Pi$ that, for a nonempty coloured graph $G=(V,E)$, gives a probability distribution $\pi_G$ on the power set of $V$ with the properties that
\begin{description}
\item{(i)} $\sum_{v \in S \subseteq V} \pi_G(S)=\sum_{w \in T \subseteq V} \pi_G(T)$ for every $v,w \in V$ such that $\d_G(v)=\d_G(w)$;
	\item{(ii)}  $\pi_G(S)=0$ if $S$ contains any vertex  whose colour is $\propto'$.
\end{description}
\end{defn}

Thus, the probability that a vertex $v$ lies in the set of selected vertices is determined by its type $\d_G(v)$. For instance, a single vertex can be selected, with probability determined by its type. Alternatively, the selection may be done by assigning probabilities $p_j$ for choosing a vertex of type $j$ and then adding each vertex to the set independently with probabability determined by its type. 

Next, in order to define how the algorithm performs its exploration step, we introduce the concept of query graph. For this, we need to introduce a wildcard for an undetermined type, which we denote $\se$.  We regard $\se$ formally as a different type,  which has no colour or degree. We call the other types {\em vertex types} if we need to distinguish them from $\se$.
\begin{defn}[Query graph]
An $h$-vertex query graph is a nonempty coloured graph $H$ with vertex set $\{1,\ldots,h\}$ such that each vertex $v \in V(H)$ is associated with  a finite (possibly empty) multiset $\ell_v$, each of whose elements is a type (either vertex type or $\se$). A query graph has depth $D$ if every vertex is at distance at most $D-1$ from the vertex with label 1, which is called the root of the query graph. 
\end{defn}
The set of all query graphs is denoted by $\mathcal{Q}$, while the set of query graphs of depth $D$ is denoted by $\mathcal{Q}_{D}$. To avoid ambiguity arising from automorphisms, a copy of a query graph $H$ in a coloured graph $G$ is formally defined as follows. 
\begin{defn}[Copy of a query graph]
Given a coloured graph $G$ and a query graph $H$, a copy of $H$ in $G$ rooted at $v \in V(G)$ is a function $\psi:V(H)\rightarrow V(G)$ with the following properties:
\begin{itemize}
\item[(i)] $\psi(1)=v$;

\item[(ii)] $\psi$ is a (graph theoretical) colour-preserving isomorphism from $H$ to a subgraph   of $G$; 

\item[(iii)] for every $w \in V\big( \psi(H)\big)$ the multiset $\ell_w$ given by the types of neighbours of $w$ in $G-E\big( \psi(H)\big)$ contains the multiset of vertex types in $\ell_{\psi^{-1}(w)}$ (respecting multiplicities). 

\end{itemize}
\end{defn}

In the description of the algorithm, a copy $\psi$ of a query graph will be used to record the information obtained after the selection of a vertex $v$ and some subsequent exploration steps in an underlying graph $G$. The names $1, 2, \ldots$ of the vertices in the query graph record the order in which the vertices are explored in the copy, so that the root of the query graph corresponds to the vertex $v$ at the start of the exploration step. The multiset $\ell_j$ associated with a vertex $j=\psi^{-1}(w)$ indicates the information that the local rule may use on the types of the neighbours of $w$ in $G$ that are not adjacent to $w$ in the  copy of the  query graph. The fact that the depth is limited by $D-1$ implies that the deletion of all vertices of the image of a query graph whose root is mapped to $v$ cannot affect the degree of any vertex whose distance from $v$ is larger than $D$. 

In Figure~\ref{f:query}, we show a query graph $H$ of depth $2$ and an embedding of this query graph into a larger graph. For simplicity, we assume that there is a single colour available, so that types are defined by their degrees. In the embedding of $H$ into the larger graph, the edges of $H$ are represented by solid lines, while the edges associated with the open adjacencies of $H$, or with the degree of such an open adjacency, are represented by dashed lines.

 \begin{figure}[hbt]
\centerline{\includegraphics[width=9cm]{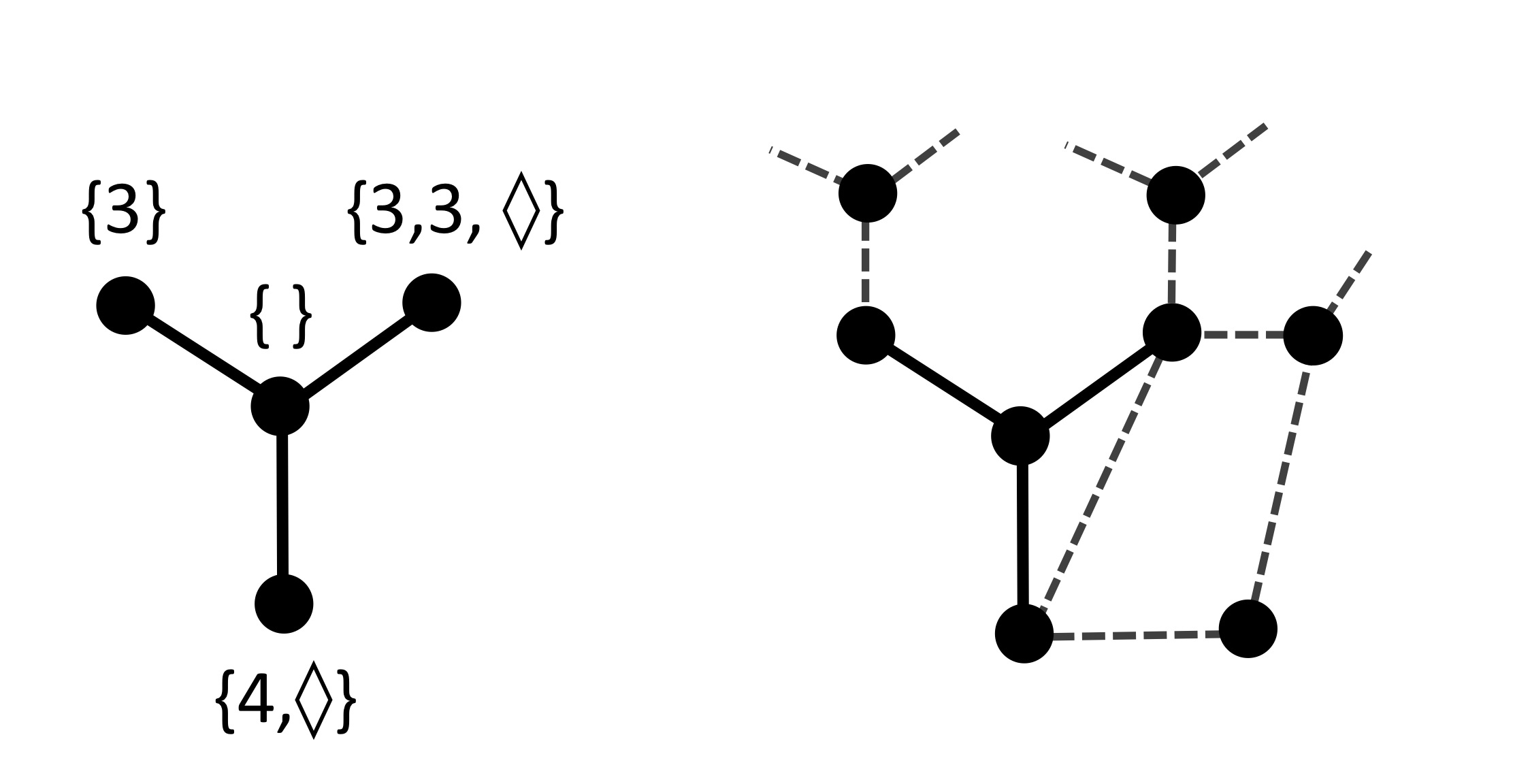}}
\caption{A query graph $H$ and its embedding in a larger graph}\label{f:query}
\end{figure}

The exploration steps mentioned above perform repetitions of a basic query operation, defined as follows.

\begin{defn}[Open adjacencies, querying, and exposing an edge] 
Let $G$ be a coloured graph, $H$   a query graph and $\psi$  a copy of $H$ in $G$ associated with the subgraph $\ \psi(H)$.  Define the set $U_i=\{(i,\d): i \in V(H) \textrm{ and }\d \in \ell_i\}$.  An element  $(i,\d) \in U_i$  is called an  {\em open adjacency}  of $i$ with type $\d$. Assume that each element of $\ell_i$  has  already been associated with a neighbour of $\psi(i)$. 
{\em Querying}  an open adjacency $(i,\d)$ consists of selecting u.a.r.\ a neighbour, $u$, of $\psi(i)$ that is associated with a copy of $\d$.  When such an open adjacency is queried, we say that the edge $\psi(i)u$ of $G$ has been \emph{exposed}.
\end{defn} 
  
 Next we define the crucial randomised rule determining how query graphs are extended.  This uses a symbol $\halt$ to denote that the exploration is terminated.
\begin{defn}[Local subrule, terminal open adjacencies, depth]
A local subrule is a function $\phi$ that, for each $H \in \mathcal{Q}$, specifies a probability distribution $\phi_H$ on the set  $U \cup \{\halt\}$, where  $U=\{(i,\d): i \in V(H) \textrm{ and }\d \in \ell_i\}$ is the set of open adjacencies. We require that $\phi_H(f)=0$ if $f$ is an open adjacency whose type has colour $\propto'$.  An open adjacency with the latter property is called {\em terminal} and all others are {\em non-terminal}. Moreover, we require that there is a natural number $D$ such that, given any query graph $H$, $\phi_H(f)=0$ for every open adjacency $f=(i,\d)$  for which the distance from $i$ to the root is at least $D-1$. The minimal such $D$ is the {\em depth} of  $\phi$.  Finally, we require that $\phi_{H}(\halt)=0$ if  $H$ is a single vertex $u$ such that   $\ell_u$ contains at least one $\se$, but no other type. 
\end{defn}
The last condition above ensures that the local subrule explores at least one open adjacency unless the initial vertex chosen has degree 0.

 We say that a local subrule is \emph{normal} if it always queries an open adjacency of the form $(i,\se)$, if there is one in the query graph. Most algorithms described in this paper use normal subrules.  Given a query graph $H$, the randomised choice of an open adjacency of $H$ (or of interrupting the exploration step) according to the probability distribution $\phi_H$ will be referred to as \emph{an application of the local subrule $\phi$ to $H$}. 

To conclude the description of the exploration step, we define how the algorithm uses a local subrule,  together with  querying, to build a copy of a query graph rooted at $v$.
\begin{defn}[Local rule]
The local rule $L=L_\phi$ associated with a local subrule $\phi$ is a function that maps an ordered pair $(v,G)$, where $G$ is a coloured graph and $v \in V(G)$, to a copy $\psi$ of a query graph $H$ in $G$. The image of $(v,G)$ under $L_{\phi}$ is defined inductively as follows:  start with the query graph $H_0$ that has a single vertex labelled $1$ associated with the multiset $\ell_1$ containing one copy of $\se$ associated with each neighbour of $v$ in $G$. Let $\psi_0$ be the bijection that maps $1$ to $v$. On subsequent steps $k>0$,  the local rule applies the local subrule for $H_{k-1}$ and then: 
\begin{itemize}
\item[(i)] If the outcome is $\halt$, the output $\psi$ is defined to be $\psi_{k-1}$. 

\item[(ii)]  If the outcome is $(i,\se)$,   the local rule queries it,
 obtaining a vertex $u$ adjacent to $w=\psi(i)$. Let $\d$ be the type of $u$. A new query graph $H_k$ is defined by replacing the occurrence of $\se$ in $\ell_i$ associated with $u$ by $\d$, while the copy $\psi_k$ of $H_k$ in $G$ is equal to $\psi_{k-1}$.  

\item[(iii)]  If the outcome is an open adjacency $(i,\d)$ with $\d\neq \se$, then the local rule queries it, obtaining a neighbour $u$ of $w$ having type $\d$. By induction, $\{w,u\}$ is not an edge of $\psi_{k-1}(H_{k-1})$. If $u \notin V(\psi(H_{k-1}))$, define the query graph $H_k$  from $H_{k-1}$ by adding a new vertex labelled $j=|V(H_{k-1})|+1$. A new multiset $\ell_j$ is created, with one $\se$ associated with each neighbour of $u$ in $G$ other than $w$. Moreover, the multiset $\ell_{i}$ is updated by removing the occurrence of $\d$ associated with $u$. The copy $\psi_k$ of $H_k$ in $G$ is the extension of $\psi_{k-1}$ obtained by mapping $|V(H_{k-1})|+1$ to $u$.  On the other hand, if $u \in V(\psi(H_{k-1}))$, say $u=\psi_{k-1}(j)$, the query graph $H_k$ is obtained from $H_{k-1}$ by adding the edge $ij$ and by removing the items corresponding to $u$ and $w$ from $\ell_{i}$ and   $\ell_{j}$ respectively. The copy $\psi_k$ of $H_k$ in $G$ is equal to $\psi_{k-1}$.  
\end{itemize}
\end{defn}

The copy $\psi$ of the query graph $H$ obtained when this process stops is called the \emph{query graph obtained by the local rule}. Furthermore, if the local subrule associated with the local rule has depth $D$, it is clear that the query graph produced cannot have a vertex at distance larger than $D-1$ from its root, so that the local rule is said to have {\em depth $D$} in this case. Note also that querying an open adjacency $(i,\se)$ does not alter the query graph, whilst querying  an open adjacency $(i,\d)$, where $\d$ is a vertex type, always adds an edge, and possibly a vertex, to the query graph.  

\begin{defn}[Recolouring rule]
A recolouring rule is a (possibly randomised) function $c$ that, given a query graph $H \in \mathcal{Q}$, assigns an output colour other than $\propto$ to each vertex $i$ of $H$ such that $\se \notin \ell_i$ and a transient colour other than $\propto'$ to each of the other vertices in $H$ and to each element $\d \ne \se$ in a multiset $\ell_i$ for $i\in V(H)$. In the latter, a type whose colour is already $\propto'$ is always assigned $\propto'$.  Optionally, $c$ also assigns an output colour to each edge of $H$.  
\end{defn}

In the algorithm, the recolouring rule defines how the algorithm uses the information given in the exploration step to update the survival graph, the output vector and the transient colours on vertices.
For many algorithms there is no need to colour edges, but it is convenient for some, so we include this as an option. Our algorithms do not assign transient colours to edges.
\begin{defn}[Local deletion algorithm]

A  local deletion algorithm consists of a triple $({\bf \Pi},L,c)$  with ${\bf \Pi}=\{\Pi_t\}_{t=1}^N$ for some $N$, where each $\Pi_t$ is a selection rule, and where $L$ is a local rule and $c$ is a recolouring rule. The {\em depth} of the algorithm is the depth of the local subrule $\phi$ associated with $L$.
Applied to an input graph $G=(V,E)$, the algorithm runs for $N$ steps and then outputs a vector $\vu$ according to the following prescription. Assume that the transient colours are in a set $\C$ and the output colours in a set $\output$ disjoint from $\C$. The algorithm starts with $G_0=G$, all of whose vertices initially have neutral transient colour, and repeats the following step $t$, $t=1,\ldots, N$.  Step $t$ is as follows:
\begin{itemize}
\item[(i)] (selection step) obtain a set $S_t \subset V(G_{t-1})$ using the selection rule $\Pi_t$;

\item[(ii)] (exploration step) for each $v \in S_t$, obtain  a copy $\psi_v$ of a query graph $H_v \in \mathcal{Q}$ through the application of the local rule $L=L_{\phi}$ to $G_{t-1}$ and $v$;

\item[(iii)]  (clash step) a vertex $u$ is called a {\em clash} if at least one of the following applies: (a) $u$ lies in $\psi_v$ for at least two vertices $v \in S_t$; (b) $u$ lies in a single $\psi_v$ but is adjacent to a vertex in some $\psi_w$, where $u \neq w$; (c) $u$ lies in a single $\psi_v$ but is adjacent to a vertex in $\psi_v$ through an edge in $G_t-E(\psi_v(H_v))$. Let $B$ be the set of all clashes. 

\item[(iv)] (recolouring step) the \emph{survival graph} $G_t$ is obtained from $G_{t-1}$  as follows. At first, the algorithm queries each open adjacency whose type is not $\se$, and the vertices found are placed in a multiset $W$.  All vertices and edges which are specified colours by $c$ are recoloured accordingly, with the exception that clash vertices and all incident edges are assigned $\propto$. Vertices in $G-\bigcup_{v \in S_t} \psi_v(H_v)$ which were placed in $W$  must appear in at least one list $\ell_i$ in a query graph.  Recolour each one with the colour assigned by $c$ to the corresponding element $\d$ in $\ell_i$, unless  they appear two or more times in $W$, in which case they are recoloured $\propto'$. All vertices and edges receiving output colours, and all exposed edges, are deleted. 
\end{itemize}

In each step, any random choices must follow the prescribed distributions conditional upon the graph $G_t$. The output of the algorithm is the vector $\vu$ whose $j$th component is the set of vertices receiving the $j$th output colour. 
\end{defn}

Actually, instead of having to clash colours $\propto$ and $\propto'$, we could have required all vertices involved in clashes to be assigned a single output colour $\propto$. However, we assign a second (transient) colour to vertices involved in the clashes described in the recolouring step as a signal to avoid querying their open adjacencies. 

The algorithms described in Section~\ref{s:introresults} did not use colours, but can be recast as local deletion algorithms. We next give an example of a local deletion algorithm that makes intrinsic use of colours, inspired by the algorithm in~\cite{Volec}.

 \bigskip
 
 \noindent {\bf A local deletion algorithm for cuts in cubic graphs}
\smallskip

Let $G=(V,E)$ be a graph. Given $A \subseteq V$, the (edge) cut induced by $A$ is the set of all edges in $E$ with one endpoint in $A$ and the other in $V \setminus A$. We are interested in the size of a maximum cut in a cubic graph. 

To describe the algorithm, we need two output colours, red and blue (as well as the clash colour $\propto$). To denote the transient colour of a vertex in the survival graph which has  $r$ of its neighbours in the input graph $G$  red, and $b$  blue,  we use $rb$, or $r,b$ if necessary for clarity. Thus, $00$ denotes the neutral transient colour, since initially every vertex has no blue and no red neighbours. For economy we dispense with separately recording the degree $3-r-b$, and we say that a vertex with transient colour $rb$ has type $rb$. At step $t$, the algorithm selects u.a.r.\ a vertex $v$ with type $rb$, where $rb$ is the type of highest priority that appears in the survival graph $G_t$, according to a given priority list. If $r>b$, $v$ is coloured blue and deleted from $G_t$, and the types of its neighbours are updated with an additional blue neighbour. If $r<b$, the same thing happens with blue replaced by red and, in the case $r=b$, $v$ is assigned red or blue uniformly at random and deleted from the survival graph, and the types of its neighbours are updated accordingly. This procedure may be described formally as a local deletion algorithm as follows. 

The selection rule chooses u.a.r.\ a vertex in $G_t$ from those whose type has highest priority according to the following prescription. Any type $rb$ with $r\ge 2$ or $b\ge 2$ has higher priority than the rest (any ordering within these gives the same result), and apart from this,
$$
 10 > 01 > 11 > 00
$$
where $ rb > r' b'  $ means that  $rb$ has higher priority than $r'b'$.

The local rule is based on a normal local subrule which, for a given query graph, always chooses to query $(1,\se)$, if it appears in $\ell_1$, and to terminate otherwise. In other words, after a vertex $v$ is selected, the local rule builds a copy of the query graph rooted at $v$ which consists of the single vertex $v$ and a multiset with the types of its neighbours in the survival graph.

If the type of $v$ is $rb$ and $r>b$, the recolouring rule $c$ assigns blue to $v$ and the colour $r',b'+1$ to each element $r'b'$ in the multiset associated with it. This accounts for the fact that after colouring $v$ blue, any neighboor of $v$ has one more blue neighbour, which will be reflected in its new colour in the recolouring step of the algorithm.  If $r<b$, it assigns red to $v$ and $r'+1,b'$ to each element $r'b'$ in the multiset associated with it. Finally, if $r=b$, the recolouring rule colours $v$ randomly with blue or red with uniform probability, and  recolours the open adjacencies  associated with it accordingly.

To refer to algorithms such as those in Section~\ref{s:introresults} that do not really need colours, we use the following.

\begin{defn}[Native local deletion algorithm]
This is a local deletion algorithm with only two transient colours, being neutral and $\propto'$, and three output colours. The set of vertices of output Colour 1 is interpreted as the {\em output set} $\U$, as usual $\propto$ denotes clash vertices, and the rest of the deleted vertices are of Colour 2.
\end{defn}

We adopt the convention of defining  native local deletion algorithms without colours in an obvious way: the local and recolouring rules simply determine which vertices and edges are deleted and which are added to the output set.

 \subsection{Chunky local deletion algorithms}\lab{s:chunky}

To facilitate analysis we can restrict the kind of algorithms under scrutiny, without sacrificing the power of the results  to the accuracy we are interested in.  A local deletion algorithm is said to be \emph{chunky} if the probability distribution associated with the selection rule $\Pi_t$ at step $t$ is defined as follows for every $t \in \{1,\ldots,N\}$: 
\begin{itemize}
\item[(1)] there are fixed real numbers $ p_{t,j} \in [0,1]$ for every type $j$, called the governing probabilities of the chunky algorithm; 

\item[(2)] 
the selection rule $\Pi_t$ is such that the probability of a nonempty set $S$ is given by $\pi_G(S)=\prod_{v \in S} p_{t,\d_G(v)}$. 

\end{itemize}
In other words, in a chunky local deletion algorithm, each vertex $v$ of $G_{t-1}$ with type $j$ is added to the set $S_t$ randomly, independently of all others, with probability $p_{t,j}$. This separation into $N$ chunks of time with a good deal of independence within each chunk is useful, especially as we will in practice keep $N$ fixed independently of the size of the input graph. If $\mathcal{T}=\{1,\ldots,\r\}$ denotes the set types, we call the matrix  $(p_{t,j})_{1\le t\le N, 1\le j\le \r}$  a {\em matrix of probabilities} of the algorithm, and  we define the {\em granularity} of the algorithm to be the maximum entry in this matrix. Such a matrix is general enough for any algorithm with a fixed number of colours that is intended to act only upon graphs with maximum degree $r$, as there is a bounded number of types in this case.

We illustrate using the procedure for finding an independent set in a graph $G$ discussed in Section~\ref{s:introresults}. This can be turned into a chunky local deletion algorithm if, instead of choosing a single vertex of minimum degree in the survival graph in each step, the selection rule produces a subset $S_t$ by including each such vertex with some probability $p$. Assuming there are no clashes, the survival graph $G_t$ is updated by incorporating the changes prescribed by each query graph separately. However, if there are clashes, some recolouring is involved. For instance, if two vertices in $S_{t }$ are adjacent, they will be coloured with $\propto$. In this way, the vertices with the output colour designating them for the final independent set will retain the property of actually being an independent set.  This is a workable approach if we ensure that clashes are rare, by choosing $p$ to be small. For large $n$, most vertices have output colours (rather than clash colours) by the end of the algorithm.

\section{Analysis of chunky local deletion algorithms}\lab{s:largegirth}

The aim of this section is to show that, for any fixed $r \geq 2$, chunky local deletion algorithms give essentially the same performance, for the random variables of interest, when applied to any $r$-regular graph of sufficiently large girth. This will be proved in two main steps. First we establish that the expected values of these random variables is the same for every $n$-vertex $r$-regular graph with girth at least $g$, provided that $g$ is sufficiently large in terms of the number of steps of the algorithm. We then show that the variables are sharply concentrated around these same expected values as long as the number of short cycles is `small,' which  holds trivially for graphs with large girth and is a well known property of random regular graphs. This establishes the desired connection between these two types of graphs.

For our purposes in this section, we are given an $r$-regular graph $G$ in which all vertices initially have neutral transient colour. We  perform the chunky algorithm making all random choices in advance, as follows. With each $t\in \{0,1,2,\ldots,N-1\}$ and each possible type  $j$, we associate a random set $S_{j}(t) \subseteq V$, where a vertex is placed in $S_j(t)$ with probability $p_{t,j}$. All these choices are made independently of each other and independently of the other sets $S_{j'}(t')$. In addition, for each copy $\phi_H$ of any query graph $H$ in $G$, we specify  an open adjacency (or $\halt$) according to the probability distribution specified by the local subrule $\phi$, and also, if  $\halt$ is specified and the recolouring rule $c$ is randomised, the values of the colours assigned by $c$. All of this  gives us a function $\Upsilon$ defined on all copies of query graphs in $G$. 

Given $\Upsilon$ and the sets $S_j(t)$, define $G_0=G$, and to determine $G_{t}$ from $G_{t-1}$ as follows. For  a vertex $v$ whose type in $G_{t-1}$  is $j$, place $v$ in the set $S_{t }$ if and only if $v \in S_j(t)$.
Then use the function $\Upsilon$ iteratively via the local subrule to obtain the query graphs (just as determined by the local rule in the chunky algorithm). Finally, apply steps (iii) and (iv) in the description of the local deletion algorithm,  using $\Upsilon$ to determine any random choices in the recolouring rule. This determines a `survival' graph $G_{t}$ and the colours of the vertices deleted.  Note that the value of $\Upsilon$ at any particular point in its domain is   used  at most  once.
This is because for each query graph produced by the local rule, there are two possibilities. Either the selected vertex has degree at least 1 and, according to the local subrule, at least one open adjacency is queried, with the corresponding edges subsequently deleted, or it has degree 0 and  is given an output colour and deleted according to the recolouring rule.

We may now   prove an upper bound on the ``speed'' at which random choices in some part of the graph may affect random choices made in another part of the graph. First, affix a label $\mathcal{L}^v=(\mathcal{L}^v_1,\mathcal{L}^v_2) $ to each vertex $v$ where $\mathcal{L}^v_1=\{t:v\in   S_j(t)$\}, and $\mathcal{L}^v_2$ consists of the action of $\Upsilon$ on all the  copies of query graphs whose root is $v$. Note that these labels are all determined before generating the process $(G_t)_{t\ge 0}$.

\begin{lemma}\lab{aux ind lab} Consider a chunky local deletion algorithm whose local rule has depth $D$. 
Let $G$ be a graph, let $G_t$ be the survival graph defined above,  and let $F$ be a subgraph of $G$. Then $F\cap G_t$, and the types of its vertices, are determined by the subgraph induced by the vertices whose distance in $G$ from $F$ is at most $2Dt$, together with the labels of those vertices, up to label-preserving isomorphisms.
\end{lemma}

\proof The nature of the change in type of a vertex $w$  at step $t$, and whether $w$ is deleted in step $t$,  is dependent only on the copies of query graphs existing in  $G_{t-1}$
 which $w$ is in or adjacent to, together with the labels on the root vertices of such copies. The set of such copies is determined by the graph $F'\cap G_{t-1}$, where $F' =\displaystyle{G[\{v:d_G(v,F) \leq 2 D \}]}$. The lemma now follows by induction on $t$ and the triangle inequality.\qed  

A \emph{rooted graph} is a graph with one vertex distinguished, and called the \emph{root}.  If the graph is a tree $T$, we call it a \emph{rooted tree}, and its height (i.e.\ maximum distance of a vertex from the root) is denoted $h(T)$. Let $T_{r,h}$ denote the \emph{balanced $r$-regular tree with height $h$}, that is, the rooted tree in which every non-leaf vertex has degree $r$ and every leaf is at distance $h$ from the root. 

For a vertex $v$ in a graph $G$, let the \emph{$s$-neighbourhood} of $v$ be the subgraph of $G$ induced by the set of vertices whose distance to $v$ is at most $s$, viewed as a rooted graph with root $v$. If a vertex does not survive in $ G_t$, we define its $s$-neighbourhood to be empty.

\begin{cor}\lab{gen ind lab}
Fix $r$ and $h$, and consider a chunky local deletion algorithm $\mathcal{A}$ whose local rule has depth $D$. Let $G$ be an $r$-regular graph  of girth greater than  $4Dt+2h$ and apply   $\mathcal{A}$  to $G$ with all vertices initially of neutral transient colour.  
\begin{description} 
\item{(i)}
Let   $F$ be a vertex-coloured version of a subgraph of $T_{r,h}$, and fix a positive integer $t$. Consider an embedding $\hat{T}$ of $T_{r,h}$ in   $G$. Then 
$
\pr(G_t\cap\hat{T}=F) $ 
(where equality requires the colours of $G_t$ to match the colours of $F$) is a constant depending on $\mathcal{A}$ and $r$  but independent of the choice of the graph $G$ and the embedding $\hat{T}$. 

\item{(ii)}
Let $T$ be a coloured rooted tree of height at most $h$ and fix a vertex $v\in V(G)$. Then  the probability that the $h$-neighbourhood of $v$ in $G_t$ is isomorphic to $T$ (where isomorphism requires the colours to be preserved) after $t$ steps is  a constant $p_{t, T,h}=p_{t, T,h}(\mathcal{A},r)$ independent of the choice of the graph $G$ and the vertex $v$.
\end{description}
\end{cor}

\proof Define   $G_{\hat{T}}$ to be the subgraph of $G$ induced by all vertices within distance $2  D  t$ of $\hat{T}$. 
By Lemma \ref{aux ind lab}, for an embedding $\hat{T}$ of $T$ in $G$, the intersection of $ G_t$ and $\hat{T}$ is determined by the label-preserving isomorphism type  of the subgraph induced by the vertices of $G_{\hat{T}}$.
By the girth hypothesis, the  graph  $G_{\hat{T}}$  is a tree whose isomorphism type (ignoring the vertex labels) is independent of the choice of $G$ and $\hat{T}$. The labels are assigned independently with the same probability distribution to each vertex. The result follows for (i). Then (ii) follows by applying (i) to all graphs $F$ whose component containing the root of $T_{r,h}$ is  isomorphic to $T$ (with matching colours), and such that $v$ is the root of $\hat T$, and summing the probabilities.\qed

Instead of proving facts about the number of vertices with each output colour, it will be convenient to consider output functions that are slightly more general.
\begin{defn}[Output function, active copy]
An {\em output function} $W$ is a function determined by fixed numbers $c_{H,L,J}$  via the recurrence  
\begin{equation}\lab{def_output}
W(t+1)=W(t)+\sum_{H,L,J} c_{H,L,J}W_{H,L,J} \quad (t\ge 0),\quad  W(0)=0. 
\end{equation}
Here $ H\in {\cal Q}$, $L\in {\cal L}_2$  where ${\cal L}_2$ is the (finite) set of possible values of the label ${\cal  L}^v_2$,  $ J\subseteq V(H)$, and 
 $W_{H,L,J}$ is determined by the set of copies of $H$ generated in  step $t+1$ of the algorithm, as follows. We call such copies of $H$ {\em active}, and we refer to the vertices $\psi(V(H))$ as the {\em vertices} of the copy $\psi$. Then $W_{H,L,J}$ is the number of active copies $\psi$ of $H$ for which ${\cal L}^v_2(\psi)=L$, where $v$ is the root vertex $\psi(1)$, and such that the set of clash vertices of $\psi$ is precisely $\psi(J)$.
\end{defn}

Note that the function $W(t)$ that counts vertices of a given output colour is an example of an output function. This is obtained by setting  $c_{H,L,J}$ equal to the number of vertices coloured with that output colour when the query graph $H$ has ${\cal L}_2^v=L$ for $v$ and the set of clash vertices is $J$.
\begin{thm}\lab{large girth}
Let $G$ be an $r$-regular graph on $n$ vertices with girth larger than $4DN+\max\{2,6D\}$, all of whose vertices initially have neutral transient colour. Let $Y_k(t)$ denote the number of vertices of type $k$ in the survival graph $G_t$ in a chunky local deletion algorithm $\mathcal{A}$ whose local rule has depth $D$, where $0 \leq t \leq N$. Also let $W(t)$ be an output function. Then, for each such $t$,  the vector
$$
s_{r,\mathcal{A}}(t)=s_{r,\mathcal{A}}(t,n) = \big(\ex Y_1(t),\ldots, \ex Y_R(t),\ex  W(t)\big)
$$
is independent of the choice of $G$.
\end{thm}

\proof  
Applying Corollary~\ref{gen ind lab}(ii) with $h=1$, we may sum $p_{t,T,1}$ over all $T$ with root of type $k$,  to obtain the probability, independent of $G$, that any given vertex of $G$ is of type $k$ in $G_t$. Multipying by $n$ gives $\ex Y_k(t)$, which proves the first part, as the girth is at least $4DN+2$.   The proof for $W(t)$ is similar, but requires girth at least at least $4DN+6D$. The expected value of $W(t+1)- W(t)$ is fixed provided that, for each copy $\psi$ of any query graph $H$, the probability that $\psi$ is active  and $L={\cal  L}^v_2$, where $v$ is the root of $\psi$, and additionally $J$ is the set of clash vertices in this copy in step $t+1$, is a fixed number given $H$, $L$ and $J$. This is true for $t\le N-1$ by Corollary~\ref{gen ind lab} with $h=3D$, since the intersection of the $(3D)$-neighbourhood of $v$ with $G_{t-1}$ determines all these things.
\qed
 
The next result shows that the key random variables of a chunky local deletion algorithm are concentrated near the values suggested by the constants given in the previous two results, provided the input graph has few short cycles.

 \begin{lemma}\lab{thm conc}
Given a positive integer $N$, let $G$ be an $r$-regular graph with $n$ vertices, at most $\Theta<n^{2/3}$ of which are in cycles of length less than or equal to $4DN+6D+2$.  Consider $N$ steps of a chunky local deletion algorithm $\mathcal{A}$ of depth $D$ applied to $G$. Then the following hold.
\begin{itemize}
\item[(i)]
 Given $t\in \{0,\ldots,N\}$ and a coloured rooted tree $T$ of height at most  $h\le  2D(N-t)+3D+1$,  let $A_{t, T,h}$ denote the number of vertices in $G_t$ whose coloured $h$-neighbourhood in $G_t$ is isomorphic to $T$. Then, for some   positive constant $C$ determined by the local and recolouring rules, the degree $r$ and the number of steps $N$, 
$$\pr\left(\left|A_{t, T,h}- np_{t,T,h}\right| \geq   C n^{2/3}\sqrt{ \Theta+1}\big)\right)  \le n^{-1/3},$$
where $p_{t,T,h}$ is defined in Corollary~\ref{gen ind lab}(ii).

\item[(ii)] Define $Y_k(t)$, $W(t)$ and $s_{r,\mathcal{A}}(t,n)$ as in Theorem~\ref{large girth}. Set $\vs(t)=\big(Y_1(t),\ldots, Y_R(t),W(t)\big)$ and, given $C>0$, consider the event $F(C)$ that $\left|\left|\vs(t)-s_{r,\mathcal{A}}(t,n)\right|\right|_{\infty}\geq  C n^{2/3}\sqrt{ \Theta+1}$ for some $t \in \{0,1,\ldots,N\}$. Then, for some constants $C$ and $C'$  determined by $r$, $N$ and the local and recolouring rules,
$$\pr\big(F(C)\big) \leq C'  n^{-1/3}.$$
\end{itemize}
\end{lemma}

\proof 
We first prove part (i). A vertex in $G$ is said to be \emph{good} if its distance from every vertex lying in a cycle  of length at most $4DN+6D+2$ is larger than $2D(N+1)$, and   \emph{bad} otherwise. Note that the set $B$ of bad vertices of $G$ has size at most $C\Theta$ where $C$ is   constant   given $r$, $D$ and $N$.

Let $t$ and  $T$ be as in the hypothesis. By Corollary \ref{gen ind lab} and linearity of expectation,
$$\E A_{t, T,h}=(n-|B|)p_{t,T,h} +|B|O(1)=np_{t,T,h} + O(\Theta)$$
for fixed $r$, $D$ and $N$.
On the other hand,  
$$
\Var(A_{t, T,h})=\sum_{v \in V} \Var(X_{v, T,h}) + \sum_{v \neq w \in V}\Cov(X_{v,  T,h},X_{w, T,h}),
$$
where $X_{v, T,h}$ is the indicator random variable for the event that $v$ has coloured $h$-neighbourhood $T$ after $t$ steps of the algorithm. 

Being an indicator, $X_{v, T,h}$ has variance at most 1. If the distance between $v$ and $w$ is greater than $4D(N+2)$, then the variables $X_{v, T,h}$ and $X_{w, T,h}$ are independent, and hence their covariance is 0. For all other pairs, the covariance has absolute value at most 1. Thus
$$
\Var(A_{t, T,h})=O(n) + O(n|B|)=O(n(\Theta+1)).
$$
By Chebyshev's inequality, 
 $\pr\big(|A_{t,T,h}-\E A_{t,T,h}| \ge   \sqrt{n^{1/3}\Var(A_{t,T,h})}\big)  \le n^{-1/3}$, giving (i).
 
For  part (ii), we argue in a fashion similar to the proof of Theorem~\ref{large girth}. 
$Y_k(t)$ is just  a sum of $A_{t,T,1}$ over appropriate trees $T$, so the required concentration these components of $\vs(t)$ follows from (a). For the final component, 
$W(t)$, one can define an indicator variable to denote that a copy $\psi$  of a query graph  $H$ is active,  with given value $L$ of  $ {\cal  L}^v_2$ where $v$ is the root of $\psi$, and given set $J$ of   clash vertices in this copy in step $t+1$. Using the fact from (i) that the numbers of vertices with given neighbourhoods are concentrated, and again arguing as for (i) using Chebyshev's inequality, the sum of  these indicators, which equals $W_{H,L,J}$ as in~\eqn{def_output}, is concentrated, to the extent claimed, except for an event with probability $n^{-1/3}$. Summing over all $H$, $L$ and $J$ and applying the union bound shows that $W(t+1)-W(t)$ is similarly concentrated, up to a constant factor in the error, with probability   $O(n^{-1/3})$.  Doing the same for $t$ gives the result required to complete (ii).
 \qed

The conclusion  of Lemma~\ref{thm conc} could easily be strengthened  by using a slightly more complicated inductive argument involving Azuma-type inequalities, but this is not needed.

Lemma~\ref{thm conc}(ii) may be applied directly to random regular graphs.
\begin{thm}\lab{conc rreg}
Let  $\mathcal{A}$  be   a chunky local deletion algorithm  performing $N$ steps applied to a random $r$-regular graph $G \in \mathcal{G}_{n,r}$. The vector $\vs(t)=\big(Y_1(t),\ldots, Y_R(t),W(t)\big)$ produced after $t$ steps of the algorithm is a.a.s.\ within
$C'_1 n^{2/3}$
of $s_{r,\mathcal{A}}(t,n)$ in $L^\infty$-norm, where $C_1'$ is a constant for given $N$, $\r$ and $D$,  uniformly over all $t \in \{0,\ldots,N\}$.
\end{thm}
\proof In a random regular graph, the expected number of cycles of length at most $4D(N+1)$ is bounded (see~\cite{BB} or~\cite{models}), so the result immediately follows from Lemma~\ref{thm conc}(ii). \qed

 Theorem~\ref{conc rreg} implies that  the output vector produced by a chunky local deletion algorithm $\mathcal{A}$ applied to a random $n$-vertex $r$-regular graph and the vector of expected output of $\mathcal{A}$ applied to an $n$-vertex $r$-regular graph $G$ with sufficiently large girth a.a.s.\ differ by $o(n)$. A useful consequence  is that, given a chunky local deletion algorithm $\mathcal{A}$, we may estimate the value of $s_{r,\mathcal{A}}(n)$ by analysing the performance of $\mathcal{A}$ on random regular graphs. This is a crucial factor in our approach, as it allows us to make use of the powerful machinery already developed for analysing random regular graphs. This will be addressed in the next section.

 To obtain our results, we need to limit the number of clashes in order to obtain a connection between the results of a chunky local deletion algorithm and other local deletion algorithms with the same local rule and recolouring function.
We define a {\em pre-clash} in a chunky local deletion algorithm of depth $D$ to be a pair of vertices of distance at most $2D$ apart, which are both chosen for inclusion in the set $S_t$ in a given step $t$. The number of clashes will be bounded above by a constant times the number of pre-clashes.
 
Recall that the  granularity of a  chunky algorithm is the maximum entry in its matrix of probabilities.
\begin{lemma}\lab{l:clashes}
The expected number of pre-clashes in step $t$ of any chunky algorithm of depth $D$ and granularity $\eps$, applied to a graph on $n$ vertices and maximum degree $r$, is $O(\eps^2n)$, where the implicit constant depends only on $D$ and $r$. 
\end{lemma}
\proof
The number of pairs of vertices of distance at most $2D$ is $O(n)$, and the probability that both vertices are chosen in $S_t$ is at most $\eps^2$, so this  follows immediately from the union bound.
\qed

%%%%%%%%%%%%%%%%%%%%%%%%%%%%%%%%%%%%%%%%%%%%%%%%%%%%%%%%%%%%%%%%%%%%%%%%%%
\section{Explicit bounds from chunky algorithms}
\lab{s:explicit} 

In light of Theorem~\ref{large girth}, to obtain explicit bounds on the  value of an output function  when a local deletion algorithm is applied to a graph with sufficiently large girth, we can try to determine the value of the vector $s_{r,\mathcal{A}}(t)$ in Theorem~\ref{large girth}.
Owing to Theorem~\ref{conc rreg}, it suffices to determine this vector for random regular graphs. In this section we address this issue, developing what will eventually be part of some machinery which lets us translate some existing results on random regular graphs to results for regular graphs of large girth.

 We will need to analyse the behaviour of a chunky local deletion algorithm $\mathcal{D}$ on an $n$-vertex random $r$-regular graph (as usual, $rn$ is assumed even and the probability space of all random regular graphs is denoted by $\mathcal{G}_{n,r}$).  To this end, we shall use a well known model introduced by Bollob{\'a}s~\cite{BB}, which can be described as follows. Consider $rn$ points in $n$ buckets labelled $1,\ldots,n$, with $r$ in each bucket, and choose uniformly at random (u.a.r.) a \emph{pairing} $P=a_1,\ldots,a_{rn/2}$ of the points such that each $a_i$ is an unordered pair of points, and each point is in precisely one pair $a_i$. We use  $\mathcal{P}_{n,r}$ to denote this probability space of random pairings. Each pairing corresponds to an $r$-regular pseudograph  (i.e., graph with possible loops and multiple edges) with vertex set $1,\ldots,n$ and with an edge for each pair. A pair with points in buckets   $i$ and $j$ gives rise to an edge joining vertices $i$ and $j$. (If $i=j$ a loop is formed.) A straightforward calculation shows that the simple $r$-regular graphs (i.e.\ those with no loops or multiple edges) on $n$ vertices are produced each with the same probability, and hence an $n$-vertex $r$-regular graph can be produced u.a.r.\  by choosing a pairing u.a.r.\ and rejecting the result if it has loops or multiple edges. The probability that a random pairing produces an $r$-regular graph tends to the positive constant $e^{(1-r^2)/4}$ as $n$ tends to infinity (Bender and Canfield~\cite{BC}). There is an obvious generalisation of $\mathcal{P}_{n,r}$ to $\mathcal{P}(\dv)$, where   bucket   $i$ contains $r_i$ points. Provided that the elements of $\dv$ are bounded, the probability that the resulting graph is simple is bounded below. Conditioning on this event gives the model $\G(\dv)$ of uniformly random graphs with degree sequence $\dv$.  (See~\cite{models} for more details.) 

The pairing model may be redefined slightly by specifying that the pairs are chosen sequentially: the first point in the next random pair can be selected using any rule whatsoever, depending only on the choices made up to that point together with some independent random input, as long as the second is chosen u.a.r.\ from the remaining points. This is called \emph{exposing} the pair containing the first point, and this property is sometimes called the \emph{independence property} of the pairing model. For instance, one can require that the next point chosen to be exposed comes from the lowest-labelled bucket available, from a bucket with fewest unpaired points, or from the bucket containing one of the points in the previous completed pair if any such points are still unpaired. In particular, for any algorithm applied to the final random graph, the process for generating the pairs can be determined by the order in which the algorithm queries the edges of the graph. We call the process that generates the pairs governed by an algorithm in this way a {\em pairing process}. For 
simplicity, we will often refer to a pairing as a (pseudo)graph and to a bucket as a vertex. The fact that we are dealing with a pairing instead of the underlying pseudograph will be clear from context.

We need to provide a  definition of local deletion algorithms, in particular local rules, so as to act upon pairings. There is an obvious way to do this, provided the set $\cal Q$ of query graphs includes all coloured pseudographs. We will assume this property of $\cal Q$ henceforth.  Recall that, at every step $t$, once a vertex $v$ is selected, the algorithm queries vertices in the survival graph $G_{t-1}$ to build a query graph 
rooted at $v$. This query graph then determines (perhaps with randomisation) which vertices are added to the output set and which are deleted from or recoloured in $G_{t-1}$ to create $G_t$.  In the context of the pairing process described in the previous paragraph,   querying  an open adjacency $(i,\se)$ is equivalent to specifying, for an unpaired point $j$ in bucket $i$, the type $\d(j)$ of the bucket containing the mate of $j$,  while querying an open adjacency $(i,\d)$  corresponds to selecting any point $j$ in bucket $i$ with $\d(j)=\d$ and choosing an unpaired point $k$ in a bucket of type $\d$ to complete a pair with $j$.  The point $k$ is chosen randomly among all unpaired points in buckets that had type $\d(j)$ at the moment that bucket $i$ was queried. The random choice of $v$ can similarly be incorporated, creating a pairing process that generates the edges of the random pseudograph just as the algorithm requires them. Instead of a survival graph, we have for pairings a survival pairing, containing just the unexposed pairs in the pairing process after $t$ steps of the algorithm, and a corresponding   survival pseudograph, which we call $G_t$. The clash step is interpreted in an obvious way. Note that a vertex with a loop is will automatically be a clash.

It is easy to check that, conditional on the pseudograph of the pairing being simple, the above description of the local rule in the pairing corresponds to the local rule as applied to the graph corresponding to the pairing. In view of the correpondence of models described earlier, this immediately gives the following useful result.
\begin{lemma}\lab{l:pairings to graphs}
Suppose that when a local deletion algorithm is applied to a random pairing with $n_k=n_k(n)$ vertices of type $k$ for each $k$, the survival pseudograph $G_t$ has property $P$ with probability   $p(n)$. If the same algorithm (i.e.\ with the same selection and local rules and recolouring function) is applied to a random graph   with $n_k=n_k(n)$ vertices of type $k$ for each $k$, then the survival graph after $t$ steps  has property $P$ with probability at most $cp(n)$, where $c$ is a constant depending only on the maximum degree of the input graph.
\end{lemma}

When applying local deletion algorithms to a pairing, we will often for simplicity speak of applying it to the associated pseudograph. The next result is our basic tool describing how the proportion of vertices of each type in the survival graph is expected to change due to an application of the local rule to a single vertex in the survival graph. For a coloured pseudograph $G$ on $n$ vertices, let $n_k$ denote the number of vertices of type $k$ in $G$, define the vector $\tilde{\nv}=(n_1/n,\ldots ,n_\r/n)$, and let the degree of a vertex of type $k$ be denoted $\dk$.

\begin{lemma}\lab{l:trend}
Let $r$ and $D$ be positive integers. Given a local rule $L_\phi$ of depth $D$ and a recolouring rule $c$, let $\r$ denote the number of types when a local deletion algorithm with these rules is applied to graphs with degrees bounded above by $r$.  Let $i,k \in \{1,\ldots,\r\}$ be types and let $H$ be a query graph. Then there exist functions $f_{k,i}$ and $g_{H,i} \colon \mathbb{R}_{\geq 0}^{\r} \rightarrow \reals$ which, for all $\eps>0$, are  Lipschitz continuous in $D(\eps):=  \{(\tilde{n}_1,\ldots,\tilde{n}_\r) \in \reals_{\geq 0}^{\r} \colon  \sum_{k=1}^{\r} \dk \tilde{n}_k \geq \eps\}$,  such that the following is true. Let $G$ be a random pseudograph with degree sequence $(r_1,\ldots,r_n)$ generated by the pairing model, where $0 \leq r_k \leq r$ for all $k$, and fix a colouring of $G$. Let $n_k$ denote the number of vertices of type $k$. Consider $i \in \{1,\ldots,\r\}$ such that $n_i>0$ and fix a vertex $v$ of type $i$ in $G$. Let $H_v$ be the (random) query graph obtained by one application of $L_\phi$ to $v$ in $G$, and let $G'$ be the survival graph obtained from $G$ after the clash and recolouring steps. Assume that $M=\sum_{j=1}^n r_j=\sum_{k=1}^{\r} \dk n_k>\eps n$ for some $\eps>0$. The following hold with the constants implicit in $O()$ terms independent of the~$r_k$.
\begin{itemize}
\item[(i)] For any query graph $H$,
$$\pr(H_v=H)= g_{H,i}(\tilde{\nv})+O(1/n).$$

\item[(ii)] The number $Y_k$ of vertices of type $k$ in $G'$ satisfies
$$\E(Y_k)=n_k+f_{k,i}(\tilde{\nv})+O(1/n).$$

\end{itemize}
Moreover, the expected number of vertices that are assigned colour $\propto'$ in the recolouring step is $O(1/n)$.
\end{lemma}
\proof Recall that the survival graph $G'$ is derived from $G$ by deleting all edges in the copy of the query graph $H_v$ obtained through an application of the local rule $L_\phi$ to $v$ in $G$, and by recolouring or deleting vertices in this copy, as well as recolouring vertices associated with the open adjacencies of $H$.  

Before addressing the claims directly, we find it useful to bound the (conditional) probability that either of the following two events occurs: (a) one of the open adjacencies in $H_v$ is incident with a vertex within $H_v$ (that vertex is then a clash); (b) two or more open adjacencies in $H_v$ are incident with the same vertex (that vertex is then assigned colour $\propto'$).   As noted before, we may assume that, when a pair is exposed, the mate of the starting point is distributed u.a.r.\ over the remaining available points. Because of this, the (conditional) probability of  event (a) is bounded above by $s^2/S$, where $s$ is the number of points in the vertices of $H_v$ and $S$ is the number of unpaired points in $G-V(H_v)$. Clearly, $s\leq r|V(H_v)|$, $S \geq M-r |V(H_v)|$ and $|V(H_v)| \leq m_{r,D-1} \le 1+r^{D}$, where $m_{r,D-1}$ denotes the number of vertices in a balanced $r$-regular tree of height $D-1$. Hence the probability of (a) is bounded above by $r^2 m_{r,D-1}^2/(M-r m_{r,D-1})=O(1/M)=O(1/n)$. Similarly, the (conditional) probability of  event (b) is $O(1/n)$, as the number of points in buckets outside $H_v$ is $\Omega(M)$ and hence the probability two open adjacencies have mates in the same bucket is $O(1/M)$. Observe that, since the recolouring step assigns colour $\propto'$ only if (b) occurs, the expected number of vertices that are assigned colour $\propto'$ in the recolouring step is $O(1/n)$, which gives the final statement of the lemma.

Next, we show that part (ii) is a direct consequence of part (i). We argue at first that, conditional  upon the  query graph $H_v$, the quantity  $\E(Y_k)-n_k$ is determined up to a $O(1/n)$ term.   Indeed, consider the vertices that have type $k$ in $G$. The type of one of these vertices may be different from $k$ in $G'$ if it is in the copy of the query graph $H_v$, or if it is the end of an open adjacency (whose type is necessarily $k$)  in this copy of  $H_v$. The expected number of these vertices that change type is determined by the action of the (possibly randomised) recolouring function $c$ on $H_v$, independently of  the degree sequence. On the other hand, consider the  vertices that have type $k$ (and whose colour is not $\propto'$) in $G'$ but do not have this type in $G$. These vertices come from two sources: those that have degree $\dk+1$ in $G$ and are the end of an open adjacency of $H_v$, which in addition is assigned the colour of $k$ by the recolouring function $c$, and those with the correct degree that lie in $H_v$ and are assigned the colour  of $k$ by $c$. Again, the expected number is a function of the recolouring rule. In this discussion we have somewhat ignored events (a) and (b), but their probabilities are determined by the degree sequence of the graph, and the expected changes in the quantities due to their occurrence are similarly determined up to a $O(1/n)$ error.

As a consequence, $\E(Y_k)-n_k=\Delta_k(H_v)+O(1/n)$ for a function $\Delta_k$ depending only on the rules, not the degree sequence. Hence the expected value of $Y_k$ is given by
\begin{eqnarray*}
\ex(Y_k) 
&=&n_k + \sum_{H \in \mathcal{Q}_{D}} \pr(H_v=H) \cdot \Delta_k(H) + O(1/n),
\end{eqnarray*}
as the number of query graphs in $\mathcal{Q}_{D}$ with maximum degree bounded by $r$ is a constant depending on $r$, $D$ and the number of colours available. The claim of part (ii) thus follows from part (i).

Let $H_{v,j}$ denote the (random) query graph obtained after $j$ steps of the local subrule when applying the local rule starting with a vertex $v$ of type $i$ in $G$.  To prove part (i), we show by  induction on $j$  that, for a given $H$, $\pr( H_{v,j}=H) =  g^{(j)}_{H,i}(\tilde{\nv})+ O(1/n)$, where $g^{(j)}_{H,i}$ is Lipschitz continuous on $D(\eps)$. After the local rule is initiated, but before any application of the local subrule, the only possible query graph $H$ is a singleton $\{1\}$ with a multiset $\ell_1$ containing $i$ copies of $\se$. This starts the induction off with $j=0$. For the inductive step, we may assume that $H$ is a tree; otherwise, at some point in the process, the query operation would have exposed an open adjacency incident with a vertex within the query graph, which has probability $O(1/n)$ as explained above. Hence, the query graph $H$ has been  derived in one of two possible ways. In the first, it is derived from a query graph $H'$ by adding the vertex that has largest label in $H$, which we may denote by $u$. At the same time, the multiset associated with the unique neighbour $w$  of $u$, is adjusted appropriately. There is a unique $H'$ with this property. The second possibility is to derive it from a query graph $H^\star$  by replacing an occurrence of $\se$ in a multiset associated with a vertex in $H^\star$  by a vertex type. 

In the first case, by induction $\pr(H_{v,j-1}=H') = g^{(j-1)}_{H',i}(\tilde{\nv})+O(1/n)$. The probability that $H$ was obtained from $H'$ is the probability that the open adjacency $(w,\d(u))$ is queried, where $\d(u)$ is the type of the image of $u$ in $G$, which is equal to $\phi_{H'}(w,\d(u))$. In the second case, for each element in the multiset $S(H)=\{(w,\tau) \colon w \in H \textrm{ and }\se \ne \tau \in \ell_w\}$, let $H^{\star}_{w,\tau}$ be the same as $H$ but with one occurrence of $\tau$ replaced by $\se$ in $\ell_w$. Then    $\pr(H_{v,j-1}=H^{\star}_{w,\tau}) =g^{(j-1)}_{H^{\star}_{w,\tau},i}(\tilde{\nv})+O(1/n)$. Moreover, the probability that the open adjacency $(w,\se)$ is queried is $\phi_{H^{\star}_{w,\tau}}(w,\se)$, while the probability that the outcome of this query is $\tau$ is, in the pairing model, equal to 
$
d(\tau)n_{\tau}\big(\sum_{s=1}^{\r} d(s)n_s\big)^{-1}+O(1/n).
$
This implies that the function $g_{H,i}^{(j)}$ defined by 
$$g_{H,i}^{(j)}(\tilde{\nv})=g^{(j-1)}_{H',i}(\tilde{\nv})\phi_{H'}(w,\d(u))+\sum_{(w,\tau) \in S(H)}g^{(j-1)}_{H^{\star}_{w,\tau},i}(\tilde{\nv})\phi_{H^{\star}_{w,\tau}}(w,\se)\frac{d(\tau)n_{\tau}}{\sum_{s=1}^{\r}d(s)n_s} $$
satisfies the required properties for the induction to go through.

To conclude the proof of (i), note that  
$$
\pr(H_v=H) = \sum_{j \geq 0}  \pr(H_{v,j}=H) \phi_H(\halt). \qed
$$

For later reference, note that the proof of Lemma~\ref{l:trend} implies that the way a local rule acts on query graphs with cycles is irrelevant asymptotically. Hence, every extension of a local deletion algorithm for graphs to one for pseudographs gives the same functions $f_{k,i}$. 

The next lemma is very useful for explicit evaluation of the functions $f_{k,i}$, and rests heavily on the fact that we only use deletion algorithms.  
\begin{lemma}\lab{l:survival} When a local deletion algorithm is applied to a random pairing  as in Lemma~\ref{l:pairings to graphs}, 
after each step the survival pairing, conditional on history of the algorithm, is distributed u.a.r.\ given the degrees and transient colours of the surviving vertices. 
\end{lemma}
\proof
The algorithm can be formulated as a pairing process in which the pairs in the survival pairing consist precisely of the pairs that are not yet exposed in the process. So this lemma  follows immediately from the independence property of the pairing process stated above. \qed

Note that one corollary of this lemma is a fact that is in any case rather obvious: the transient colours do not affect the distribution of the survival pairing.

Lemma~\ref{l:trend} refers to one step in a local deletion algorithm, and we are now ready to extend this to analysis of a complete algorithm. Assume that the selection step preceding the exploration and recolouring steps is that of a chunky algorithm where vertices of type $i$ are chosen with probability $p_i$, and we are to perform step $t$ in a local deletion algorithm with given local and recolouring rules. The expected number of vertices of degree $i$ selected for the set $S_t$ is $n\tilde{p}_i$  where $\tilde{p}_i= p_i n_i/n$, with $n_i $ being  the number of vertices of type $i$ in the current graph.  The expected change in $n_i/n$  suggested by Lemma~\ref{l:trend} is thus approximately  $ \sum_{i=0}^r\tilde p_i  f_{j,i}(\tilde{\nv})$  (for $0\le i \le r$). 

With a slight abuse of notation, we use $W$ to denote the random value of an output function, and also use a phrase such as `let $W$ be an output function' to denote the choice of  the constants $c_{H,L,J}$ in~\eqn{def_output}. In the remainder of the paper, assuming that local, selection and recolouring rules are given and letting $W$ be an output function, we analogously  define the function
\begin{equation}\lab{def_output_func}
f_{\r+1,i}(\tilde{\nv})=\sum_{H \in \mathcal{Q}_i,L} c_{H,L,\emptyset} \, a_{H,L}  \, g_{H,i}(\tilde{\nv}), 
\end{equation}
where $\mathcal{Q}_i$ denotes the set of query graphs whose root vertex has type $i$, $a_{H,L}$ is the probability that $\mathcal{L}_2^{(v)}=L$ given that the query graph is $H$, and $g_{H,i}$ is the function given in Lemma~\ref{l:trend}(i). We restrict to the case $J=\emptyset$ because the influence of clashes is negligible, as observed in Lemma~\ref{l:trend}.  

The next result asserts that the vertex degree counts, and output function, of an application of a chunky local deletion algorithm a.a.s.\ approximate the solution of the difference equation suggested by the expected changes estimated above. The error in the approximation will depend on the maximum entry $\eps$ in the probability vectors. As the proof reveals, for fixed $N$ an error of order $\eps^2 N n$ is unavoidable to the first degree of approximation due to the occurrence of clashes. We impose the condition $N\le C/\eps$ to achieve eventually an error   $O(\eps n)$. Recall that $\G(\dv)$ denotes the model of random graphs with degree sequence $\dv=(r_1,\ldots, r_n)$, with the uniform distribution. 
\begin{thm}\lab{chunkyranreg}
Consider a chunky local deletion algorithm applied to $G\in \G(\dv)$,  all of whose vertices initially have neutral transient colour, for $N$ steps with matrix of probabilities $ Q=(p_{t,j})$, and let $W$ be an output function. Let $n_k$ denote the number of vertices of type $k$ in $G$, and let $\eps = \max_{t,j}p_{t,j}$. For $C$ a fixed constant and $N\le C/\eps$, consider the quantities $z_j(t)$, where $j \in \{1,\ldots,\r+1\}$ and $t \in \{0,\ldots,N\}$, given iteratively by
$$
z_j(t)= z_j(t-1) + \sum_{i=1}^\r p_{t,i}z_i(t-1) f_{j,i}(\zv(t-1)),
$$ 
where $z_j(0)=n_j/n$ for  $1\le j\le \r$, $z_{\r+1}(0)=0$, and the functions $f_{j,i}$ are those appearing in Lemma~\ref{l:trend}. Then, a.a.s.\ the number of vertices of type $j$ in the survival graph is $nz_{j}(t)+O(\eps n)$, for all $j \in \{1,\ldots,\r\}$, and the value of the output function $W(t)$ at step $t$ is $nz_{\r+1}(t)+O(\eps n)$, uniformly for $0 \leq t\leq N$. Here the constant implicit in $O(\eps n)$ depends on $C$ and on the local and  recolouring  rules but is otherwise independent of $N$ and of the $ p_{t,i}$ (and in particular of $\eps$). Furthermore it is independent of $\dv$ provided $\max r_i\le r$ for some fixed $r$.
\end{thm}
%%%%%%%%%%%%%%%%%%%%%%%%%%%%%%%
Part of the proof of this theorem uses a deprioritised algorithm, as defined in the next section, so we defer the proof until then. We note here that, if so desired,   `a.a.s.'\ in the conclusion can be improved to `with probability $1-\exp(-n^c)$ for any $c<1/3$', by using~\cite[Theorem 6.1]{desurvey} for the main application of the differential equation method inside the proof of Theorem~\ref{chunkyranreg}. However this would not improve the  later  consequences for graphs of large girth.

Suppose that the matrix of probabilities of a chunky local deletion algorithm satisfying the statement of Theorem~\ref{chunkyranreg} is such that, for each $i$, the quantities $p_{t,i}/\eps$  can be interpolated by a fixed piecewise Lipschitz continuous function $p_i(x)$ whose domain is rescaled to $[0,C]$, that is,  $p_{t,i}/\eps = p_i(\eps (t-1))$ for all relevant $t$. Then the difference equation in Theorem~\ref{chunkyranreg} defines the approximate solution by Euler's method, with step size $\eps$, of the differential equation system
\bel{desystem0}
y_j'(x)= \sum_{i=1}^\r  p_i(x)y_i(x) f_{j,i}(y_1,\ldots, y_\r) \quad (1\le j\le \r+1)
\ee
on the interval $[0,C]$ with initial conditions
\bel{deics}
y_j(0)=n_j/n\  (1\le j\le \r),\  y_{\r+1}(0)=0,
\ee
where $n_j$ is the number of vertices of type $j$ in the input graph for the algorithm. Here and in other differential equations in our paper, differentiation is with respect to the variable $x$. 

The next result may be viewed as converse to Theorem~\ref{chunkyranreg}, as it shows that, for any system of differential equations of the form~\eqref{desystem0} arising from a local rule $L_\phi$ and a recolouring rule $c$, there is a chunky local deletion algorithm with such rules whose performance on large graphs is a.a.s.\ determined within a small error by this system of differential equations. Note that, although the functions $p_i$ were motivated by comparison with probabilities, they do not need to be bounded above by 1.

\begin{thm}\lab{chunky des}
Fix $r>0$.  Consider a local rule $L_\phi$, a recolouring rule $c$ and a set of non-negative piecewise Lipschitz-continuous functions $p_i(x)$ $(1\le i\le \r)$ on $[0,C]$ for some $C>0$.  Define $f_{j,i}$  to be the functions in Lemma~\ref{l:trend}, and let $f_{\r+1,i}$ be defined with respect to an output function $W$ as in~\eqref{def_output_func}. Let the functions $y_{j}$ be determined by the solution of~\eqn{desystem0} with initial conditions~\eqn{deics}. Then the following hold.  
 
\begin{itemize}

\item[(i)]  For all $\eps'>0$  there is an $\eps>0$ with $\eps<\eps'$, and a chunky local deletion algorithm $\mathcal{A}(\eps')$ with local rule $L_\phi$, recolouring rule $c$ and of granularity at most $\eps'$ as follows. If $\mathcal{A}(\eps')$ is applied to a random $n$-vertex graph $G\in \G(\dv)$ with $\max r_i \le r$,  then a.a.s.\ as $n\to\infty$, the number $Y_j(t)$ of vertices of type $j$ in the survival graph $G_t$ satisfies $\left|Y_j(t)-ny_{j}(\eps t) \right|<\eps' n$ and $\left| W(t)-ny_{\r+1}(\eps t) \right|<\eps' n$,  for all $1 \leq j \leq \r$ and $t \in \{0,\ldots,\lfloor C/\eps \rfloor\}$. 

\item[(ii)] For all $\eps'>0$, there  exists  a positive constant $g$ such that the chunky local deletion algorithm $\mathcal{A}(\eps'/2)$ of (i) has the following property.  Let $G_n$  be an $n$-vertex, $r$-regular graph with girth at least $g$. Let $Y_j(t)$ denote the number of vertices of type $j$ in the survival graph and $W(t)$ the value of the output function $W$, after $t$ steps of the algorithm. Then, with $\eps$ as in (i), $\left|\E(Y_j(t))-ny_{j}(\eps t)\right|<\eps'  n$ ($1 \leq j \leq \r$) and $\left|\E(W(t))-ny_{\r+1}(\eps t)\right|<\eps'  n$, for all $t \in \{0,\ldots,\lfloor C/\eps \rfloor\}$.
\end{itemize}
\end{thm}

\proof
To prove part (i), choose $\eps$ such that $0<\eps<\eps'$, $p_i(x)<1$ for all $0\le i\le r$ and $0\le x\le C$ and let $N=\lfloor C/\eps \rfloor$. This determines a chunky local deletion algorithm $\cal A$ with local rule $L_\phi$ and recolouring rule $c$.
Then define $p_{t,i} = \eps   p_i(\eps t-\eps)$ for all $1\le t\le N$. The difference equation in Theorem~\ref{chunkyranreg} gives precisely the solution of~\eqn{desystem0} by Euler's method  with step size $\eps$ (apart from a final possible partial step). Hence, applying that theorem, we conclude that   a.a.s.\  $Y_j(t)=ny_{j}(\eps t)+O(\eps n)$ ($1 \leq j \leq \r$) and $W(t)=ny_{\r+1}(\eps t) +O(\eps n)$, uniformly for $0 \leq t \leq N$. Part (i) follows upon taking $\eps$ sufficiently small.

For part (ii), given $\eps'>0$, let  $\mathcal{A}=\mathcal{A}(\eps'/2)$ defined in part (i).   Apply Theorem~\ref{large girth} to find $g=g(\eps')$ for which the vector $s_{r,\mathcal{A}}(t,n)$  of expected values of the components of   $\vs(t)=\big(Y_1(t),\ldots, Y_R(t),W(t)\big)$  is independent of the $n$-vertex $r$-regular input graph $G$ whenever the girth of $G$ is at least $g$.  By Theorem~\ref{conc rreg} the size of each component of $ \vs(t)$, when applied to a random $r$-regular graph, is a.a.s.\   $s_{r,\mathcal{A}}(t,n)+o(n)$.  Finally, part (i) implies that the size of  each component of $\vs(t)$, when applied to a random regular graph, is a.a.s\  within $\eps' n/2$ of the corresponding component $ny_{j}(\eps t)$, which implies the desired result.
\qed

This result summarises the general strategy for analysing a local deletion algorithm defined through appropriate selection, local and recolouring rules. Indeed, it tells us that, once we fix local and recolouring rules, and well-behaved selection probability functions, we may devise a chunky local deletion algorithm $\mathcal{A}$ whose behaviour for a random graph in $\mathcal{G}(\mathbf{d})$ a.a.s.\ follows the solutions of a system of ordinary differential equations, which can be derived explicitly from the local, recolouring and selection rules. Moreover, when applied to any fixed $n$-vertex $r$-regular graph with sufficiently large girth,  the expected value $c_r n$ of the output function $W$ at the end of the algorithm  is very close to the value achieved almost surely for random graphs. In particular, if the output  function counts the number of elements in a set  satisfying some particular property, this implies that every $r$-regular graph with sufficently large girth contains a set with this property whose size is bounded by $c_r n$.

%%%%%%%%%%%%%%%%%%%%%%%%%%%%%%%%%%%%%%%%%%%%%%%%%%%%%%%%%%%%%%%%%%%%%%%%%%
\section{Deprioritised algorithms and proof of Theorem~\ref{chunkyranreg}}
\lab{s:deprio} 
    
The aim of this section is to prove Theorem~\ref{chunkyranreg} and some further useful results. The main ingredient in our proofs is  the concept of  \emph{deprioritised algorithms}, that is, algorithms that resemble other algorithms in which random choices are made using a priority list of possible operations, but with the priority list replaced by an appropriate set of probabilities. For instance, the algorithms for independent and dominating sets introduced in Section~\ref{s:introresults} are \emph{prioritised} in that the selection step chooses only vertices with minimum degree in the survival graph. We say that this prioritisation is \emph{degree-governed}, as defined below. In particular, recall the procedure $P_{ind}$, which at every step chooses a vertex with minimum degree in the survival graph and deletes it along with its neighbours.  In order to deprioritise  this algorithm, one may estimate the probability $p_{t,j}$ that, at the $t$-th step of the algorithm, the minimum
  degree of $G_{t-1}$ is equal to some fixed degree $j$. Then, instead of just choosing a vertex with minimum degree, one first chooses the degree $i$ with probability $p_{t,i}$. A vertex $v$ of degree $i$ is then selected uniformly at random. In~\cite{deprio} it was shown that a certain class of algorithms yield almost the same results when deprioritised appropriately. 

This idea is easily generalized for algorithms on coloured graphs, as follows.
\begin{defn}[Type-governed or degree-governed deprioritised algorithm]\label{dd_deprio} 
A type-governed deprioritised algorithm is an iterative algorithm whose input is an $n$-vertex coloured graph $G$ with colours in a bounded set $\C$ and maximum degree $r$, and with the following properties.  Let the possible types with such colours and degrees be enumerated $1,\ldots,\r$. At each step~$k$,
\begin{description}
\item{(i)} a number $i \in \{1,\ldots,\r\}$ is chosen with probability $p_i(k,n)$ (where $\sum_ip_{i}(k,n)=1$);
\item{(ii)}  a vertex $v$ of type $i$ is chosen uniformly at random;
\item{(iii)} a (possibly randomised) rule is applied to the coloured graph, depending on $v$, which  colours some vertices with output colours, then possibly recolours some vertices with colours in $\C$, and finally deletes some vertices.
\end{description}
Such an algorithm is called degree-governed if there is only one colour and only one output colour, in which case colours can be ignored in the usual way.
\end{defn}
The random choices in such an algorithm are normally assumed to be independent of the previous steps of the algorithm, and of the graph $G$, with the exception of the usual dependencies within the application of the local rule.
%, however we will have occasion to modify this convention slightly
  To use the above definition in a precise manner, we would need to define what constitutes a `rule' in (iii), which will be done for certain deprioritised algorithms below. Note that if, in step (ii), there are no vertices of type $i$, then the algorithm cannot be followed; we say that it is {\em stuck}, and it is terminated.

Working towards the relationship between deprioritised algorithms and chunky local deletion algorithms, we immediately focus on a special class of deprioritised algorithms. Note that parts (i) and (ii) of the above definition provide a selection rule as required for native local deletion algorithms.
 
\begin{defn}[Amenable deprioritised algorithms] 
A deprioritised algorithm is said to be amenable in an interval $[0,M)$ if it is type-governed and satisfies both of the following.
\begin{enumerate}

\item   $p_{i}(k,n)=\tilde{p}_i\left(k/n\right)$ for a fixed set  of piecewise Lipschitz continuous functions $\tilde{p}_i : [0,M] \to [0,1]\ (1\le i \le \r)$, called {\em relative selection functions},  such that $\sum_{i=1}^\r \tilde{p}_i(x)=1$  for every $x \in [0,M]$.

\item the randomised rule in part (iii) of the definition of deprioritised algorithms consists of an exploration step followed by a recolouring step as in the definition of local deletion algorithms.
\end{enumerate}
\end{defn}
 
Note that an amenable deprioritised algorithm is a local deletion algorithm and is defined by its relative selection functions $\tilde{p}_i$, the local rule  and the  recolouring  function.  

We first consider the random graph $G\in \G(\dv)$.
\begin{thm}\lab{deprioranreg}
 Let $\dv$ have some probability distribution over the possible degree vectors with maximum entry $r$. Consider an amenable deprioritised  algorithm $\mathcal{D}$ with local rule $L_\phi$, recolouring function $c$, and relative selection functions $\tilde{p}_i$, and  consider an application of this algorithm to  an $n$-vertex random graph $G\in \G(\dv)$ with neutral colouring. Let $N_j$ be the (random) number of vertices of type $j$ in $G$ $(1\le j \le \r)$,  let $W$ be an output function  and set  $N_{\r+1}=0$.
For fixed  $\delta' >0$,  let $T$ be the first step of the algorithm for which there exists $j$, $1\le j\le \r$, such that $\tilde{p}_j(x)>0$ for some $x$ with $|x-T/n|<\delta'$ and the number of vertices of type $j$ in the survival graph is less than $\delta' n$.  
 Let $\yv$ be given by the solution of the differential equation system
\begin{eqnarray}\lab{desystem}
y_j'(x)&=& \sum_{i=1}^\r\tilde p_i(x) f_{j,i}(y_1,\ldots, y_\r)  \quad (1\le j \le \r+s)
\end{eqnarray}
with initial conditions $y_j(0)= N_j/n$ for $1\le j\le \r$ and $y_{\r+1}(0)=0$,
where $f_{j,i}$ is defined as in Lemma~\ref{l:trend}  and equation~\eqref{def_output_func}.
Then a.a.s.\
the  value  $W(t)$ of the output  function  after step $t$ is $ny_{\r+1}(t/n)+o( n) $, and the number $Y_j(t)$ of vertices of type $j$ in the survival graph $G_t$ is $ny_{j}(t/n)+o( n)$, uniformly for  $0\le t\le T$.
\end{thm}
\proof 
With a slight abuse of notation, let $G_t$ denote the survival graph (actually, pseudograph) obtained after $t$ steps of the algorithm applied to a random pairing in the pairing model with degree sequence $\dv$ and all vertices neutral. Similarly, define $Y_j(t)$ and $W(t)$ for this random pseudograph, and set $\Yv(t)=  (Y_1(t),\ldots ,Y_{\r}(t),W(t))$.
We show that the random vector $\Yv(t)$ is a.a.s.\ sharply concentrated near the solution of the differential equation system~\eqn{desystem}, using the differential equation method  as presented in~\cite{des} or~\cite{desurvey}.  The proof is by induction over the intervals on which the relative selection functions $\tilde{p}_i$ are Lipschitz continuous. Let $t_k$ denote the first value of $t$  for which $t/n$ lies in the $k$th such interval. The inductive hypothesis is that the claimed approximation holds for $t_{k-1}\le t\le \min\{T,t_k\}$. In particular, at $t=t_k$, each coordinate of the vector is a.a.s.\ within $o(n)$ of the differential equation solution. Initially, of course the two coincide precisely.

In the inductive step, we are only examining the part of the process from $x_k$ to $x_{k+1}$, where $x_k=t_k/n$ and $x_0=0$. Note that we can assume that the functions $p_j(t)=\tilde{p}_j(t/n)$ are Lipschitz continuous in the corresponding $t$-interval, since a discontinuity after $t_{k+1}-1$ does not affect the process before $t_{k+1}$. Also we may extend them to be continuous, say by making them constant, on intervals before and afterwards. Define $ D \subseteq \reals^{\r+2}$ to be the interior of the set of all $(x,y_1,\ldots, y_{\r+1})$ satisfying $x_k-\delta'< x< x_{k+1} +\delta'$, $-1< y_i <  C$ for all $i\le \r+1$ and $x$, and additionally $\delta'/2<y_i $  for each $ i\le \r $ such that $\tilde{p}_i(x')>0$ for some $x'$ with $|x'-x|<\delta'/2$. Here $C$ is a suitably large constant such that $Cn$ is a deterministic upper bound on the value of any of the variables $Y_i$. For $i\le \r$ any $C>1$ will suffice, whilst for the output  function, some such $C$ 
exists because this function is initially zero and its  variation in one application of the local rule is a constant that depends on the query graph that is active at this step, and hence is bounded.

Then $D$ is bounded, connected and open for $\delta'>0$ sufficiently small. To apply~\cite[Theorem 1]{des} requires verifying a few conditions, which we clarify below. We use the simple modification in~\cite[Theorem 6.1]{desurvey} that specifies that these conditions only need to be verified up to a stopping time $T'$, which we define as $\min\{T,t_{k+1}\}$. (The proof of this is quite easy, given~\cite[Theorem 1]{des}.)  See also~\cite[Theorem 5.1]{desurvey}  for some minor variations we make use of here.
Note that the vector $(t/n,\Yv/n)$ must remain inside $D$ and cannot come within distance $\delta/2$ of its boundary,  up until the time $T'$. We use the domain $D$ to give `elbow-room' within which the process behaves well, and the stopping time $T'$ for a convenient explicit stopping point beyond which we have no need of continuing the analysis. 

The conditions of the theorem, and the reasons that they hold, are  the following after some suitable trivial modifications. (For instance we shift  $0$ to $x_k$. Also, in the theorem,  $m$ is a suitable upper bound on the number of steps the process can take. Since each step deletes at least one vertex, one can choose $m>n( x_{k+1} +\delta')$ so it can be ignored.)  Here we define   $H_t$ to be the history   of the process, which consists of the sequence of coloured graphs $G_0,\ldots, G_t$. 

\smallskip

\noindent (a) There is a constant $C'$ such that for all  $t<T'$ and all $1\le j\le \r+1$,
$$|Y_{j} (t+1) - Y_j(t)|<C'$$ always. This holds since one application of the local rule alters the degrees of only a bounded number of vertices and changes the  value of the output function  by a bounded quantity.

\smallskip
\noindent
(b) For all $j$ and uniformly over all $t<T'$, 
$$\ex (Y_{j} (t+1) - Y_{j} (t)|H_t)=\tilde f_j(t/n,{\Y(t)}/n) +o(1) $$ 
for suitable functions $\tilde f_j$. As $t<T$, we have  $Y_j(t)\ge \delta' n$ for all $j$ such that $\tilde p_j(t/n)>0$. Hence, the algorithm does not terminate at this step.  By Lemma~\ref{l:survival}, the survival pairing is distributed u.a.r.\ given the types of its vertices. Hence we may apply Lemma~\ref{l:trend}(ii) to conclude that, for $j \leq \r$,  the expected change in $Y_j(t)$ resulting from the next step of the algorithm (which determines $\Yv(t)$) is equal to $\tilde f_j(t/n,\Yv(t)/n)+o(1)$ where $\tilde f_j(x,\mathbf{y}) = \sum_{i=1}^\r \tilde{p}_i(x) f_{j,i}(\mathbf{y})$.  Moreover, the expected change in the value of the output function at each step is determined by the query graphs explored at each step, and hence the expected change in $W(t)$ is given by 
$\tilde f_{\r+1}(x,\mathbf{y}) = \sum_{i=1}^\r \tilde{p}_i(x) f_{\r+1,i}(\mathbf{y})$ with $f_{\r+1,i}$ defined in~\eqref{def_output_func}. 
Note that these computations are valid even conditional upon the history of the process up to this step. 
\smallskip

\noindent
(c) For each $j$ the function $f_j$ is Lipschitz continuous  on $D$.   This follows from Lemma~\ref{l:trend}, since amenability implies that the relative selection functions $\tilde{p}_i$ are Lipschitz  continuous, and also the number $M$ of Lemma~\ref{l:trend} is at least $Y_r\ge  \delta'n/2$ on $D$. 
Moreover, $D$ contains
the closure of
 $$ \{(0,y_1, \ldots , y_{\r+1}):
\pr(Y_j(0)=y_j n, 0\le j \le \r+1)\ne0 \mbox{ for some}\: n\}.$$

\smallskip

The conclusion of the theorem is as follows.

\smallskip
\noindent
(i)  The system of
differential equations
$$  y_j'(x)  = \tilde f_j(x,y_1, \ldots , y_{\r+1}),
\qquad j = 1, \ldots, \r+1$$ 
  has a unique solution in $D$ for $y_j: \reals \to \reals$, which we denote by $\tilde y_j$,
passing through  
$$\tilde y_j(x_k) =Y_j(t_k)/n \quad \mbox{$1\le j\le \r+1$},$$
and which extends to points arbitrarily close (in Euclidean distance say) to the boundary of $D$. (Here $Y_j(t_k)/n$ is deterministic if $k=0$.)
 
\smallskip
\noindent
(ii)  Asymptotically almost surely 
$$Y_j(t) = n \tilde y_j(t/n) + o(n)$$
uniformly for $t_k\le t \le \min \{\sigma n,T'\}$ and for each $j$, where
$\sigma =\sigma(n)$ is the supremum of those $x$ to which the solution can be
extended.   

Suppose that $Y_j(t)$ has not yet reached a point such that the condition in the definition of $T'$ holds.    
Immediately, by the definition of $T'$, if $  Y_j(t) $ is sufficiently close to $ n \tilde y_j(t/n)$, the point $\tilde y_j(t/n)$ has distance at least  $\delta'/4$ (say) from the boundary of $D$. Thus, the approximation in (b) a.a.s.\ holds  for $t_k\le t \le   T' $. 
 
By the inductive hypothesis, we have a.a.s.\ 
$$Y_j(t_k)/n=  y_j(x_k) +o(1)  \quad \mbox{$1\le j\le \r+1$},$$
and hence  $\tilde y_j(x_k) =  y_j(x_k) +o(n)$ a.a.s. Since the derivatives in the differential equation system are Lipschitz, the standard property of solutions of differential equations   implies that $\tilde y_j(x) =  y_j(x) +o(n)$ uniformly for all $x_k\le x < \min\{\sigma, x_{k+1}\}$. Thus, from (b) above and the ensuing conclusions, a.a.s.\
$$Y_j(t) = n   y_j(t/n) + o(n)$$
uniformly for $t_k\le t \le   \min\{T',t_{k+1}\} = \min\{T,t_{k+1}\} $ and for each $1\le j \le \r+1$. The differential equation system in (i) becomes that of~\eqn{desystem}. This proves the inductive hypothesis.
\qed
\smallskip

\noindent
{\bf Note.}  We may permit the following variation of type-governed deprioritised algorithms: according to part (ii) of the definition of the algorithm, the next vertex of type $i$ is chosen u.a.r. If the algorithm is local, the next vertex of type $i$ may be chosen by any rule whatsoever that depends only on the outcomes of the previous applications of the local rule. (Similarly, somewhat more restrictively, it is valid to permit the rule to depend only on the set of edges incident with vertices that have so far been deleted by the algorithm.) The conclusion of Theorem~\ref{deprioranreg} holds in this case also. This is because then the choice of the next vertex depends only on the history $H_t$ and the survival graph $G_t$ occurs u.a.r.\ among the possible pseudographs generated by the pairing,  given the types of its vertices, and Lemma~\ref{l:trend} permits any rule for choosing the next vertex, provided its type $i$ is selected with the correct probability $p_i$.
\smallskip

\noindent
{\bf Proof of Theorem~\ref{chunkyranreg}\ }
Let $\mathcal{A}$ denote the given chunky local deletion algorithm. As discussed before, it is enough to apply the algorithm to the random pseudograph with degree sequence $\dv$ obtained from the pairing model. We will define a deprioritised algorithm whose behaviour is similar to that of $\mathcal{A}$, within errors $O(\eps   n)$ in the values of the cardinalities of the relevant sets.  
Define $Y_j$ ($1\le j \le \r+1$) as in the proof of Theorem~\ref{deprioranreg}. 
By induction on $t$, $0\le t\le N-1$,   we will show 
for each $j$   that  
\bel{indclaim}
Y_j(t)=n z_j(t)+   O(\eps^2tn)  \quad\mbox{a.a.s.}
\ee
 (where, of course, the constant in $O()$ does not depend on $t$).
We may assume that this is true for some $t-1\ge 0$, and   consider step $t$.   For $1\le i\le \r$, let $S_i$ denote the set of  vertices of type $i$ in the survival graph $G_{t-1}$ chosen by $\mathcal{A} $ in this step, and put $Z_i=|S_i|$. Then $Z_i$ has a binomial distribution with expected value $p_{t,i}  Y_i(t-1)\le \eps n $, and so by Chernoff's bound, with probability $e^{-\Omega(n)}$,
\bel{zi}
  |Z_i- p_{t,i}  Y_i(t-1)|\le  \eps^2 n,\qquad |Z_i|< 2\eps n;
\ee 
and we may assume these henceforth.

Now consider a process $\widetilde {\cal A}$ which consists of repeatedly applying the following step, beginning with   $G_{t-1}$. 

At each step, let $S_i'$ denote the set of vertices of $S_i$ that remain, and still have the degree $d(i)$ associated with type $i$, in the current survival graph. Let $i$ be minimum such that  $S_i'$ is nonempty, select a  vertex $v$ u.a.r.\  from $S_i'$, apply the local rule to $v$, and update the output sets and the survival graph according to the query graph obtained. We say that this step {\em processes} $v$.

The process $\widetilde {\cal A}$ naturally breaks up into $\r$ phases, where phase $i$ consists of processing vertices in $S_i'$ ($1\le i \le \r$), normally until none of these are left. Vertices of $S_i$ that remain in the survival graph after phase $i$ must then have lost at least one neighbour during the process $\widetilde {\cal A}$. However, there is also an abnormal variation. Conditional upon the history of $\widetilde {\cal A}$ to the end of phase $i-1$, we may stop phase $i$ prematurely, before it even begins, if fewer than $\eps^2 n$ vertices  lie in $S_i'$ at the end of the previous phase. In this case we say the phase is {\em skipped} and the process proceeds to the next phase.
 
Phase $i$ of the  process $\widetilde {\cal A}$ is essentially a deprioritised algorithm, the only difference being that it contains the following two variations. Firstly, the choice of the next vertex of type $i$ is restricted to vertices in the set $S_i'$. By the note after Theorem~\ref{deprioranreg}, such restriction does not alter the applicability of that theorem.  Secondly, there is a stopping rule (stop prematurely when  no vertices remaining in $S_i$ are still of degree $d(i)$, or if $S_i'$ is too small at the start). The conclusion of the theorem still applies at this stopping time. (To see why, see~\cite[Section 4.2]{desurvey} and the proof of Theorem~6.1 of that paper.)

Let $\widetilde Y_j(\xi)$ denote the number of vertices of type $j$  (or,  in the case $j=\r+1$,  the  value of the output function $W$) after $\xi$ steps of the process $\widetilde {\cal A}$.
We may now apply Theorem~\ref{deprioranreg} with $\delta'=\eps^2$ (and the note as discussed above), and with 
\bel{pi}
\tilde{p}_j = \left\{
\begin{array}{rl}
1&\mbox{if $j=i$}\\
0& \mbox{otherwise},
\end{array}
\right.
\ee
to conclude that phase $i$ of $\widetilde {\cal A}$  a.a.s.\ finishes with  
\bel{enderror} 
\widetilde  Y_j(\xi_i)/n=  \tilde y_j(\xi_i/n) +o(n)  \quad \mbox{$1\le j\le \r+1$}, 
\ee
where $\xi_i$ denotes the last step of phase $i$ (which is random and we will examine   shortly), $\xi_0=0$,
and 
the variables $ \tilde y_j$ satisfy the equations~\eqn{desystem} with initial conditions 
$$
\tilde y_j(\xi_{i-1}/n)=   \widetilde  Y_j(\xi_{i-1})/n.
$$
Of course in the case $i=0$,  $\widetilde  Y_j(\xi_{i-1})$ is defined to be $Y_j(t-1)$, the random variable inputted to  the process ${\cal A}$, which by induction on $t$ satisfies $Y_j(t-1)=n z_j(t-1)+   O(\eps^2tn)$. 
Recall that $p_i=1$ in phase $i$, the functions $f_{j,i}$ are Lipschitz, and the length $\xi_i-\xi_{i-1}$ of phase $i$ is $O(\eps n)$   by~\eqn{zi}.  Hence 
\bean
\tilde y_j(\xi_i/n)&=&    \frac{\widetilde Y_j(\xi_{i-1})}{n}  +\frac{\xi_i-\xi_{i-1}}{n} f_{j,i}\big(\tilde \yv (\xi_{i-1}/n)\big) + O(\eps^2)\\
&=&   \frac{\widetilde Y_j(\xi_{i-1})}{n} +\frac{\xi_i-\xi_{i-1}}{n} f_{j,i}\big(\tilde \yv (0)\big) + O(\eps^2)  
\eean
where $\tilde \yv=(\tilde y_1,\ldots, \tilde y_{\r+1})$.
Thus, by induction using~\eqn{enderror}, the final values at the end of the process $\widetilde {\cal A}$ are a.a.s.\ given by
\bel{zchanges}
\tilde y_j(\xi_\r/n) = \tilde y_j(0) +\sum_{i=0}^\r \frac{\xi_i-\xi_{i-1}}{n} f_{j,i}\big(\tilde \yv (0)\big) + O(\eps^2) 
\ee
for $1\le j\le \r+1$.

Before proceeding, we need to bound the effect of the   vertices being deleted during  the process   $\widetilde {\cal A}$.  We call a vertex {\em susceptible} if it is a member of  $S=\bigcup_{i=1}^{\r}S_i$  and has distance at most $2D$ from some other vertex in $S$.  {In the terminology of Section~\ref{s:largegirth}, a vertex of $S$ is susceptible if it forms a pre-clash with some other  vertex of $S$. Let $A_t$ denote the set of susceptible vertices.

We next  claim that, for some constant $C_0>0$ depending only on $r$ and on the number of colours, conditional on the graph $G_{t-1}$,  $\pr(|A_t|>C_0\eps^2n ) = e^{-\Omega(n)}$. To establish this claim first note that by~\eqn{zi}, the probability that a given vertex is contained in $S$ and is furthermore susceptible is $O(\eps)$ (noting that $D$ is fixed, determined by the local rule, and the maximum degree of $G_{t-1}$ is at most the fixed number $r$). Hence
$\ex |A_t| = O(\eps^2 n)$. 

On the other hand, sharp concentration of $|A_t|$ can easily be proved as follows. Note that a simple switching (replacing two pairs on four points in the pairing by any two other pairs on the same points) can only change $|A_t|$ by an amount bounded by a function of $D$ and $r$. For the case of $r$-regular graphs,~\cite[Theorem 2.19]{models} immediately implies that $\pr (|A_t|>\ex   |A_t| +\eps^2 n)= e^{-\Omega (n)}$. In this case we have $G_{t-1}$, which is not regular, but it is easily checked that the proof  of that theorem still applies, and the constants in the bounds are independent of the degree sequence of $G_{t-1}$ as long as the maximum degree is $r$. 
Hence, we have shown the claim holds, and we may proceed assuming $ |A_t|\le C_0\eps^2n$.

We now turn to consideration of the lengths of the phases, i.e.\ $\xi_i-\xi_{i-1}$.  
 Note that a non-susceptible vertex in $S$ of type $i$ will definitely retain degree $\d(i)$ in the survival graph until at least the start of phase $i$. Hence, if the phase is not skipped, the number of vertices of degree $d(i)$ in $S_i$ available for processing must be between $Z_i$ and $Z_i- A_t$, and the phase cannot finish until there are at most $\delta' n = \eps^2n$ of them left. 
It follows that 
\bel{length}
\xi_i-\xi_{i-1} = Z_i+O(\eps^2 n).
\ee                                
On the other hand, if the phase is skipped then $\xi_i-\xi_{i-1}= 0$  and $Z_i- A_t<\eps^2n$, so~\eqn{length} holds in all cases. Thus, by~\eqn{zi} and~\eqn{zchanges}, we may assume that
$$
\tilde y_j(\xi_\r/n) = \tilde y_j(0) + \sum_{i=1}^\r p_{t,i}\frac{Y_i(t-1)}{n} f_{j,i}\big(\tilde \yv (0)\big) + O(\eps^2).
$$
Recall that $\tilde \yv(0) = \Yv(t-1)/n$, so by the inductive claim~\eqn{indclaim} for $t-1$, we may assume
$
\tilde y_j(0)=  z_j(t-1)+   O_1(\eps^2(t-1) )  
$
for each $j$. We use $O_1()$ to distinguish the constant here, which is bounding the inductive error, from the other error bound constants. By the Lipschitz property of each $f_{j,i}$, it follows that 
 $f_{j,i}\big(\tilde \yv (0)\big)= f_{j,i}\big(  \zv (t-1)\big) +   O(\eps^2(t-1))
 = f_{j,i}\big(  \zv (t-1)\big) +   O(\eps )
$ since $t\le N=O(1/\eps)$ by hypothesis in the theorem.
 Hence, using $p_{t,i}\le \eps$ and $Y_i\le n$, together with~\eqn{enderror} for $i=\r$, and then~\eqn{indclaim} again for the second step, we obtain  
\bean
\tilde Y_j(\xi_\r)/n &= &  z_j(t-1)   + O(\eps^2(t-1))+ \sum_{i=1}^\r p_{t,i}\frac{Y_i(t-1)}{n} \big(f_{j,i}\big(\zv (t-1)\big) + O(\eps)\big)\\
&= &   z_j(t-1)   + O_1(\eps^2(t-1))+ \sum_{i=1}^\r p_{t,i}z_i(t-1) (f_{j,i}\big(\zv (t-1)\big) + O(\eps^2)\\
&=&    z_j(t) + O_1(\eps^2(t-1))   + O (\eps^2  ).
\eean

To establish~\eqn{indclaim} inductively, it only remains to show that $\tilde Y_j(\xi_\r)= Y_j(t) +  O (\eps^2  n)$ a.a.s.
This fact is easy to see, because, in view of~\eqn{length}, the processes $\widetilde {\cal A}$ and ${\cal A}$ process the same vertices except for a set of  $O(\eps^2 n)$ vertices, each of which has a bounded effect on the value of $Y_j$ or $\tilde Y_j$. We now have~\eqn{indclaim}, so it is true for $t=N$. Recalling that $\eps t\le \eps n=O(1)$, we deduce that the  value of the output function  at the end of the algorithm is equal to $nz_{\r+1}(N)+O(\eps n)$, and the conclusion of the theorem follows. 
\qed

Before proceeding, we note that the conclusion of Theorem~\ref{deprioranreg} can be restated as follows.\smallskip

\noindent
{\bf Alternative Conclusion for Theorem~\ref{deprioranreg}.} \ Let $\delta>0$ and let $x_0$ be the infimum of all $x\ge 0$ for which there exists $j$ such that $\tilde{p}_j(x)>0$  and $y_j(x) < \delta$. Then a.a.s.\ $Y_j(t)=ny_j(t/n)+o(n)$ and $W(t)=ny_{\r+1}(t/n)+o(n)$ for all $t \leq \lfloor (x_0-\delta)n\rfloor$. \smallskip

To see why this follows, define $\delta'=\delta/2$ and note that the conclusion of Theorem~\ref{deprioranreg} implies that for $0 \leq t \leq T$ the value of $Y_j(t)$ differs from $ny_{j}(t/n)$ by at most $o(n)$, so that step $T$ a.a.s.\ does not occur before $n(x_0-\delta)$. Note also that the proof of the theorem itself gives the alternative conclusion quite directly.
 
\begin{thm}\lab{chunkydeprio}
Assume the hypotheses of Theorem~\ref{deprioranreg}, let   $\delta>0$ and define $x_0$ as in the alternative conclusion given above. Then there exists a number $\eps$, $0<\eps<\delta$, and a chunky local deletion algorithm $\mathcal{A}$ with the same local  and recolouring rules as  $\cal D$,  and granularity at most $\delta$, such that,  when both algorithms are applied to $\Gnr$, a.a.s.\ $\left|\left|\Yv_{\mathcal{A}}(t)-\Yv_{\mathcal{D}}(\lfloor\eps n t\rfloor) \right|\right|\leq  \delta n$ for all $0 \leq t \leq \lfloor (x_0-\delta)/\eps\rfloor$, where $\Yv_{\mathcal{A}}(t),\Yv_{\mathcal{D}}(t) \in \mathbb{R}^{\r+1}$ denote the vectors whose components $j \leq \r$ and $\r+1$ are the number of vertices of type $j$ in the survival graph, and the value of the output function, respectively, after $t$ steps of the algorithm. 
\end{thm}

\proof 
Put $\delta'=\delta/2$ and let ${\bf y} = (y_1,\ldots , y_{\r+1})$ denote the solution of the differential equation system~\eqn{desystem} in the domain $D$ used in the proof of Theorem~\ref{deprioranreg} up until the endpoint is reached.

Define $p_i(x)=\tilde p_i(x)/y_i(x)$, where $y_i$ is precisely the solution function just determined, for all $0\le x\le x_0$. Then, on this interval, the differential equation system~\eqn{desystem} is identical to~\eqn{desystem0}, with the same initial conditions. Hence, this theorem follows from Theorem~\ref{chunky des}(i) applied with $\eps'=\delta$ and $C=x_0-\delta$, and the alternative conclusion of Theorem~\ref{deprioranreg} with $t$ replaced by $\lfloor \eps n t\rfloor$. Note that the fact that the functions $\tilde p_i(k,n)$ are bounded implies that the floor function has no significant effect. \qed

Theorem~\ref{chunkydeprio} has an immediate consequence for regular graphs with sufficiently large girth which will be useful in the next section. 

\begin{thm}\lab{chunkydeprio large girth}
Fix integers $r,D \geq 1$. Consider  an amenable deprioritised  algorithm $\mathcal{D}$ with  local rule $L_\phi$ of depth $D$, recolouring function $c$ and relative selection functions $\tilde{p}_i$.  Let $W$ be an output function.  Let $\yv$ be given by the solution of the differential equation system~\eqn{desystem} with initial conditions $y_1(0)=1$ and $y_j(0)=0$ for $j>1$,
where $f_{j,i}$ is defined as in Lemma~\ref{l:trend}  and~\eqref{def_output_func}. Let $\delta>0$ and define $x_0=x_0(\delta)$ as the infimum of all $x\ge 0$ for which there exists $j$ such that $\tilde{p}_j(x)>0$  and $y_j(x) < \delta$.  Then there exists some $\epsilon>0$ with $\eps<\delta$, a constant $g$, and a chunky local deletion algorithm $\mathcal{A}$ with local rule $L_\phi$, recolouring function $c$ and granularity at most $\delta$ such that, for every $n$-vertex $r$-regular graph $G$ with girth at least $g$,
$$\left|\left|s_{r,\mathcal{A}}(G,t) - n\yv(\eps t)\right|\right|_\infty<\delta n,$$
 for $0\leq t \leq \lfloor (x_0-\delta)/\eps\rfloor$, where $s_{r,\mathcal{A}}(G,t)$ is the vector of expected values for algorithm  $\mathcal{A}$ defined in Theorem~\ref{large girth}.
\end{thm}

\proof Let $\delta>0$ and $x_0$ be as stated. Theorem~\ref{deprioranreg} implies that the performance of $\mathcal{D}$ applied to a random $r$-regular graph is a.a.s.\ determined by the solution of~\eqref{desystem} up until the quantity $T$ defined there. Note that the definition of $x_0$ ensures that $T \geq (x_0-\delta)n$ a.a.s.  By Theorem~\ref{chunkydeprio} (applied with $\delta$ replaced by $\delta/2$), there is a chunky local deletion algorithm $\mathcal{A}'$ whose performance in $\Gnr$ is tracked by a vector $\Y$ that differs from its counterpart in an application of $\mathcal{D}$ by at most $(\delta/2)n$. The `steps' of that algorithm are of size $\eps<\delta$. The algorithm $\mathcal{A'}$ is obtained in the proof of Theorem~\ref{chunkydeprio} by appealing to Theorem~\ref{chunky des}(i).  The girth $g$ and the algorithm $\mathcal{A}$  may be obtained from the further information given in Theorem~\ref{chunky des}(ii), using $\eps'=\delta/2$. \qed

For the performance of the chunky algorithm on a graph of large girth, the bound on the expected number of susceptible vertices was enough in~\cite{LW} to be able to say that the effect of these vertices is negligible (since there are $O(1/\eps)$ chunks). However in that case the analysis required some `independence lemmas' showing that certain appropriate events were independent in the branches of the trees induced by the vertices near a given vertex, which we will not need in the current approach. In addition to the obvious advantage of not needing to prove this independence, our strategy extends to algorithms whose actions depend on what happens in two branches simultaneously, for which no such independence can be expected to hold. The (small) price we had to pay for avoiding independence lemmas is that we needed to analyse the effect of the vertices that have small distance apart on the variables $Y_i$, and this requires sharp concentration arguments.

%%%%%%%%%%%%%%%%%%%%%%%%%%%%%%%%%%%%%%%%%%%%%%%%%%%%%%%%%%%%%%%%%%%%%%%%%%%%%%%%%%%%%%
\section{Applications}\lab{sec:applications}

There is a long list of algorithms which fit in the definition of local deletion algorithms that have already been analysed to obtain properties of random regular graphs, either with direct arguments in the spirit of~\cite{des,desurvey}, or through a general-purpose theorem using deprioritised algorithms~\cite[Theorem 1]{deprio}.  In this section, we consider the latter theorem, combining its results with those of the present paper, in particular Lemma~\ref{l:trend} and Theorem~\ref{chunkydeprio large girth}. We show that the analysis of deprioritised algorithms acting on random regular graphs in~\cite{deprio} may be carried over  to the large girth setting, at least for a general class of algorithms including basically all those to which it has been applied. We consider various concrete instances, showing that several bounds that are known to hold asymptotically for random regular graphs, by virtue of their use of~\cite[Theorem 1]{deprio}, immediately imply deterministic bounds for graphs whose girth is sufficiently large. Other  applications of our general results, including new instances of applicability of~\cite[Theorem 1]{deprio}, are found in Section~\ref{s:more}.

In all the applications we will use Lemma~\ref{l:clashes} to show that the output set of the chunky algorithm is (for the properties being considered) within $\eps n$ size of the desired set. Recall that this error term arises from possible clashes of query graphs in an application of a local deletion algorithm.

\subsection{Extension of a general theorem on deprioritised algorithms}\lab{sec:ext}
 
In this subsection we use the results from the previous sections to build on  the results from~\cite{deprio}.} In that paper, there are shown to exist a general class of deprioritised algorithms whose  performances on random regular graphs are determined (a.a.s.)\ by their basic operations. For such an algorithm, if the input graph is $r$-regular, each basic operation is one of several types, called ${\rm Op}_i$ for $1\le i\le r$, where ${\rm Op}_i$ consists of selecting a vertex $v$ of degree $i$ in the survival graph u.a.r., and then applying  a specified sequence of randomised tasks (including deletion of $v$). Actually, these operations are defined in terms of pairings, but we may extend any operations for graphs   to pairings in the obvious way, and as explained in the comment just after Lemma~\ref{l:trend} the way they act on query graphs with loops or multiple edges is irrelevant. Provided that the operations satisfy certain conditions, the conclusion of~\cite[Theorem~1]{deprio} is that there is a randomised algorithm using these operations, such that a.a.s.\ at some point in the algorithm the output set size is  asymptotic to $\rho n$ and, for $i=1,\ldots,r$, the number of vertices of degree $i$ is   asymptotic to $\rho_i$. Here $\rho$ and $\rho_i$ are specific constants that can be computed by solving a system differential equations.
 See the end of this subsection for a statement of a special case  of~\cite[Theorem 1]{deprio}. Its hypotheses are fairly technical, and are mainly used to ensure that a system of differential equations associated with the random process has a solution that can be extended arbitrarily close to the boundary of the region in which the variables are meaningful (for instance, in the region where variables representing probabilities lie in $[0,1]$). We omit the precise conditions here in the general case, since our present concern is to make further deductions in some cases where the conditions have already been verified by existing papers in the literature, so we rely on these references for verification of the hypotheses. 
 The general result relates to the pairing model $\mathcal{P}_{n,r}$ discussed in Section~\ref{s:explicit}. As explained there, this implies corresponding results for random graphs. Note that the result does not mention the number of vertices of degree 0 in the final graph.  In general, these are irrelevant,  because the operations can normally be designed such that all vertices of degree 0 are immediately deleted in any step in which they are formed, so their number in the survival graph is always 0.

In this section, we restrict ourselves to native local deletion algorithms (so, in particular, the type of a vertex is just its degree), and consider a local rule $L$ and recolouring function $c$ from such an algorithm. We can specify the sequence of random tasks for ${\rm Op}_i$ (after selection of $v$) to consist  of two steps. First,   obtain   a copy $H$ of a query graph via $L$. Next, add $c(H)$ to the output set $\U$, and delete  the vertices in $H$.  Thus,  ${\rm Op}_i$ consists precisely of parts (ii) and (iii) in the definition of a certain kind of degree-governed deprioritised algorithm.  We call such an operation ${\rm Op}_i$ a {\em degree-governed query} operation, and say that it {\em arises} from $L$ and $c$. Note that we continue to assume that every local rule is defined on all coloured pseudographs, though only its action on graphs will matter, by the note after Lemma~\ref{l:trend}.
 
Given appropriate operations ${\rm Op}_i$ for $1\le i\le r$,  let  $\rho({\rm Op}_1,\ldots, {\rm Op}_r)$ and $\rho_i({\rm Op}_1,\ldots, {\rm Op}_r)$ denote the constants referred to above which are determined in the conclusion of~\cite[Theorem 1]{deprio}.  

 \begin{thm}\lab{t:fromDeprio}
Fix $r\ge 1$, a local rule $L$ and an recolouring function $c$. Suppose that ${\rm Op}_i$, $1\le i\le r$, are degree-governed query operations arising from $L$ and  $c$, and that the pairing versions of these operations satisfy the assumptions of~\cite[Theorem 1]{deprio}. Let $\rho,\rho_1,\ldots,\rho_r$ be the constants defined above. Then, for all $\delta>0$, there exists a  chunky native local deletion algorithm $\A$ with the very same local rule $L$ and recolouring function $c$, and granularity at most $\delta$, such that when $\A$ is applied to any $n$-vertex $r$-regular graph $G$  of sufficiently large girth,  there is a step $T$ in which the expected size of the set $\U_1$ of vertices with output colour 1 is within $\delta n$ of $\rho n$, and where the expected number of vertices of degree $i$ in the survival graph $G_T$ is within $\delta n$ of $\rho_i n$ for $1 \leq i \leq r$. Moreover, the number of vertices with colour $\propto'$ or $\propto$ at step $T$ is at most $\delta n$.
\end{thm}
\proof The conclusion of~\cite[Theorem 1]{deprio} asserts the existence of a randomised algorithm, and  it is straightforward to check that the proof of~\cite[Theorems 1 and 2]{deprio} consists of constructing an amenable deprioritised algorithm whose steps (ii) and (iii) consist of ${\rm Op}_i $ for vertices of degree $i$. More precisely, it is a parametrised family of algorithms. (Note in particular that the definition of $p_i$ in~\cite[Equation (4.2)]{deprio} depends only on fixed functions of $x=t/n$ where $t$ is the step number of the algorithm.) The first  assumption  of~\cite[Theorem 1]{deprio} is that each operation ${\rm Op}_i$ satisfies~\cite[Eq.\ (2.2)]{deprio} when (in the notation of the present paper) the number, $Y_r$, of vertices of degree $r$, is at least $\eps n$ (for some predetermined $\eps>0$). In terms of the present paper, this requires that for some fixed functions
$f_{k,i}\left(x,\yv\right)=f_{k,i}\left(x,y_1,\ldots ,y_{r+1}\right)$, 
the expected increase in $Y_k$ in one step of the algorithm involving ${\rm Op}_i$ is  $f_{k,i}(t/n,Y_1(t)/n, \ldots , Y_{r+1}(t)/n) +o(1)$ for $i=1,\ldots, r$, $k=1,\ldots, r+1$. Here $Y_{r+1}$ denotes the size of the output set, and the convergence in $o(1)$ must be uniform over all states of the algorithm with $Y_r\ge \eps n$.

Since ${\rm Op}_i$ arises from $L$ and $c$,  we may apply  Lemma~\ref{l:trend}(b) and (c), and we deduce that these   functions $f_{k,i}$ exist and must be exactly the  same as the functions appearing in Lemma~\ref{l:trend}, at least where all arguments are nonnegative  and $Y_r>\eps n$.  (Note that $j$ in~\cite{deprio} corresponds to $r+1$ in the present paper).   Theorem~\ref{chunkydeprio large girth} refers to the same functions $f_{k,i}$. 

We pause to summarise the proof of~\cite[Theorem~2]{deprio}. Probability functions  
$\tilde{p}_i$ are constructed to complete the definition of a deprioritised algorithm that uses the operations ${\rm Op}_i$ mentioned above. 
The probabilities are strategically altered slightly from some that are naturally associated with the original native local deletion algorithm which always selects a random minimum degree vertex. 
The alteration is introduced in order to keep the required variables strictly positive, and it uses an arbitrarily small parameter $\eps_1$. With these probabilities, the solution of the differential equations given by~\eqn{desystem} (c.f.~\cite[Equation (3.5)]{deprio}) is named  ${\bf \tilde y^{(1)}}=(\tilde y_1^{(1)},\ldots,\tilde y_{r+1}^{(1)})$ to differentiate it from the solution ${\bf \tilde y}$ which arises when the natural probabilities are used. The initial condition  in both cases is $(0,\ldots,0,1,0)$, so in particular $\tilde y_{r}^{(1)}(0)=1$. Keeping $\eps$ fixed, these two vector solution functions are shown to be arbitrarily close to each other up to a point $x_m$, when  $\eps_1$ is arbitrarily close to 0 
(see~\cite[(4.25)]{deprio}). By definition of $x_m$, the solution ${\bf \tilde y}$ is such that $(x,{\bf \tilde y})$  lies inside a domain ${\cal D}_\eps$ for $0\le x\le x_m$. As part of the definition of this domain, $\tilde y_{r }\ge \eps$ at all points inside it. The proof also makes use of other quantities $x_1,\ldots, x_{m-1}$, being points at which the derivatives associated with the non-deprioritised algorithm are discontinuous.

Given $\delta>0$, to apply Theorem~\ref{chunkydeprio large girth}  we need to determine $x_0$, which is the infimum of all $x\ge 0$ for which there exists $j$ such that $\tilde{p}_j(x)>0$  and $\tilde y^{(1)}_j(x) < \delta$.  We claim (shortly to be justified) that for $\delta$ sufficiently small, the functions $\tilde{p}_j$ in the proof of Theorem~2 of~\cite{deprio} are defined in such a way that guarantees that $x_0$ is at least as large as the quantity $x_m^{(1)}$ defined in~\cite{deprio}. As stated at the end of that proof, $x_m^{(1)}$ can be taken arbitrarily close to $x_m$. (Warning: there is a different quantity $x_0$ in that paper, but it equals 0 and does not play a significant role.)  Actually, this is a necessary feature of that proof, by its very nature: if $\tilde{p}_j(x)>0$  and $\tilde y^{(1)}_j(x)$ became arbitrarily small, there would be no guarantee that there are any vertices available for selection and application of ${\rm Op}_j$. To see why the claim is true, first
observe that in Part 2 of the proof of~\cite[Theorem  2]{deprio}, the  probability vector in the deprioritised algorithm is defined so that $\tilde p_r=1$ and $\tilde p_i=0$ for $i<r$ on an interval $[0,\eps_1]$ where $\eps_1>0$ is one of the parameters of the family of algorithms. This interval can be regarded as a ``burn-in'' period for the deprioritised algorithm (in~\cite{deprio} it is called the preprocessing subphase): vertices of degree $r$ are selected each time in the selection step, and the fact that they are abundant (for sufficiently small $\eps_1$) guarantees that the steps of the algorithm can (a.a.s.) be carried out as prescribed. (The burn-in period is the basic reason that ${\bf \tilde y^{(1)}}\ne {\bf \tilde y}$.)   On $[0,x_1]$, the function $\tilde y_r^{(1)}$ is bounded away from 0, and on $[\eps_1,x_m]$,  the functions $\tilde y_j^{(1)}$ ($1\le j\le r$) are all bounded away from 0. These claims follow from some of the observations in the proof in~\cite{deprio}, as we
  describe in more detail below. This then implies that $x_0\ge x_m^{(1)}$ as desired.

The claim about $\tilde y_r^{(1)}$ follows directly from the facts mentioned above, which imply that $ {  \tilde y_r}(x)\ge \eps )$   for $0\le x\le x_m$, and that the difference between ${\bf \tilde y^{(1)}}$ and ${\bf \tilde y}$ can be made arbitrarily small by choice of the parameters $\eps_1,\ldots$.
 It remains to check that the functions $y_i$ remain positive on $[\eps_1,x_m]$ ($1\le i\le r-1$). Being continuous, they will then be bounded below by a positive constant.  First, in the interval $[\eps_1,x_1]$, called phase 1, this is shown  in~\cite[(4.16), (4.18), (4.19)]{deprio}. The general argument for the later phases, covering the interval $[x_1,x_m^{(1)}]$, is a minor variation of the argument for the first phase. In summary, the main features work like this. In phase $k$, there is a preprocessing subphase which gets all $\tilde{y}_i$  strictly positive. Then, $\tilde y_{d-k}^{(1)}$ must remain positive by definition of the phases (i.e.\ the number $x_k$) except perhaps in the final phase, where it can reach 0 at $x_m$. The definition of $x_m^{(1)}<x_m$ excludes that case. By definition of the probabilities, $\tilde y_{d-k}^{(1)}$ has zero derivative in phase $k$ and hence remains positive. For all other $1\le i\le r-1$,   $\tilde y_{i}^{(1)}$ cannot reach 0 because by hypothesis (B) in~\cite{deprio}, its derivative is at least $-C_2\tilde y_{i}^{(1)}$ for some constant $C_2$ (so it is bounded below by an exponentially decaying function). In the case $i=m$, $\tilde y_{i}^{(1)}$ is prevented from reaching 0 by the definition of the domain $\D_{\eps,M}$ of~\cite[(3.3)]{deprio}, which requires $y_d>\eps$, and the fact that $x_m$ cannot exceed the first point that the solution leaves the domain $\D_{\eps,M}$.

Given this,  we may apply Theorem~\ref{chunkydeprio large girth} to deduce that, for all $\delta'>0$, there is $\eps>0$, a constant $g$ and a chunky native local deletion algorithm $\mathcal{A}$ with local rule $L_\phi$ and granularity at most $\delta'$ such that
$$\left|\left|s_{r,\mathcal{A}}(G,t) - n\tilde{y}_{i}^{(1)}(\eps t)\right|\right|_\infty<\delta' n,$$
for $0\leq t \leq \lfloor (x_0-\delta')/\eps\rfloor$, where $s_{r,\mathcal{A}}(G,t)$ denotes the vector whose components $j \leq r$ and $r+1$ are the expected number of vertices of type $j$ in the survival graph and the expected cardinality of the output set of colour 1 when $\mathcal{A}$ is applied for $t$ steps to an $r$-regular input graph $G$ with girth at least $g$. 

 Since $x_m^{(1)}\le x_0$, we can replace $x_0$ by $x_m^{(1)}$ in this statement. By the fact mentioned above that ${\bf \tilde y^{(1)}}$ and ${\bf \tilde y}$ can be made arbitrarily close by appropriate choice of parameters, we can   substitute one for the other in this conclusion, at the expense of replacing  $\delta n$ by $2\delta n$. Since $\tilde y_{r}$ has bounded derivative 
on ${\cal D}_{\eps}$, and $x_m^{(1)}$ can be made arbitrarily close to $x_m$, we deduce that by terminating $\mathcal{A}$ at an appropriate point corresponding to $x_m^{(1)}$, we can achieve the expected size of  output colour 1 to be at most  $C\delta' n$ different from $ \tilde y_{r+1} (x_m)n$, for some constant $C$ independent of $\delta'$. Taking $\delta'<\delta/C$ gives the required result for output colour~1. A similar argument applies to all other variables. 
\qed

\newcommand{\op}{{\rm op}}
\newcommand{\yvtil}{{\bf \tilde y}}
\newcommand{\ytil}{{\tilde y}}
\newcommand{\kbar}{{m}}
 \newcommand{\Pnd}{{\cal P}_{n,r}} 
\newcommand{\aas}{a.a.s.}  
 
For cubic graphs, if~\cite[Theorem 1]{deprio} applies then the deprioritised algorithm exhibits at most two  ``phases.'' There is often only one, as in the example in  the next section, and then the hypotheses are simpler, so we state this special case here. The condition $\xi>0$  ensures that there is only one phase (i.e.\ $m=1$, in the terminology of~\cite{deprio}). 
\begin{thm}\lab{t:cubic} (\cite[Theorem 1]{deprio}, special case) 
 Consider any algorithm consisting of repeated operations $\Op_i$, $i=1,\ldots, 3$ as described above,  applied to $\mathcal{P}_{n,3}$.  For $0\le k\le 3$, let $Y_k(t)$ denote the number of buckets of degree $k$ in the survival graph after $t$ steps of the algorithm. Denote the graph induced by the removed pairs by $G_t$, which therefore has $Y_k(t)$ vertices of degree $3-k$.  Similarly, let $Y_{4}(t)$ denote the cardinality of the output set. Assume that for some fixed functions $f_{k,i}(x,y_1(x),\ldots,y_{4}(x))$ and
for
$k=1,\ldots, 4$, $i=1,\ldots, 3$,
$$
\ex\Big(Y_k(t+1)-Y_k(t) \mid G_t \wedge \{\op_t=\Op_i\}\Big) =
f_{k,i}(t/n,Y_1/n, \ldots , Y_{4}/n) +o(1)
$$
where, for some fixed $\eps>0$, the convergence in $o(1)$ is uniform over all $t$ and $G_t$ for which $Y_i(t)>0$ and $Y_3(t)>\eps n$.  Assume furthermore
that 
\begin{description}
\item{(i)} there is an upper
bound, depending only upon  the number of pairs exposed, and on the
number of elements added to the output set, during any one
operation;
\item{(ii)} the functions $f_{k,i}$ are rational
functions of $x$, $y_1, \ldots, y_{4}$ with no pole in
 the domain
$\D_\eps=\{(x,\yv):
0\le x
\le 3,\ 0\le y_j\le 3 \mbox{ for } 1\le j\le 4,\  y_3\ge\eps\};
$
\item{(iii)} there exist  positive constants $C_1$, $C_2$ and $C_3$ such
that   for
$1\le i < 3$, everywhere on $\D_\eps$, $f_{k,i}\ge C_1y_{k+1}-C_2y_k$ when
$k\ne i$, and
$f_{k,i}\le C_3y_{k+1}$ for all $i$.
\end{description}
Consider $\yv=(y_1,y_2,y_3,y_4)$ and the differential equation 
\bel{de}
\frac{{\rm d}{y_k}}{ {\rm d}x}= 
 \frac{\beta }{\beta +\alpha} f_{k,2}\left(x,\yv\right)
  +\frac{\alpha }{\beta +\alpha } f_{k,1}\left(x,\yv\right)
  \ee
where
$\alpha=  f_{ 1,2}\left(x,\yv\right)$ and $\beta=  -f_{1,1}\left(x,\yv\right)$,
 with initial conditions $\yv(0)=1$ for $k=3$ and $0$ otherwise.
 Define $\yvtil$ to be the solution
of this differential equation, extended for $0\le x\le x_1$, where $x_1$ is defined as the infimum of those $x>0$ for which
at least one of the following holds:  $\beta +\alpha \le \eps$; 
$\ytil_{2}\le 0$; or the solution $(x,\yvtil)$ leaves  $\D_\eps$. Finally, assume that $\beta\left(x,\yvtil\right)> 0$  for $0\le x\le x_1$, and $f_{2,2}(0,\yvtil(0))>0$. 
 
Then there is a randomised algorithm on
${\cal P}_{n,3}$ for which \aas\ there exists  $t$ such that  $Y_4(t)=n\ytil_{4}(x_1)+o(n)$ and $Y_k(t)=n\ytil_{k}(x_1)+o(n)$  for
$1\le k \le 3$. Also 
$\ytil_{1}(x)\equiv 0$  for $0\le  x \le x_1$. \end{thm}

%%%%%%%%%%%%%%%%%%%%%%%%%%%%%%%%%%%%%%%%%%%%%%%%%%%%%%%%%%%%%%%%%%%%%%%%%%%%%%%%%%%%%%
\subsection{Previously analysed deprioritised algorithms}
\lab{s:previousapps} 
In this section we discuss quick and direct applications of Theorem~\ref{t:fromDeprio} using  deprioritised algorithms whose applications to random regular graphs were previously analysed via~\cite[Theorem 1]{deprio}.
  
First we need to address the problem that the output of an algorithm might not be in exactly the right form to solve a problem. For instance,   an algorithm may be aiming to create an output set that is a dominating set, but some of the clash vertices might need to be added to the output set in order for it to dominate all other vertices of the input graph. We say that a set $S\subseteq V(G)\cup E(G)$ in an $n$-vertex graph $G$ is {\em $\eps$-far} from a property $P$ if some set $S'$ in $G$ has $P$ and $|S\Delta S'|\le \eps n$. The number of clash vertices is restricted in Theorem~\ref{t:fromDeprio}. 

 In the remainder of this subsection, we list some of these results.

\subsubsection{Maximum independent set}\lab{Isets thesis}  Asymptotic lower bounds on the size of a maximum independent set (see Section~\ref{s:introresults}) in a random $r$-regular graph were obtained by Wormald~\cite{des}. The simpler of the two algorithms analysed in~\cite{des} was chunkified and then analysed directly for large-girth graphs in~\cite{LW}. 
With Theorem~\ref{t:fromDeprio} in hand, we can easily deduce that the stronger results of the second algorithm also carry over to the large-girth case.

However, slightly better results for the random case were obtained by Duckworth and Zito~\cite{DZ2} using a prioritised native local deletion algorithm that achieves its improvement by looking more carefully into the neighbourhood of the selected vertices. More precisely, given $i \in \{1,\ldots,r\}$, the operation ${\rm Op}_i$ consists of picking a vertex $v$ of degree $i$ u.a.r.\ in the survival graph, adding it to the independent set and removing it along with its neighbours from the 
 survival graph, \emph{unless} the minimum degree $j \in N(v)$ is such that:
\begin{itemize}
\item[(1)] there are at least two vertices of degree $j \leq i$ among the neighbours of $v$, or

\item[(2)] there is a single vertex $u$ of degree $j \leq i$ among the neighbours of $v$, and either
\begin{itemize}
\item[(a)]   $i=2$ and the minimum degree in $N(u)\setminus \{v\}$ is smaller than the minimum degree in $N(v)\setminus \{u\}$, or

\item[(b)] $2<i<r-1$, the minimum degree of $N(u)\setminus \{v\}$ is larger than $i$, and the sum of all degrees in $N(u)\setminus \{v\}$ is smaller than that in $N(v)\setminus \{u\}$.
\end{itemize}
\end{itemize}
If this happens, a vertex $u \in N(v)$ of degree $j$ is added to the independent set, and   is deleted from the survival graph along with its neighbours. It is easy to see that this can be expressed as a native local deletion algorithm, and we leave the definition of the local rule as an exercise for the reader. 

The argument in~\cite{DZ2} applies~\cite[Theorem 1]{deprio}, and the bound on the size of the independent set is   the bound on the size   attained a.a.s.\ by the output set of the algorithm. By Theorem~\ref{t:fromDeprio}, for any fixed $\delta>0$, there is $g>0$ and a chunky local deletion algorithm with the same local and recolouring rules, whose output on $r$-regular graphs of girth at least $g$ is an independent set whose expected size differs from the output of the original algorithm by at most $\delta n$. By letting $\delta\to 0$, we deduce that essentially the same bounds hold in the large girth case.

To state the resulting bounds, we may use the same figures as in~\cite{DZ2} since the quoted figures are, by their nature, some positive constant less than the real number determined from the solution of the associated differential equations. These give lower bounds $\alpha(r)n$ on     the maximum independent set in an $r$-regular graph of sufficiently large girth, which are the best known for all $r\ge 4$. Table~\ref{tableindsets}  shows $ \alpha(r)$ for a few small values of $r$.  For comparison, we provide the best known upper bounds $\alpha_u(r)$ on the best possible values of $ \alpha(r)$. That is, $\alpha_u(r)\ge \alpha$ for all $\alpha$ such that there exists $g$ for which all $r$-regular graphs of girth at least $g$ have an independent set of cardinality at least $\alpha n$.
The values of $\alpha_u(r)$ were obtained by McKay~\cite{mckay} using random regular graphs as described for $c_u(r,P)$ in Section~\ref{s:introresults}.

\begin{table}
\begin{center}
\begin{tabular}{|c||c|c|c|c|c|}
\hline $r$ & 3 & 4 & 5 & 6 & 7\\
\hline $\alpha(r)$ & 0.43475 & 0.39213 & 0.35930 & 0.33296 & 0.31068\\
\hline $\alpha_u(r)$ & 0.45537 & 0.41635 & 0.38443 & 0.35799 & 0.33567\\
\hline $\gamma(r)$ & 0.27942 & 0.24399 & 0.21852 & 0.19895 & 0.18329\\
\hline $\gamma_\ell(r)$ & 0.2641 & 0.2236 & 0.1959 & 0.1755 & 0.1596\\
\hline
\end{tabular}
\end{center}
\caption{Lower bounds $\alpha$  for maximum independent sets, and upper bounds $\gamma$ for minimum independent dominating sets. Best possible values are $\alpha_u$ and $\gamma_\ell$ respectively.}
\label{tableindsets}
\end{table}
The  lower bounds $\alpha(r)$ are now the best known for large girth $r$-regular graphs when $r\ge 4$. Improvements for the case $r=3$ (cubic graphs) are discussed separately in Section~\ref{s:isets}. 
\subsubsection{Minimum independent dominating set}\lab{mids}
Duckworth and Wormald~\cite{DW3,DW} gave upper bounds on the size of a minimum independent dominating set (see Section~\ref{s:introresults}) in a random $r$-regular graph through a deprioritised version of the algorithm described in the introduction. It so happens that these are still the best known upper bounds on the size of {\em any} dominating sets in these graphs. The bounds are based on the size $\rho n$ of the output set of this algorithm (which is an independent set) at the point when the number of vertices in the survival graph falls to $\xi n$ for some preselected $\xi>0$.  An independent dominating set may be obtained by judiciously adding some of the surviving vertices greedily. Clash vertices can also be treated in this way after we apply Theorem~\ref{t:fromDeprio}. We deduce that,  for all $\eps>0$,  there is a deprioritised algorithm is applied to a graph $G$ with sufficiently large girth it a.a.s.\ produces an independent  set that is $\epsilon$-close to an independent dominating set of size at most $\rho n+\xi n$. Since the numerical bounds on $\rho+\xi$ reported in~\cite{DW3,DW} are strictly greater than the true values, we deduce the same upper bounds $\gamma(r)n$ on the minimum independent dominating set in an $r$-regular graph of sufficiently large girth. See Table~\ref{tableindsets} for some specific values of $\gamma(r)$ (more are supplied in~\cite{DW}), and also lower bounds $\gamma_\ell(r)$ obtained via random regular graphs (see Zito~\cite{zito}).
 
The upper bound $\gamma(3)$ improves that obtained by  Kr\'{a}l, \v{S}koda and Volec~\cite{KSV} ($0.299871$). The only existing bounds we know for $r\ge 4$ come from L\"{o}wenstein and Rautenbach~\cite{LR}, who  showed that $n$-vertex graphs with minimum degree at least 2 and girth $g \geq 5$ have a dominating set of size at most $\big(1/3+2/(3g)\big)n$.

\subsubsection{Maximum $k$-independent set}\label{maxkind}

For any positive integer $k$, a $k$-independent set of a graph is a set of vertices such that the minimum distance between any two vertices in the set is at least $k + 1$. Let $\alpha_k(G)$ be the size of a largest $k$-independent set in the graph $G$. Duckworth and Zito~\cite{DZ1}, and  Beis, Duckworth and Zito~\cite{beis_duckworth}  obtained lower bounds on $\alpha_k(G)$ for random $r$-regular graphs through the analysis of an algorithm which satisfies the requirements  for Theorem~\ref{t:fromDeprio}, and whose operations one can easily check to be degree-governed query operations.  Whenever a vertex $v$ is selected, the algorithm tries to find in $G$ one or more  balanced $r$-regular trees of height $\lceil k/2 \rceil$ of which $v$ is a leaf, in a precise way which we do not need to describe here. If it succeeds, it deletes all vertices in the tree and adds its root to the $k$-independent set. Otherwise, it deletes all edges corresponding to open adjacencies that were queried. In this case, the local rule is not normal, as it might be that only some edges of a vertex were queried before finding a new vertex of degree less than $r$, at which point the local exploration is terminated. Similarly to the case of independent sets in Section~\ref{Isets thesis}, Theorem~\ref{t:fromDeprio} shows that the lower bounds  given in~\cite{beis_duckworth} must also be valid lower bounds for the size of the largest $k$-independent set in $r$-regular graphs of sufficiently large girth.
 
\subsubsection{Minimum $k$-dominating set} 

For any positive integer $k$, a $k$-dominating set of a graph $G=(V,E)$ is a set of vertices $S$ such that every vertex $v \in V$ is at distance at most $k$ from $S$. Duckworth and Mans~\cite{duckworth_mans} analysed an algorithm that creates small $k$-dominating sets in random $r$-regular graphs, which may be stated as a native local deletion algorithm as follows. At each step, the selection rule chooses a vertex $v$ with minimum degree u.a.r.\ in the survival graph. The local rule builds a query graph using a normal local subrule: starting from $v$, it unveils the vertex types of its neighbours and chooses an open adjacency with maximum degree, which leads to a new vertex (or closes a cycle, in which case we say that there is a clash). This is repeated iteratively for each new vertex until a vertex $u$ at distance $k$ from $v$ is reached, unless there is a clash or the rule selects a vertex of degree 1, interrupting the process. Finally, the local rule finds the $k$-neighbourhood of $u$ in the survival graph. In the recolouring step, $u$ is added to the output set and the remaining vertices in the copy of the query graph are deleted. 

Using~\cite[Theorem 1]{deprio}, the authors of~\cite{duckworth_mans} analysed a deprioritised version of this algorithm. Implicit in their argument is the requirement that the leftover vertices do not affect their bounds, which they must have checked numerically as for the dominating set example in~\ref{mids}. Applying Theorem~\ref{t:fromDeprio} as for the above examples, these are also valid bounds for $r$-regular graphs with large girth. 

\subsubsection{Maximum  $k$-separated matching}
 For any positive integer $k$, a $k$-separated matching in a graph is a set of edges such that the minimum distance between any two edges in the set is at least $k$. For $k=1$, this is just the maximum matching problem, and random $r$-regular graphs a.a.s.\ have perfect matchings if $r\ge 3$ for $n$ even as they are $r$-connected a.a.s.~\cite{BB,models}. For $k\geq 2$, Beis, Duckworth and Zito~\cite{beis_duckworth} found lower bounds on the size of a maximum $k$-separated matching in a random $r$-regular graphs by analysing an algorithm which is equivalent to one satisfying the requirements of Theorem~\ref{t:fromDeprio}. For $k=2$, this is the maximum induced matching problem, and the algorithm is analogous to the algorithm for independent sets. In each round, it selects  a vertex $u$ of minimum positive degree and then chooses a neighbour $v$ of $u$ with minimum degree among all such neighbours. The edge $uv$ is added to the output set and both $u$ and $v$ are deleted from the survival graph, along with all adjacent vertices. 

For $k \geq 3$, the algorithm resembles what has been done for $k$-independent sets (see Section~\ref{maxkind}). In step $t$, the algorithm selects u.a.r.\ a vertex $v$ with minimum degree and it tries to find one or more balanced $r$-regular trees of height $\lceil k/2 \rceil$ of which $v$ is a leaf. If it succeeds, one edge of each such tree is added to the $k$-separated matching, otherwise all edges corresponding to open adjacencies that have been queried are deleted. The arguments used for independent sets in Section~\ref{Isets thesis} imply that the lower bounds given in~\cite{beis_duckworth} must also be valid lower bounds for the size of a largest $k$-separated matching in $r$-regular graphs of sufficiently large girth.
 
\section{More applications}\lab{s:more}
 
 In this section we present additional applications of our work. The applications in the first and second subsections use Theorem~\ref{t:fromDeprio}, this time for algorithms that had not been analysed before for random regular graphs. The other applications deal with local deletion algorithms directly, and use Theorem~\ref{chunky des} or Theorem~\ref{chunkydeprio large girth}. 

\subsection{Independent sets in large girth cubic graphs}\lab{s:isets} 
The main result in~\cite{VolecFrac} is expressed in terms of the fractional chromatic number, but the basic result proved  is that a certain algorithm produces an independent set with a given (large) expected size. In this section we apply our new analysis to the essentially the same algorithm, achieving basically the same result, but more precise analytically.  A \emph{fractional colouring} of a graph $G$ is an assignment of weights to the independent sets in $G$ such that, for each vertex $v$ of $G$, the sum of the weights of the independent sets containing $v$ is at least 1. The \emph{fractional chromatic number of $c_f(G)$ of $G$ is the minimum sum of weights of independent sets forming a fractional colouring.}  The result on independent sets from~\cite{VolecFrac} was essentially almost the same as the following, but the upper bound on $\alpha$ was obtained by approximate computations. The proof here provides a simple application of the theory presented in the present paper, and also paves the way for Section~\ref{sec:improved}.
\begin{thm}\lab{t:Volec} For all $\alpha<3+\frac32 \log 2 -\frac{15}{2}\sqrt{2}\arctan(\sqrt{2}/4)= 0.43520602\cdots$,  every cubic graph $G$ of sufficiently large girth (depending on $\alpha$) contains an independent set of size at least $\alpha n$, where $n$ is the number of vertices of $G$.
\end{thm}

\proof  The algorithm analysed in~\cite{VolecFrac} is basically a local deletion algorithm of unbounded depth. 
Its local rule can be obtained by the following modification of the greedy (prioritised) independent set algorithm in Section~\ref{s:introresults}, in which minimum degree vertices are randomly placed into the growing independent set.  First, recall  the presentation  of that algorithm in Section~\ref{Isets thesis}. For the case that the minimum degree is 2, Op$_2$ is modified as follows. Instead of placing the randomly chosen vertex $v$ of degree 2 immediately into the independent set, explore along a path $P$ of  degree 2 vertices starting with $v$ (choosing either neighbour at random to start the exploration) until reaching a vertex, $u$, whose other neighbour $w$ has degree 3 (or is $v$).   Insert $u$, together with all vertices at even distance from $u$ backwards along the path $P$, into the independent set, and delete all vertices inserted and their neighbours (which includes $w$).  Although this  would translate into an unbounded depth version of a local rule, its performance can be approximated by terminating after exploring $d$ vertices, and considering $d$ arbitrarily large (but fixed).  Note that a path  of degree 2 vertices whose last vertex is joined to a vertex of degree 1 is immediately ``gobbled up'' by the prioritised algorithm, which via the prioritisation picks on the  vertex of degree 1 and then the other vertices along the path before 
picking on a degree 2 vertex. This causes the above definition of the local rule to be a little simpler than the original.

We claim   that an obvious `chunkified' version of this algorithm is equivalent to the algorithm described in~\cite{VolecFrac}. We do not prove this or even make the statement rigorous, because we will see shortly that analysis of the prioritised algorithm  defined above gives the same result numerically as in~\cite{VolecFrac}.  The present approach  has the advantage of providing a new analytic explanation of that result in terms of the solutions of natural differential equations. 

We will apply Theorem~\ref{t:fromDeprio}. (We could alternatively apply   Theorem~\ref{chunkydeprio large girth} but that would boil down to the almost the same thing, after unwinding the notation.) For cubic graphs, we may refer to the hypotheses of Theorem~\ref{t:cubic}. The main computation required is to find the functions $f_{k,i}$ determined by  $\Op_i$, $i=1$ and $2$.  

In computing the functions $f_{k,i}$, we consider the application to pairings of a local rule  which agrees with the current rule for simple graphs. By our remark after Lemma~\ref{l:trend}, the way in which it is modified to deal with loops and multiple edges does not affect the equations beyond the error term. We can assume $n_k$ is the number of vertices of degree $k$ in the random pairing, and $n_3>\eps n$ (noting condition (ii) in Theorem~\ref{t:cubic}). Then, during an application of the local rule, during which a bounded number of pairs are been exposed, the probability that a randomly chosen point is in a vertex of degree 2 is
$$
p:=\frac{2n_2}{n_1+2n_2+3n_3}  +O(1/n).
$$
In the following, we take $p$ to stand for the expression above, rather than a fixed function. Note that $n_3>\eps n$ implies $p$ is bounded away from 1.

We can assume that $n_1=0$ in all computations of expected changes. This follows from the fact that the solution of the differential equations we eventually derive will have $\tilde y_1=0$ as given by Theorem~\ref{t:fromDeprio}. (The reader unfamiliar with the method may prefer to include the terms for $n_1$, and then note that $n_1$ is so small as to be negligible.)

For $k=1$, 2 and 3, let $\phi_k$ denote the expected change in the number of vertices of degree $k$ in the survival graph due to exposing the second (randomly chosen) point   of a pair and deleting that point. If $x$ denotes the vertex containing the second point, the probability that $x$ has degree 2 is $p$, and for degree 3, it is $1-p$ (assuming $n_1=0$). Hence,  
$$\phi_1=p,\quad \phi_2= -p +(1-p)=1-2p, \quad \phi_3= -(1-p)=p-1.$$

To find $f_{k,i}$, we need to compute the expected change in the number of vertices of  degree $k$ conditional on the history of the process and on the application of operation Op$_i$. We first consider the more complicated case, when Op$_2$ is performed. Let $X$ denote the number of degree 2 vertices in the path $P$ between $v$ and $w$ (not counting $v$ or $w$) and let $p_j=\pr(X=j)$ (all this is conditional on the history of the process). Then the expected change to the vector $(n_1,n_2,n_3)$ in one operation is
(ignoring terms $O(1/n)$, and counting only the event that at most $d$ vertices are explored)
\begin{eqnarray*}
 -\delta_3+2\phi +p_e\big(  p(\phi-\delta_2) +(1-p)(-\delta_3+2\phi)\big)  + (1-p_e)\phi+\sum_{j=0}^{d-1} -(j+1)\delta_2p_j
\end{eqnarray*}
where $\delta_i$ is the unit vector with 1 in the $i$th position,  $\phi=(\phi_1,\phi_2,\phi_3)$ and  $p_e=\pr(X\mbox{ is even})$. Here the term $-\delta_3$ is from deletion of $w$, $2\phi$ is for the two other pairs emanating from $w$ that are deleted,  $-(j+1)\delta_2$ from the deletion of the vertices in the path $P$ as well as $v$, the factor multiplying $p_e$ is for the deletion of the other vertex $z$ adjacent to $v$ (as, in the case $X$ is even, $v$ is inserted into the independent set and $z$ is deleted) whilst the term $(1-p_e)\phi$ accounts for the deletion of just one pair emanating from $v$ in the case $X$ is odd. (The event $w=v$ has probability $O(1/n)$.)

Note that $p_e$ is the sum of $p_j$ for $j<d$ and $j$ even,  and $p_j=p^j(1-p)$. So $p_e=(1+p)^{-1}+o_d(1)$ and, recalling that $p$ is bounded away from 0, $\sum (j+1)p_j=  (1-p)^{-1}+o_d(1)$, where  $o_d(1)$   denotes a function of $d$ (different at different occurrences) that tends to 0 as $d\to\infty$.   Hence, the probability of exploring more than $d$ vertices is $O(p^d)= o_d(1)$, and finally the  expected change to   $(n_1,n_2,n_3)$, ignoring  $o_d(1)$ terms, is
$$
 -\delta_3  +\frac{-p\delta_2 - (1-p)\delta_3  +  2(2+p)\phi }{1+p}   -\frac{\delta_2}{1-p}.
$$
This vector is (using Maple for this and later calculations)
\bean
&&\left( \frac{  2p(2+p) }{1+p}, \
\frac{-p + 2(2+p)(1-2p) }{1+p}   -\frac{1}{1-p},\
-1 +\frac{ - (1-p)  +  2(2+p)(p-1) }{1+p}   
\right)\\
&&=\left( \frac{  2p(2+p) }{1+p}, \
\frac{3-12p+3p^2+4p^3 }{1-p^2} ,\
 \frac{ 2(-3+p+p^2) }{1+p}   
\right),
\eean
which determines $(f_{12},f_{22},f_{32})$ up to $o_d(1)$, as long as for this interpretation we   drop the $O(1/n)$ error term from the definition of $p$ in order to make the functions $f_{i,j}$ independent of $n$.
For Op$_1$, just one vertex of degree 1 is deleted, and the expected change vector is 
$$-\delta_1+p(-\delta_2+\phi) +(1-p)(-\delta_3+2\phi)
= (-1+2p-p^2,2-6p+2p^2, -3+4p-p^2),
$$
giving $(f_{11},f_{21},f_{31})$ (subject to the same interpretation as above). It only remains to compute the expected change to the size of the output set. In   $\Op_2$, the probability that at least $j$ vertices are added to the output set is  $p^{2j}$ ($j<d/2$), and hence the expected number added is 
$f_{42}=1/(1-p^2)+o_d(1)$. Clearly $f_{41}=1$. It is now relatively straightforward that the hypotheses of Theorem~\ref{t:cubic} hold with these functions $f_{k,i}$ provided the solution leaves $\D_\eps$ before the other termination conditions come into force. This can be checked   (at least, for sufficiently large $d$) by first solving the ``pure" differential equation system
%\bean
%y_2'&=& \frac{(1-3p)(p^2+2p+3)}{p^3+ p^2+3p+1 },\\
%y_3'&=& \frac{2(p-1)(p^2+2p+3)}{p^3+ p^2+3p+1 },\\
%y_4'&=& \frac{(1+p)(1+2p)}{p^3+ p^2+3p+1 } 
%\eean
$$
y_2'= \frac{(1-3p)(p^2+2p+3)}{p^3+ p^2+3p+1 },\quad
y_3'= \frac{2(p-1)(p^2+2p+3)}{p^3+ p^2+3p+1 },\quad
y_4'= \frac{(1+p)(1+2p)}{p^3+ p^2+3p+1 } 
$$
where $p=2y_2/(2y_2+3y_3)$ (which is the interpretation noted above). Recall that $y_1'$ must be identically 0. We deal with the $o_d(1)$ terms later.

Putting $s=2y_2+3y_3$, and noting that $p'= 2y_2'/s-s'p/s$, we find
$$
p'=\frac{2(1-p)(p^2+2p+3)}{s(p^3+ p^2+3p+1)}
$$
and then that ${\rm d}p/{\rm d}s =   p'/s'= (p-1)/2s$. Hence (using $s=3$ when $p=0$, which holds when $x=0$) we find $s=3(1-p)^2$.  Thus
\bel{y4deriv}
\frac{{\rm d}y_4}{{\rm d}p}=\frac{ 3(1-p)(2p^2+3p+1)}{ 2(  p^2+2p+3)}.
\ee
As  $s=3(1-p)^2$, we have $y_3=(1-p)^3$, and we may now deduce, using $y_3'<0$ for $p<1$, that $(x,\yvtil)$ leaves $\D_\eps$ when $y_3=\eps$ and $p=1-O(\eps)$. At this point, $x=x_1$.
Since the system is well behaved (the derivative functions are Lipschitz) inside  $\D_\eps$, the solutions of the ``true'' system, which contains the extra $o_d(1)$ terms, differ from the solutions of the pure system by $o_d(1)$. Hence, since $y_3'$  is bounded away from 0 when $p$ is bounded away from 1, the point of exiting $\D_\eps$ is at $x=\bar x_1=x_1+o_d(1)$. 

Let $\hat x_1$ denote the value of $x$ at which $p$ reaches 1. 
Then equation~\eqn{y4deriv} can be integrated from $0$ to 1 to yield
 $$y_4(\hat x_1)=  3+\frac32 \log 2 -\frac{15}{2}\sqrt{2}\arctan(\sqrt{2}/4)= 0.43520602\cdots.$$
By taking $\eps$ sufficiently small,   $x_1$ can be made arbitrarily close to $\hat x_1$.
      The constant $\rho$ referred to in Theorem~\ref{t:fromDeprio} is $y_4(\bar x_1)$, and we have seen that this can be made arbitrarily close to $y_4(\hat x_1)$ by taking $d$ large and $\eps$ small. Theorem~\ref{t:fromDeprio} now implies the stated result. \qed

%%%%%%%%%%%%%%%%%%%%%%%%%%%%%%%%%%%%%%%%%%%%%%%%%%%%%%%%%%%%%%%%%%
\subsection{Independent sets: improved result}\lab{sec:improved}

As mentioned in~\cite{VolecFrac}, the approach used for Theorem~\ref{t:Volec} can be extended. We show here that this is relatively painless with our new machinery.

After exploring locally in the neighbourhood of a vertex, it is not at all obvious what the best choice is for placing vertices into the output set. Observe that the previous theorem  improves the result from the earlier algorithm, suggesting the paradigm of inserting $u$ into the independent set (along with every second vertex along the path) in preference to the vertex adjacent to $u$ in $P$. If this is attempted  for both paths of degree 2 vertices starting at $v$,   and the two vertices placed in the output set have odd distance apart, this strategy is clearly non-optimal, and the other choice should be made in one of the paths. Given the information so far, either path is as good for the adjustment as the other.  However, exploring the paths emanating from $w$ gives more information, and the same paradigm can be applied to those paths.  One way of extending this idea using the paradigm gives the following result. (The constants arising in this theorem are approximations of constants determined analytically.)
\begin{thm}\lab{t:Volec:imp}
The constant in Theorem~\ref{t:Volec} can be improved to $\alpha_0=0.43757463$. Moreover, the fractional chromatic number of a cubic graph with sufficiently large girth is at most $1/\alpha_0\approx 2.2853$. 
\end{thm}
\proof
Modify the local rule in the proof of Theorem~\ref{t:Volec} as follows, in the case of $\Op_2$. Explore the maximal path $P$ of vertices of degree 2 containing $v$ (so $P$ can extend in two directions from $v$). Let  $w_1$ and $w_2$ denote two vertices outside $P$ which are adjacent to the ends of $P$. See Figure~\ref{f:cubic}.  The case that one of these vertices $w_i$ is equal to $v$, and indeed any cases in which cycles are discovered, can be ignored, as seen in the proof of Theorem~\ref{t:Volec}. So we may assume that both $w_1$ and $w_2$ have degree 3.

 \begin{figure}[hbt]
\centerline{\includegraphics[width=4.5cm]{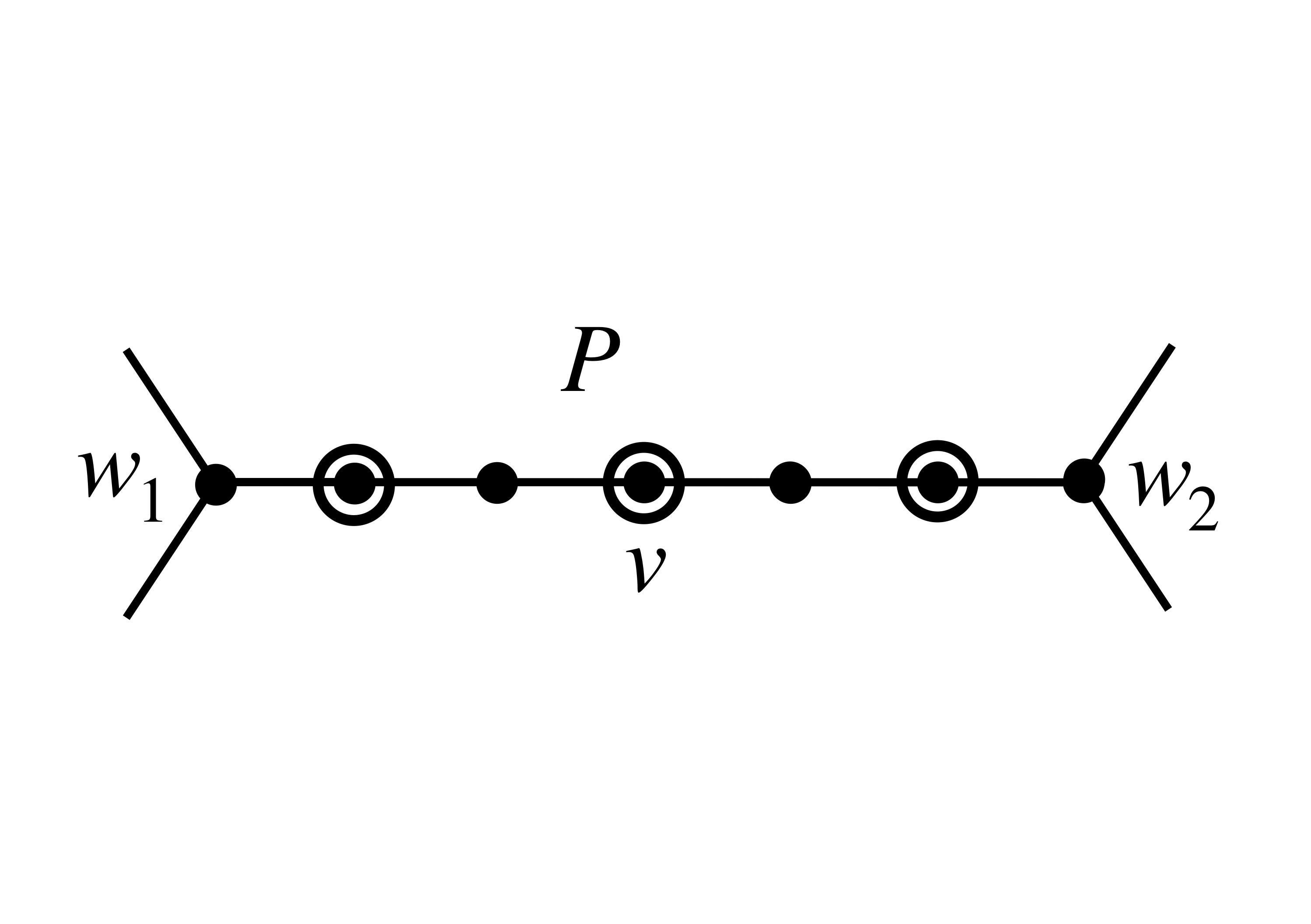}\hspace{2.7cm}   \includegraphics[width=4.5cm]{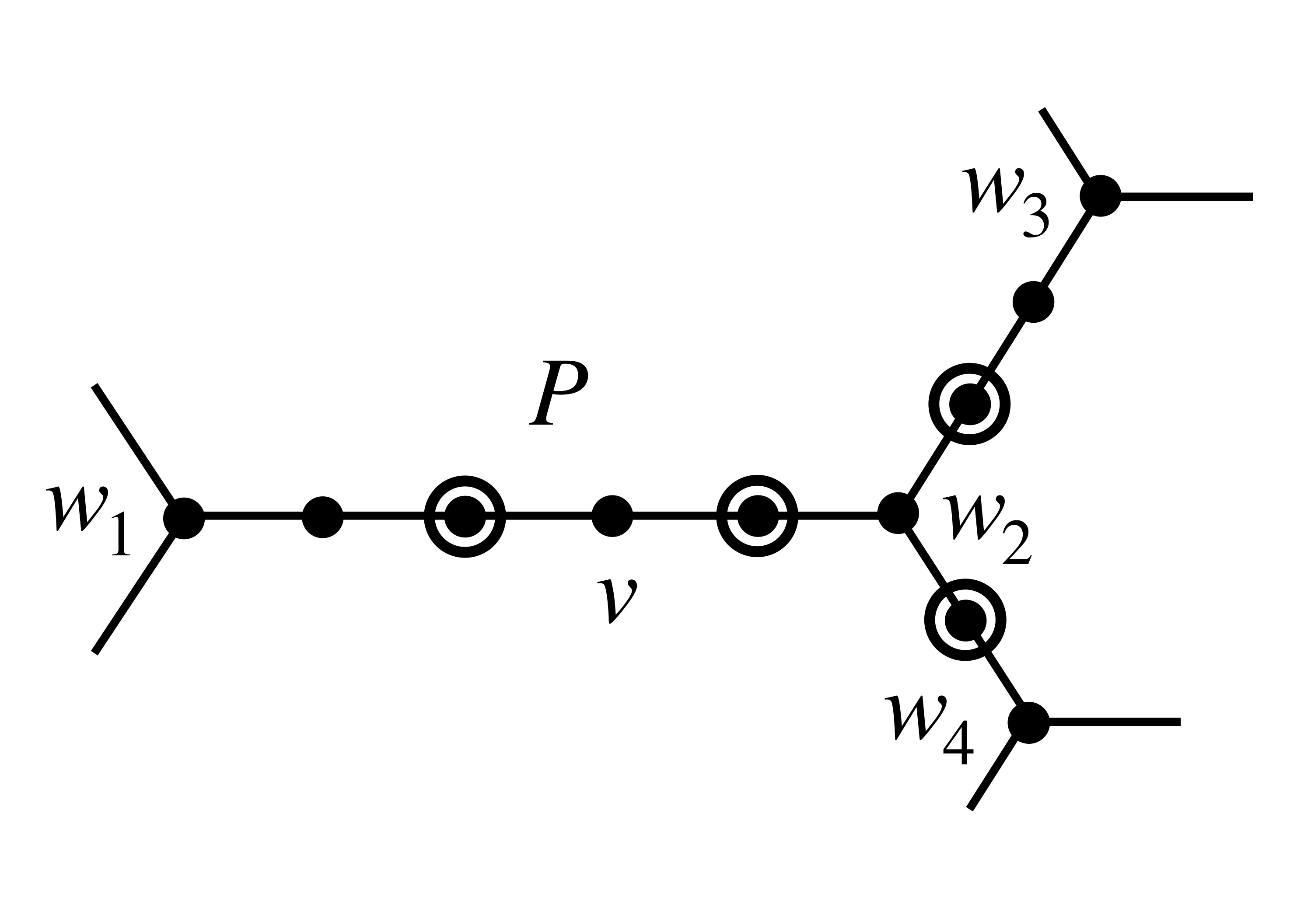}}
\vspace{-0.4cm}

\hspace{4.2cm}(a)\hspace{6.7cm}(b)

\centerline{\includegraphics[width=4.5cm]{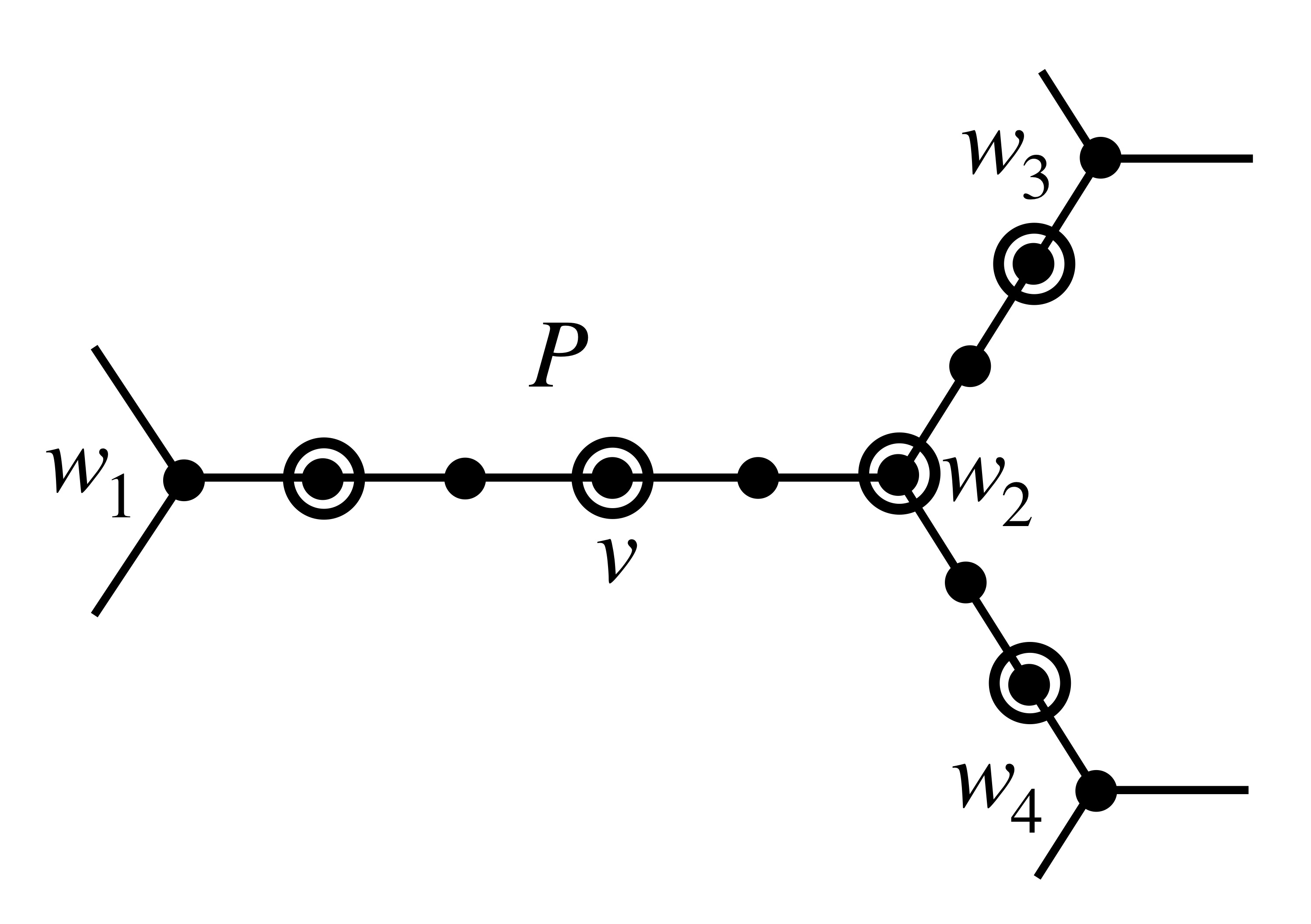}\hspace{2.7cm}   \includegraphics[width=4.5cm]{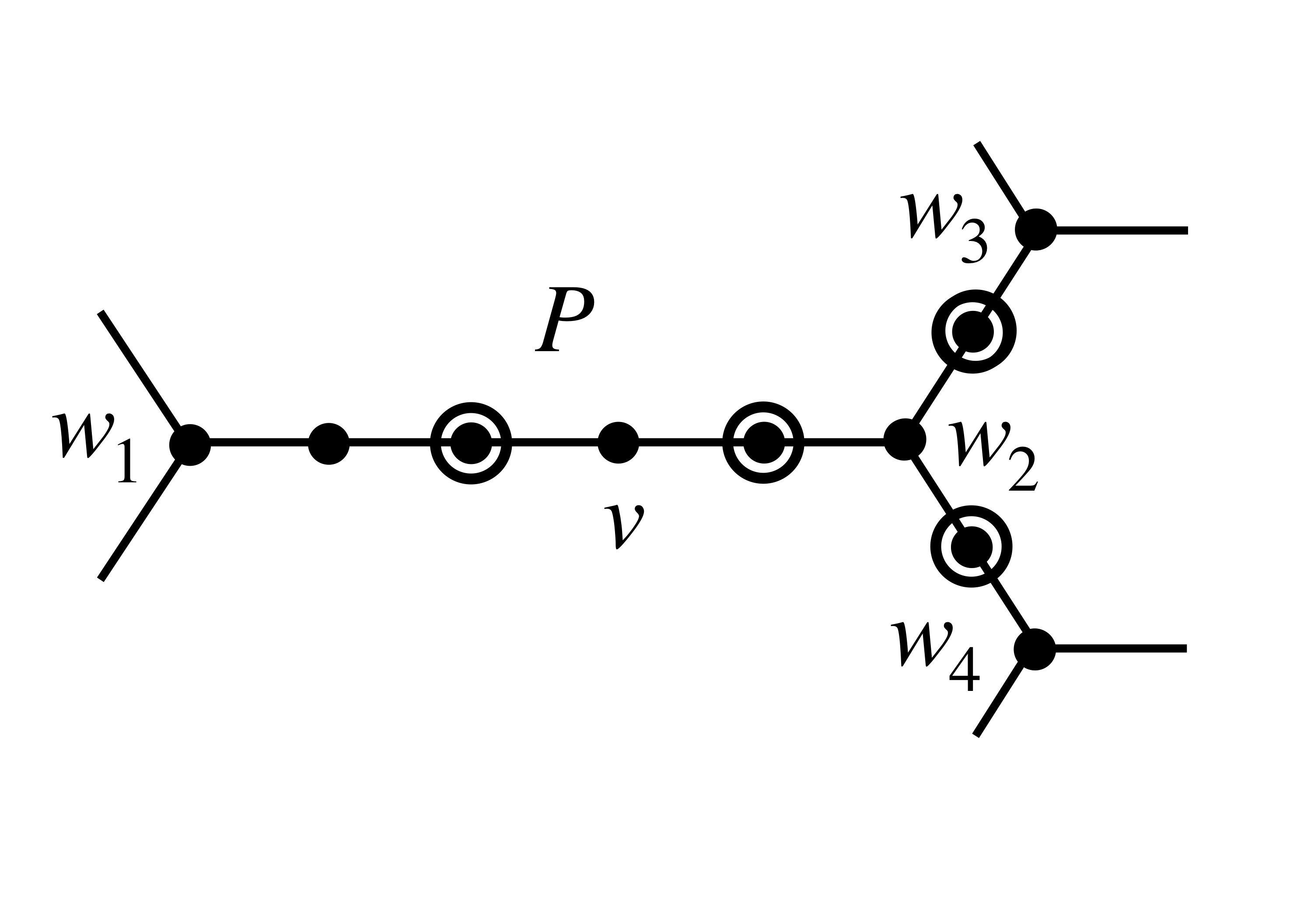}}
\vspace{-0.6cm}

\hspace{4.2cm}(c)\hspace{6.7cm}(d)
\caption{Exploration}\label{f:cubic}
\end{figure}

In each case, vertices are placed alternately in and not in the output set, and we only need to mention the extremal ones placed in it, which we will call terminal vertices. 
If $P$ has an odd number of vertices, as in Figure~\ref{f:cubic}(a), the vertices of $P$ adjacent to $w_1$ and $w_2$ are the terminal vertices. Otherwise,   explore similarly the path $P'$ containing $w_2$ in the graph with $P$ removed. Let $w_3$ and $w_4$ be the vertices of degree 3 adjacent to its ends. If $P'$ has an even number of vertices, as in Figure~\ref{f:cubic}(b), the choice of terminals on $P'$ clearly has no effect at $w_3$ and $w_4$ (one will be adjacent to a terminal, and one not) so it is best to place vertices of odd distance from $w_2$ into the output set (restricting to vertices in $P$ and $P'$). If $P'$ has an odd number of vertices,  as in Figure~\ref{f:cubic}(c) and (d), then its end vertices will be nominated as terminals (following the paradigm described above), and $w_2$ (and hence the vertex adjacent to $w_1$) will be in the output set  if and only if  $w_2$ has odd distance from $w_3$ (and consequently from $w_4$).

We can derive the functions $f_{i,j}$ for this version of $\Op_2$ using the same argument as in the previous theorem, ignoring $o_1(d)$ terms for the same reason. The probability $P$ has $2i-1$ vertices is
$$
h_{2i-1}=\sum_{j=1}^{2i-1} p_{j-1}p_{2i-j-1} = (2i-1)(1-p)^2p^{2i-2}.
$$
In this event the expected change in $(n_1,n_2,n_3)$ is (ignoring those negligible terms)
$$
-(2i-1)\delta_2 -2\delta_3+4\phi.
$$
Summing over all $i$ gives the part of the expected change due to the case of $|V(P)|$  odd:
\bel{caseo}
\left(\frac{4p(p^2+1)}{(1+p)^2},\frac{8p^5-5p^4-6p^2-8p+3}{(1-p) (1+p)^3}, 
\frac{2 (2p-3 ) (1+p^2)}{(1+p)^2} \right).
\ee
The expected change in the output set size from this case is 
$$\sum_{i\ge 1}ih_{2i-1}= \frac{1+3 p^2}{(1-p) (1+p)^3}.
$$

Similarly, the probability $P$ has $2i$ vertices is $ 
h_{2i }= 2i(1-p)^2p^{2i-1}$; in this case, let the number of degree 2 vertices  between $w_2$ and $w_3$ be $j$, and between $w_2$ and $w_4$,  $k$. The probability of this event is $h_{2i}p_jp_k$.

There are three subcases. If $j$ and $k$ are both even, let them be $2i_1-2$ and $2i_2-2$ respectively (it is convenient to keep $i_1$ and $i_2\ge 1$). Then  the expected  change in $(n_1,n_2,n_3)$ is $-(2i+2i_1+2i_2-4)\delta_2-4\delta_3+6\phi$, and in the output set it is $i+i_1+i_2-1$ as $w_2$ is included. See Figure~\ref{f:cubic}(a).  Multiply these expressions by $h_{2i}p_jp_k$ and sum over all relevant $i$, $i_1$ and $i_2$ to get 
\bel{caseeee}
\left(\frac{12p^2}{(1+p)^4},\frac{8 p (3 p^3-3 p^2-3 p+1)}{(1+p)^5 (1-p)},
\frac{4p(3p-5)}{  (1+p)^4} 
 \right) 
\ee
and
$$
\frac{4 p ( p^2  +1)}{(1+p)^5 (1-p)}.
$$
If just one of $j$ and $k$ is even and the other odd, let them be $2i_1-2$ and $2i_2-1$. See Figure~\ref{f:cubic}(b). The expressions are $-(2i+2i_1+2i_2-5)\delta_2-4\delta_3+2\phi$ (note that $2i+2i_1+2i_2-5$ vertices of degree 2 are deleted and two are created, $w_1$ and either $w_3$ or $w_4$), and   $i+i_1+i_2-1$. Multiplying by $2h_{2i}p_jp_k$ (accounting for two cases, $j$ or $k$ being even) and summing gives  
\bel{caseeeo}
\left(\frac{8p^3}{(1+p)^4},\frac{4p^2(4  p^3-9 p^2-4 p+1)}{(1+p)^5 (1-p)},
\frac{8p^2( p-3)}{  (1+p)^4} 
 \right) 
\ee
and
$$
\frac{8 p ( p^2  +1)}{(1+p)^5 (1-p)}.
$$
If both $j$ and $k$ are odd, let them be $2i_1-1$ and $2i_2-1$. See Figure~\ref{f:cubic}(c). The expressions are    $-(2i+2i_1+2i_2-3)\delta_2-4\delta_3+4\phi$ (here $w_3$ and $w_4$ are deleted), and   $i+i_1+i_2$. Multiplying by $ h_{2i}p_jp_k$  and summing gives  
\bel{caseeoo}
\left(\frac{8p^4}{(1+p)^4},\frac{2p^3(8 p^3-9 p^2-8 p+1)}{(1+p)^5 (1-p)},
\frac{8p^3( p-2)}{  (1+p)^4} 
 \right) 
\ee
and
$$
\frac{2 p^3 ( p^2  +3)}{(1+p)^5 (1-p)}.
$$
Summing the vectors in~(\ref{caseo}--\ref{caseeoo}) gives the vector $(f_{12},f_{22},f_{32})$ up to $o_d(1)$. Continuing with the method as before, we find once again that
${\rm d}p/{\rm d}s = (p-1)/2s$ and then 
$$
\frac{{\rm d}y_4}{{\rm d}p} = \frac{ 3 (1-p)(2p^5+9p^4+18p^3+22p^2+8p+1)}{ 2 (p^4+4p^3+8  p^2+14p+3)(1-p^2)}.
$$
The integral of this can be expressed as a sum over the roots of a quartic polynomial, and Maple evaluates the integral from 0 to 1 as $0.43757463\cdots.$

As remarked in~\cite{VolecFrac}, if we can define a random independent set $X$ in a graph $G$ such that each vertex has probability at least $\alpha$ of appearing in $X$, then
the fractional chromatic number $c_f(G)$ of $G$ is at most $1/\alpha$. Theorem~\ref{t:Volec:imp} shows that $\alpha_0$  (or any constant less than the constant $y_4(1)$ defined in its proof) is a lower bound on the size of a maximum independent set in a cubic graph of sufficiently large girth. The proof via Theorem~\ref{t:fromDeprio} uses a chunky local deletion algorithm which produces a random independent set whose expected size is at least $\alpha_0$.  The subgraph induced by the vertices of $G_t$ of distance at most $4D$ from $v$, where the local rule has depth $D$,  determines precisely the probability that $v$ is added to the independent set in the next step of the algorithm. (In the extreme case, $v$ has distance $D$ from the root of the query graph,   some other vertex of which at distance $2D$ from $v$ is a clash vertex.) Corollary~\ref{gen ind lab}(i) implies that this probability  does not depend on $v$.   Hence every vertex is in the independent set with probability at least $\alpha_0$, and so by the conclusion above, $c_f(G) \leq 1/\alpha_0\approx 2.285325$. 
\qed

\subsection{Large cuts in random cubic graphs}\lab{s:cuts}

In~\cite{Volec} it is shown that every cubic graph with $n$ vertices and girth at least $g_0 = 16353933$ contains an edge cut of size at least $1.33008n$. This used a computation that numerically iterated a recurrence $g_0$ times. Since the method was based on~\cite{hoppen_thesis}, which employed some of the ideas in the present paper (including some use of deprioritised algorithms), it is not too surprising that the algorithm in~\cite{Volec} is a chunky local deletion algorithm that derives in a natural way from the prioritised algorithm $\mathcal{A}$ described at the end of Section~\ref{s:localCol}. (However,  $\mathcal{A}$ is not  directly described in~\cite{Volec}.) We analyse an amenable deprioritised algorithm whose behaviour approximates that of  $\mathcal{A}$  to find a lower bound on the size of a maximum cut in a cubic graph. The algorithm used in~\cite{Volec} can be viewed as just one example of a family of such deprioritised algorithms that we use.
  
 In order to describe our result, we define the system of differential equations
\bel{decut}
u'=  \frac{-6u(6u+s)}{d},\ v'=  \frac{ 36u^2-12uv-s^2-4sv}{d},\ 
w'= \frac{2 (-6uw+2 sv-sw)} {d} 
\ee
where $s= 3u+2v+w$ and $d=s^2+6su+12uv$. 
\begin{thm}\lab{t:cut}
Let $(u,v,w)$ denote the solution of~\eqn{decut} as functions of  $x\ge 0$ with initial conditions $(u,v,w)=(1,0,0)$  at $x=0$. Let $\hat x_0$ denote the minimum $x>0$ for which $v(x)=0$. Let $z(x)$, $0\le x\le \hat x_0$, satisfy
$$
z' = \frac{ 6su+24uw+s^2+4sw+36uv}{d}, \quad z(0)=0.
$$
Then, for all $\delta>0$,  every $n$-vertex cubic graph with sufficiently large girth has an edge cut of size at least $cn$ where
$c=z(\hat x_0)+\frac32  u(\hat x_0)+\frac32w(\hat x_0)-\delta$.
\end{thm}
\noindent{\bf Note:}  For small $\delta$, using a standard second order Runge-Kutta solution method we get $c\approx 1.330209040$. (Here $\hat x_0\approx 0.8274171475$.)
\medskip

\proof
We will define an appropriate amenable deprioritised algorithm related to $\mathcal{A}$ for use with Theorem~\ref{chunkydeprio large girth}, beginning with a small $\delta>0$ and   $\eps>0$ whose relationship will be specified later. 

 We can streamline the argument with the foresight  that, because of the chosen prioritisation, only the types 00, 01 and 11 will be present in any great numbers at a given time. The local subrule is normal, so after a vertex is probed, the types of all its neighbours are determined. Only one vertex is probed in the local rule, with the exception that if a vertex is going to be given an output colour, and has a neighbour $u$ with transient colour $r'b'$ where $r'+b'=2$,  then $u$ is probed and included in the query graph. Its (output) colour is determined and none of its other neighbours are affected; since they have nontransient colours they are not in the survival graph. That is, each $\Op_{rb}$ with $r+b=3$ is amalgamated with the previous operation. 
 
Before defining the relative selection functions, we consider the $f_{j,i}$ defined in Lemma~\ref{l:trend}. In the present case these are   $f_{k\ell,ij}$ where $k\ell$ is the type of vertices whose expected changes are being tracked, and $ij$ is the type of the operation ($i+j\le 2$ given the amalgamation described above). We use $y_{k \ell}$ for the variables in the argument of an $f$ relating to the (scaled) number of vertices of type $k\ell$. Directly considering the pairing model as in Section~\ref{s:isets},  and assuming that $i>j$, which applies when the operation colours the deleted vertex blue, it is easy to derive
\bel{fcuts}
f_{k\ell,ij}(\y)= -\delta_{ k\ell,ij}+(3-i-j)\left(\frac{(3-k-\ell+1)y_{k,\ell-1}}{s}-\frac{ (3-k-\ell)y_{k  \ell}}{s}\right)
\ee
where  $s= \sum_{ij}(3-i-j)y_{ij}$, and $\delta_{ k\ell,ij}$ is the Kronecker delta: 1 if $k=i$ and $\ell=j$, and 0 otherwise. Of course, $s$ represents the number of unpaired points in the pairing. When $i>j$ the formula is the same as~\eqn{fcuts}, except that $y_{k,\ell-1}$ is replaced by $y_{k-1,\ell}$.

We are at liberty to define any relative selection functions, but we provide some motivation without full justification. Following the general argument in~\cite{deprio}, we note that in the prioritised algorithm the operations, at least reasonably early in the process, after an initial one of type $00$, will consist of types  $10$, $01$ and $02$.  We want all other operations to have probability 0, and we desire   the derivatives of $y_{10}$ and $y_{02}$ to be 0 (so these variables remain 0).  We therefore write the equations 
$$
\tilde p_{10}f_{\tau,10}  + \tilde p_{01}f_{\tau,01} + \tilde p_{02}f_{\tau,02}=0,
\quad(\tau = 10 \mbox{ and } 02)
$$
and solve these simultaneously with $\tilde p_{10}   + \tilde p_{01}  + \tilde p_{02}=1$, to obtain
\bel{cutselections}
\tilde p_{01}= \frac{s^2-12y_{00}y_{01}}{d}, \quad \tilde p_{10}= \frac{6sy_{00} }{d},\quad   \tilde p_{02}= \frac{24y_{00}y_{01}}{d}, 
\ee 
where $s=3y_{00}+2y_{01}+y_{11}$ and $d=s^2+6sy_{00}+12y_{00}y_{01}$.

Start with a `burn-in' phase of length $\eps>0$; that is, for $x\in [0,\eps]$,   define the relative selection functions $\tilde{p}_{rb}(x)=1$ if $rb=00$ and $\tilde{p}_{rb}(x)=0$  otherwise. The differential equation system~\eqn{desystem} then has a unique solution $\tilde {\bf y}$ on this interval.   For $x\ge \eps$, we need to use $\tilde {\bf y}(\eps)$ to define the selection functions. To this end, consider the differential equation system~\eqn{desystem} where the {\em functions} $\tilde p_{ij}$ are defined as in~\eqn{cutselections} as functions of $\bf y$, and with initial conditions being ${\bf y}(\eps)=  \tilde {\bf y}(\eps)$. All $\tilde p_{ij}$ not appearing in~\eqn{cutselections} are defined to be 0.  The  solution $\tilde {\bf y}$ then extends until it approaches a singularity of the derivative functions. With this definition of $\tilde {\bf y}$, we may now define the relative selection functions $\tilde p_{10}$ as the functions of $x$ determined from~\eqn{cutselections}. To distinguish these from the functions $\tilde {\bf y}(\eps)$ of $\bf y$ referred to in the computation, we denote them $\tilde p_{ij}^{(1)}(x)$. Then the system~\eqn{desystem}, with relative selection functions $\tilde p_{ij}^{(1)}$, and the same initial conditions, still has unique solution $\tilde {\bf y}$. 

We may now apply Theorem~\ref{chunkydeprio large girth}, for appropriate $\delta$. The crucial question issue is the determination of $x_0$. By definition this is determined by the first point at which some variable $y_{ij}$ is below $\delta$ when the corresponding relative selection function is positive. The only such functions that are ever positive are for $ij= 00$ for $x\le \eps$, and $10$, $01$ and $02$ thereafter. By the nature of the burn-in phase it is easy to see that there is some $\eps>0$ such that $y_{00}>1/2$ on $[0,\eps]$, and the remaining variables exceed $\delta$ at $x=\eps$, where $\eps\to 0$ as $\delta\to 0$. For $x\ge \eps$, for $ij \in \{10,02\}$,  the derivative of $y_{ij}$ is precisely 0 by the equations determining the relative selection functions. So $x_0$ is determined by the first point after $\eps$ such that $y_{01}$ falls below $\delta$. Noting that $y_{01}'>0$ for small $x\ge 0$, the point of interest is when $y_{01}$ falls back close to 0.

The differential equation system described above for the interval $x\ge \eps$, involving $\tilde p_{ij}$ as functions of $\bf y$, would be exactly the system~\eqn{decut} with $(u,v,w)=(y_{00},y_{01}, y_{11})$ if we set all other variables $y_{ij}$ to zero.  If the  initial conditions are shifted slightly so that $(1,0,0)$ applies to both systems, it follows that the trajectories of the systems would be identical: the only other variable we have not explicitly guaranteed to have zero derivative is $y_{20}$, and this remains zero because $y_{10}$ does. Moreover, until $y_{01}$ returns to 0 it is easy to prove that $s$ and $d$ remain bounded away from 0 and so the derivative functions are all Lipschitz in a neighbourhood of the solution trajectories. Hence, by well known properties of ordinary differential equations, by taking $\eps$ and $\delta$   sufficiently small, the two systems' solutions can be brought arbitrarily close to each other. It follows that, in addition, $x_0$ can be made arbitrarily close to $\hat x_0$.

Let $y_R$ be the variable corresponding to the output edges. It is straightforward to check that $y_R'$ corresponds to the derivative of $z$ given in the theorem, as follows. An operation of type 10 or 01 in expectation increases the output edge set by $1+4y_{11}/s$, where 1 is for the edge to the coloured neighbour and $4y_{11}/s$ is twice the expected number of white neighbours of type $br$ (there being two such neighbours potentially). All other neighbours can be neglected by the argument above. Similarly, type $02$ increases the output edge set by  $2+2y_{11}/s$ in expectation.

We deduce from  Theorem~\ref{chunkydeprio large girth}  that for all  $\delta>0$   there is a chunky local deletion algorithm $\mathcal{A}_1$ with the same local and recolouring rules, which, applied to any cubic graph of sufficiently large girth,  for which  the output set of edges, after a number of steps corresponding to $\big(\hat x_0+O(\delta)\big)n$, has  size   at least  $z(\hat x_0)n-\delta n$ and the numbers of vertices of types $00$ and $11$ are similarly close to $u( \hat x_0)n$ and $w(\hat x_0)n $ respectively. Moreover, the numbers of vertices of any other type os at most $\delta n$. From numerical solution of~\eqn{decut}, we find  $ \hat x_0 \approx 0.8274$, $u( \hat x_0)\approx 0.00279$ and $w( \hat x_0)\approx 0.0511$.

 In~\cite{Volec}, the algorithm runs a second phase where, in the absence of any other types besides 00, a vertex of type 11 is randomly coloured red with probability $1/2$ and blue otherwise. However, we can easily avoid analysis of this phase (and also explain why the choice of $1/2$ for this probability makes no difference to the final result) using the results in Section~\ref{s:explicit}. By Lemma~\ref{l:pairings to graphs}, we may start with a random pairing instead of a random graph, from the beginning of the present argument, to reach the same conclusions as above (indeed, this is the method of the proof  of Theorem~\ref{chunkydeprio large girth}). Then Lemma~\ref{l:survival} allows us to analyse the remaining part of the process by taking a random pairing with the given numbers of vertices of types $00$ and $11$, and few other vertices. Since the number of type $00$ vertices is approximately $u( \hat x_0)n$ and type $11$ approximately $ w( \hat x_0)n$, there are approximately 17 times as many of the latter as the former. So at this point we may abandon the local deletion algorithm and consider instead the structure of this random pairing. Starting at any vertex and exposing an unexposed pair, the expected number of unexposed points in the new vertex, after exposing this pair, is at most $2/17$. Continuing the exploration of a component in the random pairing using the standard analysis for random graphs using branching processes, we find the expected size of each component is bounded, and the probability there are more than $M$ vertices in the same component as a given start vertex is at most $e^{-CM}$ for some constant $C$. The probability of finding a cycle before reaching $M$ vertices in a component is $O(M^2/n)$. Thus, letting $X_M$ denote the number of vertices in components that are not trees of size at most $M$, we have $\ex X_M = O( e^{-C'M})$. Markov's inequality implies $\pr (X_M/n >\delta ) =  O(  e^{-C'M}/\delta)$. Lemma~\ref{l:pairings to graphs} allows us to transfer this statement to the random graph, rather than pairing, if we write the bound as  $e^{-C''M}$  where $C''$ is  another constant.
 
Every tree component in the graph can be bicoloured. Since all but $O(\delta n)$ vertices are of type either $11$ or $00$, which have degree 1 and 3 in the tree respectively, we get $(3/2)u( \hat x_0)n+ (1/2)w( \hat x_0)n +O(\delta n)$ edges of type red-blue internally in such components. Additionally, each vertex of type $11$ gives one more red-blue edge externally. The total is  $(3n/2)\big( u( \hat x_0)+w( \hat x_0)\big) +O(\delta n)$ edges which can be added to the output set. As $\delta$ can be taken arbitrarily small, the result follows.
 \qed

We note that the analysis of the ``second phase'' in the above proof shows that all ways of performing that phase by properly bicolouring all components of size at most $M$ will achieve the stated bounds for the output set, provided the granularity is small enough to make clash vertices insignificant. We can easily define such chunky algorithms,  for instance based on the rule referred to in~\cite{Volec} where  choosing type 11 has the lowest priority, and one could even use the rule that all such edges are coloured blue.
The argument used for independent sets in Section~\ref{sec:improved}, is easily applied to show that the probability that any given edge lies in the cut produced by such an algorithm is at least $0.88680$. This leads to the upper bound of $1.12765$ (almost the same as in~\cite{Volec}) on the (suitably defined) fractional cut covering number of a cubic graph with large girth. 

\subsection{Minimum and maximum bisection}\label{s:bisection}
 
 Given a graph $G=(V,E)$ with $|V|$ even, a bisection of $G$ is a partition of its vertex set into two parts with the same cardinality, and the size of the bisection is the number of edges crossing between the parts. D\'{i}az,   Do,  Serna and Wormald~\cite{DDSW,DSW} obtained upper bounds on the size of a minimum bisection, also known as the bisection width, and lower bounds on the size of a maximum bisection in a random $r$-regular graph. Table~\ref{biswidth} shows the results, where $\beta(r)n$ is an upper bound a.a.s.\ on the bisection width of a random $r$-regular graph, as $n\to\infty$. In the case of $r=4$, $\eps$ is any positive quantity. By a quirk of the proof method, $(r/2-\beta)n$ provides a corresponding asymptotically almost sure lower bound  on the size of the maximum bisection. For $r=3$, the bisection width upper bound is weaker than the deterministic $(1/6+\eps)n$ that  Monien and Preis~\cite{MonPreis01} showed for {\em all} 3-regular graphs, but the complementary bound $1.3259n$ on maximum bisection is still interesting.

 We argue that, with our approach, the applicability of these bounds extends to all regular graphs of sufficiently large girth. Since it is easy to adapt the arguments of Section~\ref{s:cuts} to this situation, we only give a general description of the proof.
 
\begin{table}
\begin{center}
\scalebox{0.9}{\begin{tabular}{|l|l|l|l|l|l|l|l|l|l|l|}
\hline
        r & 3& 4 & 5& 6 &7 & 8 &9 & 10 & 11 & 12 \\
  \hline
$ \beta(r)$&0.1741&$1/3+\eps$ & 0.5028  & 0.6674 & 0.8502  & 1.0386  & 1.2317         & 1.4278
& 1.624 & 1.823
       \\
\hline
\end{tabular}} 
\ \\
\end{center}
\ \\
\caption{Upper bounds on bisection width\label{biswidth}}
\end{table}

To derive the original bounds, D\'{i}az \etal~analysed a randomised procedure that resembles the algorithm for large cuts of Section~\ref{s:cuts}. Indeed, it may be described as a local deletion algorithm with the same local rule, but different selection and recolouring rules. There are again two output colours, denoted red and blue, to mark vertices already assigned to each class of the partition, while the remaining vertices have transient colours indicating the number of neighbours in each of these classes. Let $rb$ denote the type of a vertex with $r$ red neighbours and $b$ blue neighbours. At each step, the procedures in~\cite{DDSW,DSW} choose two vertices, which are called a symmetric pair, with types $xy$ and $yx$, respectively. The values of $x$ and $y$ are determined according to a priority ordering for types, and vertices are chosen u.a.r.\ amongst those with highest priority in the survival graph. The way in which each chosen vertex is coloured depends on the application: assuming that $x \geq y$, the vertices of type $xy$ and $yx$ are coloured red and blue if the algorithm aims for a small bisection, and the other way round if the algorithm aims for a large bisection. The survival graph is then updated just as in the case of large cuts. This prescription make sure that the two classes of the partition are balanced and that symmetric types have the same distribution in the survival graph if the input graph is random regular. The algorithm proceeds in this way until there is a small proportion of vertices remaining, which are then split evenly between the two classes in a final cleanup step, with negligible influence on the size of the cut.

To obtain the same bounds with the current approach, the algorithm needs to be slightly modified. This is because our definition of a local deletion algorithm does not allow vertices to be chosen in pairs. To address this issue, given $x>y$, we modify the prioritised algorithm so that instead of a symmetric pair of vertices of types $xy$ and $yx$, it chooses a vertex with type $xy$ or $yx$, with equal probabilities. Using this idea we can arrange  that the expected number of vertices of any type processed in a given chunk is essentially the same as if symmetric pairs were used, and this symmetry is preserved in the differential equations we derive. Theorem~\ref{chunkydeprio large girth} could then be used to show that   a.a.s.\   the numbers of vertices of types $xy$ and $yx$ are equal to the desired accuracy, and this inequality may be corrected in a final cleanup step, in which leftover vertices and vertices with clash colours are assigned colours arbitrarily, and then vertices may be swapped between red and blue to create the desired balance. 

As an alternative way to deal with the same issue, it would be very straightforward to extend our method to permit a finite number of vertices of specified  types, rather than just one, to be chosen in each step of a local deletion algorithm.

\subsection{Minimum connected and weakly-connected dominating sets}

Let $G=(V,E)$ be a graph and let $S \subseteq V$ be a dominating set. The set $S$ is a \emph{connected dominating set} if $G[S]$ is connected and a \emph{weakly-connected dominating set} if the subgraph of $G$ whose edge set consists of all edges incident with vertices in $S$ is connected. Duckworth and Mans~\cite{DM2}   determined upper bounds on the size of minimum connected and weakly-connected dominating sets in random $r$-regular graphs. Using the terminology of the current paper, this was done by analysing local deletion algorithms.

One of their algorithm for connected dominating sets, RAND ONE CDS, may be described as follows. In the first step, it selects a vertex $v$ u.a.r.\ and adds it to the output set, deleting it from the survival graph. In subsequent steps, the algorithm chooses a vertex $v$ u.a.r.\ among all vertices whose degree is less than $r$, then selects a neighbour $u$ of $v$ u.a.r. If $u$ has degree $r$, $v$ is added to the output set. If the degree of $u$ is less than $r$, the algorithm simply removes the edge $uv$ from the survival graph. This is easy to rephrase as a local deletion algorithm; note that the local rule is not normal.

 The algorithm for weakly-connected dominating sets is the same except for one modification:   if $u$ has degree $r$, then $u$, rather than $v$, is added to the output set.

These algorithms were analysed in~\cite{DM2} using the differential equation method as presented in~\cite[Theorem~6.1]{desurvey}.  Here, we may analyse their performance by considering the system of differential equations~\eqref{desystem} with $\tilde{p}_i(x,\y)=y_i(x)/(\sum_{k \neq r}y_k(x))$ and the functions $f_{j,i}$ naturally induced by the rules of the algorithms. We also include a short initial  ``burn-in" phase to ensure a good supply of vertices of degree less than $r$,  as in the proof of Theorem~\ref{t:cut}, without affecting the resulting bounds.  For fixed $\delta'>0$, Theorem~\ref{deprioranreg} implies that the random variable tracking the number of vertices of degree $r$ in the survival graph is governed by $ny_r(x)$ and that the size of the output set approaches $C(\delta')n$, where $C(\delta')=T/n$ with $T$ given by the theorem. It is easy to see that as $\delta'\rightarrow 0$ the value of $C(\delta')$ is achieved as $y_r(x)$ approaches zero, and hence the result coincides with the result obtained in~\cite{DM2}, as they both analyse the same process. We apply Theorem~\ref{chunkydeprio large girth} to conclude that the same holds for regular graphs with large girth.  (Note that we do not need to deprioritise these algorithms because there is no prioritisation between vertices of degree less than $r$, and the number of these vertices a.a.s.\ stays positive until the algorithm finishes.)
As for the dominating set example in Section~\ref{mids}, we rely on the numerical solutions of the differential equations to deduce that a negligible proportion of the  vertices remain uncovered at the end of the algorithm, and clash vertices are also sufficiently small in number. In this case, however, we need to ensure that the final set of vertices has the required connectivity properties. There is also an additional problem:  since the algorithm given by Theorem~\ref{chunky des} is chunky, instead of starting with a single connected (resp.\ weakly-connected) component, it creates many of them (in fact $O(\eps' n)$ components where $\eps'$ is given in that theorem). 
Because of this and the other vertices not yet dominated, to produce a connected (resp.\ weakly-connected) dominating set from the output produced by the chunky algorithm, we still need to greedily add vertices to the output set. However, any connected component can be merged with some other connected component with the addition of at most two vertices to the output set, and hence the bounds reported are not affected. We deduce that the   upper bounds shown in~\cite{DM2} to hold for random $r$-regular graphs using these two algorithms also hold for $r$-regular graphs with sufficiently large girth. (Here we again use the observation that the bounds were of necessity generous by some positive additive constant.)

There were further algorithms analysed in~\cite{DM2}, but it is a little harder to apply our present methods to them, and we refrain from doing so in the present paper.

\subsection{Maximum induced forest and minimum feedback vertex set}

 An induced forest in a graph $G=(V,E)$ is an acyclic induced subgraph of $G$. The size of the induced forest $G[S]$ is defined to be $|S|$.  Using a chunky local deletion algorithm, and analysing it in the large-girth setting, the authors~\cite{hoppen_wormald1} obtained lower bounds on the size of a maximum induced forest in $r$-regular graphs with large girth. 
The set of vertices {\em not} in such a forest is a minimum feedback vertex set, or decycling set, which is defined as a minimum set of vertices whose removal destroys all cycles. That work  followed~\cite{LW} quite closely, by directly considering graphs with large girth and proving independence lemmas as required. 

It would be straightforward to obtain the same results by analysing the algorithm in~\cite{hoppen_wormald1}, which is a local deletion algorithm with a normal subrule, via the machinery in the present paper. There are six colours available: three output colours, which we call purple, yellow and $\propto$, and three transient colours, which we call neutral, blue and $\propto'$. The intuition is that purple vertices represent the set $\U$ of vertices added to the forest, blue vertices have exactly one neighbour in the forest (and therefore are available to be added to $\U$ without inducing cycles in the input graph), yellow vertices are adjacent to at least two vertices in the forest (and therefore have been deleted to avoid the creation of cycles) and neutral vertices are not adjacent to any vertex in the forest. More precisely,  there is a first step where the algorithm chooses neutral vertices u.a.r., adds each such vertex $v$ to the set of purple vertices, colours its neutral  neighbours blue (and of course deletes $v$ from the survival graph). In the subsequent steps, the algorithm selects blue vertices u.a.r.\ and adds them to the set of purple vertices. Their neutral neighbours become blue, while their blue neighbours become yellow and hence are deleted from the survival graph. 

For the chunky local deletion algorithm associated with the local, selection and recolouring rules naturally derived from the above discussion, we may compute the functions $f_{k,i}$ defined in Lemma~\ref{l:trend} using calculations similar to those in Section~\ref{s:isets}, and then apply  Theorem~\ref{chunky des}. The resulting differential equations are equivalent to those in~\cite{hoppen_wormald1} by a straightforward substitution of variables, and hence the same results follow. A degree-governed version of this algorithm will be analysed using Theorem~\ref{chunkydeprio large girth} in~\cite{hoppen2}. (See also~\cite{forests_extended_abs} where this is described in a less economical direct approach.) 

\section{Final remarks}\lab{s:final}

The kinds of bounds we obtain for $r$-regular graphs with large girth sometimes carry over to graphs with large girth and \emph{maximum} degree $r$. This is true, for instance, if $P$ is \emph{monotone}, that is, if $H$ satisfies $P$ then every subgraph of $H$ also satisfies $P$ in $G$.  Fix a monotone property $P$ and suppose that there exists
$g>0$ such that every $r$-regular graph $H$ with girth at least
$g$ contains an induced subgraph $H_0$ satisfying $P$ such that  $|H_0| \geq
\gamma |V(H)|$ (and hence $f_P(H) \geq \gamma$). It is not hard to see that, given an $n$-vertex graph $G$ with maximum degree $r$ and girth at least $g$, we may construct a graph $G'$ by taking copies of $G$ and joining vertices in different copies so as to make $G'$ $r$-regular. This can be done without decreasing the girth if sufficiently many copies of $G$ are used. In particular, given a monotone property $P$, we have the inequality $f_P(G)/n \geq f_P(G')/|V(G')|$, since this holds for the copy of $G$ containing the largest number of vertices in a maximum induced subgraph with property $P$ in $G'$. Thus bounds for $r$-regular graphs with large girth imply bounds for graphs with large girth and maximum degree $r$. For instance, the property of having chromatic number at most $k$ is monotone, and  for $k=1$, this means that the graph is an independent set. Hence, we may deduce results about the size of a largest  independent set in a graph  in terms of its maximum degree, using the results of Sections~\ref{Isets thesis} and~\ref{sec:improved}. 

We note that our methods can readily be extended so as to apply to   $r$-regular bipartite graphs, or indeed multipartite graphs, of large girth.  Most of the technical results would go through with little effort.  In particular, the sharp concentration result via switchings in  the proof of Theorem~\ref{chunkyranreg}    is easily modified appropriately.
 
A general open problem is to determine how much the definition of local deletion algorithms can be relaxed, without losing the connection between random regular graphs and  regular graphs of large girth.

\end{document}